\documentclass[a4paper]{amsart}
\usepackage{amsmath}
\usepackage{amssymb}
\usepackage{mathtools}
\usepackage[a4paper,hmargin=26mm,vmargin=28mm]{geometry}
\usepackage{booktabs}
\usepackage{enumitem}
\setlist[enumerate]{label=\textup{\alph*)}}
\usepackage[svgnames,table]{xcolor}
\usepackage{colortbl}
\usepackage{expl3}

\ExplSyntaxOn

\clist_const:Nn \l_beta_numbers_clist {23,21,18,17,16,15,8,6,5,4,3,2}
\newcommand\betanumbers{(\clist_use:Nn \l_beta_numbers_clist {, } )}
\newcommand\halpha[3]{
  \int_eval:n { \clist_item:Nn \l_beta_numbers_clist{#1}
               - \clist_item:Nn \l_beta_numbers_clist{#2}-#3}
}
\clist_new:N \g_beta_h_clist
\int_new:N \l_i_int
\int_new:N \l_j_int
\int_new:N \l_sign_int
\cs_new:Nn \__compute_sign_of_beta_h_clist: {
  \int_set:Nn \l_sign_int {0}
  \int_set:Nn \l_i_int {1}
  \int_until_do:nNnn {\l_i_int} = {\clist_count:N \g_beta_h_clist}
  {
    \int_set_eq:NN \l_j_int \l_i_int
    \int_until_do:nNnn {\l_j_int} = {\clist_count:N \g_beta_h_clist}
    {
      \int_incr:N \l_j_int
      \int_compare:nT { \clist_item:Nn\g_beta_h_clist{\l_i_int} < \clist_item:Nn\g_beta_h_clist{\l_j_int} }
      {
        \int_incr:N \l_sign_int
      }
    }
    \int_incr:N \l_i_int
  }
  \int_if_even:nTF {\l_sign_int} {+}{-}
}

\tl_new:N \g_blam_row_tl
\tl_new:N \l_blam_col_tl
\int_new:N \l_beta_int
\newcommand\blamsth[3]
{
  \clist_gclear:N \g_beta_h_clist
  \int_zero:N \l_beta_int
  \tl_gclear:N \g_blam_row_tl
  \int_do_until:nn { \l_beta_int = 12 }
  {
    \int_incr:N \l_beta_int
    \int_case:nnF {\l_beta_int }
    {
      {#1}
      {
        \tl_gput_right:Nn \g_blam_row_tl {\color{red}}
        \tl_set:Nx \l_blam_col_tl {\int_eval:n{\clist_item:Nn\l_beta_numbers_clist{#1}-#3}}
      }
      {#2}
      {
        \tl_gput_right:Nn \g_blam_row_tl {\color{blue}}
        \tl_set:Nx \l_blam_col_tl {\int_eval:n{\clist_item:Nn\l_beta_numbers_clist{#2}+#3}}
      }
    }
    {
      \tl_set:Nx \l_blam_col_tl {\clist_item:Nn\l_beta_numbers_clist {\int_use:N \l_beta_int}}
    }
    \tl_put_right:Nx \g_blam_row_tl { \tl_use:N \l_blam_col_tl }
    \tl_gput_right:Nn \g_blam_row_tl {&}
    \clist_gput_right:NV \g_beta_h_clist \l_blam_col_tl
  }
  &(&\tl_use:N\g_blam_row_tl)& #1 & #2 & #3 & \halpha{#1}{#2}{#3}
  & \__compute_sign_of_beta_h_clist:
}

\ExplSyntaxOff

\usepackage{tikz}
\usetikzlibrary{matrix,arrows.meta,decorations.markings}

\tikzset{
  wall/.style={black!80, thick},
  innerwall/.style={rectangle, draw=gray, thin},
  pics/diagram/.style = {  
    code = {
    \begin{scope}[scale=0.5]
      \xdef\lastMu{0}
      \xdef\lastR{0}
      \foreach \mu [count=\r] in {#1} {
        \foreach \col in {1,...,\mu} {
          \coordinate(-\r-\col) at (\col, -\r);
          \draw[innerwall](\col, -\r) rectangle ++(-1,1);
        }
        \draw[wall](0,-\r)--++(0,1) (\mu, -\r)--++(0,1);
        \ifnum\lastMu>\mu
           \draw[wall](\mu, 1-\r)--(\lastMu, 1-\r);
        \else\ifnum\lastMu=0
           \draw[wall](0,0)--++(\mu, 0);
          \fi
        \fi
        \xdef\lastMu{\mu}
        \xdef\lastR{\r}
      }
      \draw(0,-\lastR) -- + (\lastMu,0);
    \end{scope}
    }
  }
}

\newcounter{RowNumber}
\renewcommand\theRowNumber{\alph{RowNumber})}
\newcommand\RowNumber{%
  \makebox[2em][l]{%
    \ifnum\arabic{RowNumber}>0\theRowNumber\fi%
    \stepcounter{RowNumber}%
  }%
}
\newenvironment{NumberedArray}[2][2]
{%
  \begin{equation*}
    \rowcolors{#1}{white}{yellow!20!brown!5!white}
    \setcounter{RowNumber}{0}
    \array{@{\RowNumber} #2}\toprule
}{\bottomrule\endarray\end{equation*}}

\usepackage{array,collcell}
\newcolumntype{T}{>{\collectcell\textrm}{r}<{\endcollectcell}}

\title{Positive Jantzen sum formulas for cyclotomic Hecke algebras}
\author{Andrew Mathas}
\address{School of Mathematics and Statistics F07, University of Sydney, NSW 2006, Australia.}
\email{andrew.mathas@sydney.edu.au}

\subjclass[2000]{20G43, 20C08, 20C30, 05E10}
\keywords{Cyclotomic Hecke algebras, graded cellular algebras,
          Khovanov--Lauda--Rouquier algebras}

\usepackage{cite}
\usepackage[pagebackref=true]{hyperref}
\usepackage[hyperpageref]{backref}
\renewcommand*{\backref}[1]{}
\renewcommand*{\backrefalt}[4]%
  {{[\tiny\ifcase #1 Not cited.\relax\or Page~#2.\else Pages #2.\fi]}}

\hypersetup{
  pdftitle={Positive Jantzen sum formulas for cyclotomic Hecke algebras},%
  pdfauthor={Andrew Mathas},%
  colorlinks = true,%
  linkcolor =Green,%
  anchorcolor = red,%
  citecolor = Green,%
  urlcolor = Green%
}

\usepackage{xparse,etoolbox}
\DeclarePairedDelimiterX{\set}[1]{\{}{\}}{\setargs{#1}}
\NewDocumentCommand{\setargs}{>{\SplitArgument{1}{|}}m}{\setargsaux#1}
\NewDocumentCommand{\setargsaux}{mm}
{\IfNoValueTF{#2}{#1} {#1\,\delimsize|\,\mathopen{}#2}}

\newtheorem*{MainTheorem}{Main Theorem}

\usepackage{aliascnt}
\def\NewTheorem#1{%
  \newaliascnt{#1}{equation}%
  \newtheorem{#1}[#1]{#1}%
  \aliascntresetthe{#1}%
  \expandafter\def\csname #1autorefname\endcsname{#1}%
}
\def\equationautorefname~#1\null{(#1)\null}

\numberwithin{equation}{section}
\swapnumbers

\NewTheorem{Proposition}
\NewTheorem{Theorem}
\NewTheorem{Corollary}
\NewTheorem{Lemma}
\NewTheorem{Definition}
\theoremstyle{definition}
\NewTheorem{Example}
\AtEndEnvironment{Example}{\null\hfill$\Diamond$}%
\theoremstyle{remark}
\NewTheorem{Remark}

\newlist{Cases}{enumerate}{1}
\setlist[Cases]{label=Case~\alph*., font=\bfseries, wide=0pt,
                leftmargin=*, align=left, itemindent=!
                }

\let\<=\langle
\let\>=\rangle
\def\({\big(}
\def\){\big)}
\def\Prod{\displaystyle\prod}
\def\Sum{\displaystyle\sum}

\def\pmod#1{\text{ }(\textrm{mod } #1)\,}

\def\map#1#2{\,{:}\,#1\!\longrightarrow\!#2}

\newcommand{\C}{\mathbb{C}}
\newcommand{\N}{\mathbb{N}}
\newcommand{\Q}{\mathbb{Q}}
\newcommand{\Sym}{\mathfrak S}
\newcommand{\F}{\mathbb{F}}
\newcommand{\Fp}{{\F_p}}
\newcommand{\R}{\mathbb{R}}
\newcommand{\Z}{\mathbb{Z}}

\newcommand\Acal{\mathcal{A}}
\newcommand\Kcal{\mathcal{K}}
\newcommand\Ocal{\mathcal{O}}
\newcommand\Pcal{\mathcal{P}}
\newcommand\Zcal{\mathcal{Z}}

\newcommand\AIq{\Acal[I^n]}

\newcommand{\Klesh}[1][n]{\Kcal^\Lambda_{#1}}
\newcommand{\Kleshf}[1][n]{\Kcal^{\Lambda_f}_{#1}}
\newcommand{\Parts}[1][n]{\Pcal^\Lambda_{#1}}

\let\gedom=\trianglerighteq
\let\gdom=\vartriangleright

\newcommand\charge{{\boldsymbol\kappa}}

\NewDocumentCommand\blamst{ D(){h} }{\bbeta^\blam_{st}(#1)}
\newcommand\dkst{\delta_{l_sl_t}(n)}

\newcommand\bbeta{{\boldsymbol\beta}}
\newcommand\blam{{\boldsymbol\lambda}}
\newcommand\bmu{{\boldsymbol\mu}}
\newcommand\bnu{{\boldsymbol\nu}}

\def\s{{\mathfrak s}}
\def\t{{\mathfrak t}}
\def\u{{\mathfrak u}}
\def\v{{\mathfrak v}}

\newcommand\tlam{\t^\blam}

\newcount\tableauRow\newcount\tableauCol
\NewDocumentCommand\Tableau{ O{-2} D<>{0.4} m}{%
  \begin{tikzpicture}[baseline=#1mm, xscale=#2, yscale=0.4,
                      draw/.append style={thick,black},]
    \tableauRow=0
    \foreach \Row in {#3} {
    \tableauCol=1
    \foreach\k in \Row {
      \draw(\the\tableauCol,\the\tableauRow)+(-.5,-.5)rectangle++(.5,.5);
      \draw(\the\tableauCol,\the\tableauRow)node{$\k$};
      \global\advance\tableauCol by 1
    }
    \global\advance\tableauRow by -1
    }
  \end{tikzpicture}%
}

\newcommand\Bigger[2][5]{\left#2\rule{0mm}{#1truemm}\right.}

\def\Tritab(#1|#2|#3){\Bigger(%
             \Tableau[-3]{#1}\,\Bigger|\,%
             \Tableau[-3]{#2}\,\Bigger|\,%
             \Tableau[-3]{#3}%
             \Bigger)%
}

\newcommand\bi{\mathbf{i}}
\newcommand\bj{\mathbf{j}}
\newcommand\bk{\mathbf{k}}

\newcommand{\DeclareMyOperator}[1]{%
  \expandafter\DeclareMathOperator\csname #1\endcsname{#1}%
}
\forcsvlist{\DeclareMyOperator}{
  Add,End,Hom,Ind,pdeg, rad,Rem,Rep,soc,Std,Shape,diag%
}

\DeclareMathOperator\Char{char}
\newcommand\con{\textsf{c}}
\newcommand\res{\textsf{r}}

\DeclareMathOperator\noarrow{\:\rlap{\hspace*{0.25em}/}\text{---}\:}

\newcommand\Chq{\mathop{\rm ch}\nolimits_q}
\newcommand\ch{\mathop{\rm ch}\nolimits}

\newcommand\grMod{\text{-}\mathop{\rm grMod}}

\usepackage{mathrsfs}

\newcommand\Hkn[1][v]{\mathscr{H}^\charge_{#1}}
\newcommand\HK[1][n]{\mathscr{H}^\Lambda_{#1}(\Kcal)}
\newcommand\HO[1][n]{\mathscr{H}^\Lambda_{#1}(\Ocal)}
\newcommand\Hn[1][n]{\mathscr{H}^\Lambda_{#1}}
\newcommand\Hnf[1][n]{\mathscr{H}^{\Lambda_f}_{#1}}

\newcommand\Ln[1][n]{\mathscr{L}_{#1}}

\newcommand\Rn[1][n]{\mathscr{R}^\Lambda_{#1}}
\newcommand\Rnf[1][n]{\mathscr{R}^{\Lambda_f}_{#1}}
\newcommand\Rbeta[1][\beta]{\Rn[#1]}
\newcommand\Rlam[1][\blam]{\mathscr{R}^{\gdom#1}_n}

\newcommand\Df{\mathsf{A}_{f,e}}
\newcommand\Lf{\mathsf{L}_{f,e}}

\newcommand\Hlam[1][\blam]{\mathscr{H}^{\gdom#1}_n}

\newcommand\Slam[1][\blam]{\mathbb{S}^{#1}}
\newcommand\SOlam[1][\blam]{S^{#1}_\Ocal}
\newcommand\Dmu[1][\bmu]{\mathbb{D}^{#1}}
\newcommand\Emu[1][\bmu]{\mathbb{E}^{#1}}

\begin{document}
\begin{abstract}
  We prove a ``positive'' Jantzen sum formula for the Specht modules of
  the cyclotomic Hecke algebras of type~$A$. That is, in the
  Grothendieck group, we show that the sum of the  pieces of the Jantzen
  filtration is equal to an explicit non-negative linear combination of
  modules $E^\bnu_{f,e}$, which are modular reductions of simple modules
  for closely connected Hecke algebras in characteristic zero. The
  coefficient of $E^\bnu_{f,e}$ in the sum formula is determined
  by the graded decomposition numbers in characteristic zero, which are
  known, and the characteristic of the field. As a consequence we see
  that the decomposition numbers of a cyclotomic Hecke algebra at an
  $e$th root of unity in characteristic $p$ depend on the decomposition
  numbers of related cyclotomic Hecke algebras at $ep^r$th roots of
  unity in characteristic zero, for $r\ge0$.
\end{abstract}
\maketitle

\begin{quote}
  \centerline{\textit{\large Dedicated to Professor Gordon Douglass James, 31/12/1945--5/12/2020}}

  This is the first paper that I have written about Jantzen sum formulas
  without my close collaborator and mentor Professor Gordon James.
  Gordon was a powerful and insightful mathematician who was a key
  player in shaping our current knowledge of the representation theory
  of the symmetric groups, Hecke algebras and Schur algebras. I feel
  privileged to have worked with him and I am proud to count him among
  my friends and colleagues. Even though Gordon did not contribute
  directly to this paper it is infused with his influence because it
  answers questions that Gordon and I discussed when we were working
  together.
\end{quote}



\section{Introduction}

  Jantzen filtrations have played a dominant role in representation
  theory since they were introduced by Jantzen in
  1977~\cite{Jantzen:thesis}. In particular, stronger forms of the
  Lusztig and Kazhdan-Lusztig conjectures
  \cite{KL,Lusztig:QuantumGroupConjectures} say that the Jantzen
  filtrations and weight filtrations of Weyl modules coincide, results
  that are equivalent to (questions or) conjectures of
  Jantzen~\cite{Jantzen:HightestWeightModules}. For quantum
  groups, complex semisimple Lie algebras,  and
  Soergel bimodules these conjectures are now known to be
  true~\cite{BeilinsonBernstein:JantzenConjectures,
             BrylinskiKashiwara:KLConjecture,
             GaberJoseph:JantzenConjectures,
             KazhdanLusztig:KLQuantumGrroupConjecture,
             Soergel:AndersenFiltration,
             Williamson:LocalHodge}.

  Jantzen's observation for defining the Jantzen filtrations is
  extremely simple and elegant. Let $p$ be a rational prime and $M_\Z$ a
  finitely generated free $\Z$-module that comes equipped with a
  non-degenerate bilinear form $\<\ ,\ \>_\Z$. Then $M_\Z$ has a filtration
  $M_\Z=J_0(M_\Z)\supseteq J_1(M_\Z)\supseteq J_2(M_\Z)\supseteq\dots$
  defined by
  \[
      J_k(M_\Z)=\set{m\in M_\Z|\<m,x\>\in p^k\Z\text{ for all }x\in M_\Z}.
  \]
  Reducing modulo~$p$, let $M_\Fp=M\otimes_\Z\Fp$. The
  \textbf{Jantzen filtration}
  $M_\Fp=J_0(M_\Fp)\supseteq J_1(M_\Fp)\supseteq J_2(M_\Fp)\supseteq\dots$
  of $M_\Fp$ is defined by:
  \[
        J_k(M_\Fp)=\bigl(J_k(M_\Z)+pM_\Z\bigr)/pM_\Z,\qquad\text{ for }k\ge0.
  \]
  Let $G_M$ be the Gram matrix of $\<\ ,\ \>_\Z$ with respect to some basis
  of~$M_\Z$. By looking at the Smith normal form of $G_M$, Jantzen observed
  the trivial but important fact that
  \[
            \nu_p(\deg G_m) = \sum_{k>0}\dim_{\Fp}J_k(M_\Fp),
  \]
  where $\nu_p$ is the usual $p$-adic valuation. In category
  $\mathscr{O}$, a Weyl module $\Delta^\lambda$ is uniquely determined
  by the dimensions of its weight spaces. Applying this machinery to the
  weight spaces $\Delta^\Lambda_\mu$ of $\Delta^\lambda$ determines the
  dimensions of the $\mu$-weight spaces in the Jantzen filtration of a
  Weyl module $\Delta^\lambda$. This leads to the \textbf{Jantzen sum
  formula} of $\Delta^\lambda$, which describes the sum
  $\sum_{k>0} J_k(\Delta^\lambda)$ in the Grothendieck group as an
  explicit $\Z$-linear combination of more dominant Weyl
  module~\cite[II.8.19]{Jantzen:book}.

  One limitation of all of the Jantzen sum formulas in the literature is
  that they express $\sum_{k>0}J_k(M)$ in the Grothendieck group as a
  $\Z$-linear combination of more dominant Weyl modules or Specht
  modules, where many of the coefficients are negative. This is
  ``wrong'' in the sense that $\sum_{k>0}J_k(M)$ can definitely be
  written as a \textit{non-negative} linear combination of the images of
  the simple modules in the Grothendieck group. On the other hand, these
  negative coefficients in the sum formulas have to appear because, in
  general, the simple modules cannot be written as non-negative linear
  combinations of the Weyl modules in the Grothendieck group. It is
  natural to ask for a new version of the Jantzen sum formula that
  expresses $\sum_{k>0}J_k(M)$ as a non-negative linear combination of
  simple modules.

  Soon after Jantzen proved his sum formula, his student Schaper used
  Jantzen's sum formula for the general linear groups to prove a sum
  formula for the Specht modules of the symmetric groups~\cite{Schaper}.
  Analogues of this result have since proved for the Specht modules of
  the Ariki-Koike
  algebras~\cite{JamesMathas:JantzenSchaper,JamesMathas:cycJantzenSum},
  which are cyclotomic Hecke algebras of type~$A$ that include the
  symmetric group and their Iwahori-Hecke algebras as special cases. In
  the literature, rather than working in the module category of the
  Hecke algebras, this result is always deduced by first proving a sum
  formula the Weyl modules of a related Schur algebra and then applying
  a Schur functor to deduce the sum formula for the Specht modules. In
  particular, the literature does not contain the proof of a sum formula
  for these algebras that takes place entirely inside the module
  category of Ariki-Koike algebras or symmetric groups.

  This paper proves a new \textit{positive} Jantzen sum formulas for the
  Specht modules of the cyclotomic Hecke algebras of type~$A$ that
  writes $\sum_{k>0}J_k(S^\blam)$ as an explicit non-negative linear
  combination of certain modules $E^\bmu_{f,e}$. This sum
  formula is proved entirely within the module categories of the Hecke
  algebras.  The modules $E^\bmu_{f,e}$ that appear in our sum formulas
  are not, in general, simple and instead these modules are modular reductions of
  the simple modules of related cyclotomic Hecke algebras in
  characteristic zero. The coefficients that appear in the positive sum
  formula are described explicitly in terms of the (derivatives of
  graded) decomposition numbers.

  To state our main result we quickly introduce the notation that we need
  and refer the reader to later sections for the full definitions.  Let
  $F$ be a field of characteristic $p\ge0$ and let $\Hn(F)$ be a
  cyclotomic Hecke algebra over~$F$ with a Hecke parameter $\xi$ of
  quantum characteristic~$e$ (\autoref{D:CyclotomicHecke}). Let $\Parts$
  be the set of $\ell$-partitions. For $\blam\in\Parts$ let $S^\blam_F$
  be the corresponding Specht module for $\Hn(F)$. If $M$ is an
  $\Hn(F)$-module let $[M]$ be its image in the Grothendieck group
  $\Rep\Hn(F)$ of $\Hn(F)$.

  For $f=ep^r$, for $r\ge0$, let $\zeta_f\in\C$ be a primitive $f$th
  root of unity and let $\set{E^\bnu_{\C,f}|\bnu\in\Kleshf}$ be the set
  of simple modules for a corresponding Hecke algebra $\Hnf(\C)$
  over~$\C$ that has Hecke parameter~$\zeta_f$; see \autoref{E:Lambdaf}.
  By \autoref{P:fAdjustment}, there is a map of Grothendieck groups
  $\Df\map{\Rep\Hnf(\C)}\Rep\Hn(F)$ when $p>0$. Let
  $E^\bnu_{f,e}=\Df\bigl(E^\bnu_{\C,f}\bigr)$ and set $E^\bnu_{f,e}=0$
  if $p=0$.  Finally, let $d^{\C,f}_{\blam\bmu}(q)$ be the graded
  decomposition numbers for $\Hnf(\C)$, for $\blam\in\Parts$ and
  $\bmu\in\Kleshf$; see \autoref{E:GradedDecomp}.  Then
  $d^{\C,f}_{\blam\bmu}(q)$ is a parabolic Kazhdan-Lusztig polynomial
  that, in principle, is known. Let $(d^{\C,f}_{\blam\bmu})'(1)$ be the
  derivative of $d^{\C,f}_{\blam\bmu}(q)$ evaluated at $q=1$. Our main
  result is then the following:

  \begin{MainTheorem}\hypertarget{T:MainTheorem}{}
    Suppose that $F$ is a field of characteristic $p\ge0$ and let $\blam\in\Parts$.
    In $\Rep\Hn(F)$,
    \[
    \sum_{k>0}[J_k(S^\blam_F)] = \sum_{\bmu\in\Klesh}
        (d^{\C,e}_{\blam\bmu})'(1)[E^\bmu_{F,e}]
        + \sum_{r>0}(p-1)p^{r-1}\sum_{\substack{f=ep^r\\\bnu\in\Kleshf}}
               (d^{\C,f}_{\blam\bnu})'(1)[E^\bnu_{f,e}].
    \]
  \end{MainTheorem}

  Strikingly, in characteristic $p$ the right hand side of this sum
  formula depends on modules for Hecke algebras at $ep^r$th roots of
  unity. Moreover, by \autoref{L:positivity} below, all of the
  coefficients on the right hand side of the sum formula are
  non-negative integers, so this really is a positive Jantzen sum
  formula.

  When $F$ is a field of characteristic zero the second sum over $r>0$
  on the right hand side of our main theorem vanishes vanishes.  In this
  case, the new sum formula can be proved by assuming the deep fact that
  the Jantzen filtrations of graded Specht modules coincide with their
  grading filtrations. In fact, we do not prove our main result this way
  and, indeed, we cannot prove it this way because the Jantzen
  filtrations are only known to coincide with the grading filtrations
  for cyclotomic Hecke algebras of level~$1$, where this follows from
  work of Shan~\cite[Theorem~0.1]{Shan:JantzenFiltrations}.  (Shan
  actually proves that the Jantzen filtrations of the Weyl modules of
  the $q$-Schur algebras coincide with their grading filtrations. Using
  the Schur functor, it is possible to prove the corresponding result for
  the graded Specht modules.)

  The proof of our main theorem takes place entirely inside the module
  categories of the cyclotomic Hecke algebras. The main idea is to work
  with the formal characters of $\Hn(F)$-modules
  (\autoref{D:character}). The key point is that character map is
  injective so the images of all modules in the Grothendieck group,
  including the images of the Jantzen filtration, are determined by
  their characters. This is essentially the same idea sketched above for
  computing the Jantzen sum formula of Weyl modules

  Using seminormal forms it is relatively easy to determine the
  ``Jantzen characters'' (\autoref{P:JantzenDetCharacter}). The proof of
  our main theorem follows by applying some combinatorial tricks using
  the formal characters of the graded Specht modules of the cyclotomic
  KLR algebras type~$A$, which are isomorphic to the cyclotomic Hecke
  algebras.  It is quite remarkable that the graded Specht modules enter
  this story, but their appearance in our arguments is almost accidental
  because instead of marking a direct connection between the Jantzen
  filtrations and the graded representation theory we only exploit some
  combinatorial shadows that are common to both settings.

  Here is a brief outline of the paper. Section~2 reviews the results
  that we need from the graded representation theory of the cyclotomic
  KLR algebras of type~$A$. In particular, we introduce the
  $q$-character map for these algebras and use them to give a new proof
  that the graded decomposition matrices in positive characteristic
  factor through those in characteristic zero. Section~3 introduces the
  cyclotomic Hecke algebras and their character map. We define
  seminormal bases and use these to give a new factorisation of the Gram
  determinants that is compatible with taking characters.  Section~4
  contains the heart of the paper. We properly introduce the Jantzen
  filtrations and the Jantzen characters and then prove our main theorem
  by combining the results from the previous sections.  Finally,
  section~5 shows how to prove the ``classical'' Jantzen sum formula
  using the framework of this paper. As the classical Jantzen sum
  formula is more explicit it takes more effort to prove than our
  positive sum formula, but it is worth the effort because the new
  ``classical'' sum formula that we obtain is nicer than the existing
  sum formula in the literature~\cite{JamesMathas:cycJantzenSum}.

  The definition of cyclotomic Hecke algebras used in this paper
  includes the \textit{degenerate cyclotomic algebras}
  of~\cite{Brundan:degenCentre} as the special case when the Hecke a
  parameter is one. Jantzen filtrations for these algebras have not
  appeared in the literature previously, so even the ``classical''
  Jantzen sum formula that we obtain in the generate case is completely
  new.

\subsection*{Acknowledgements}
We thank Jun Hu for useful discussions about the results of this paper.
This research was supported, in part, by the Australian Research
Council.

\section{Cyclotomic KLR algebras and their modules}

  In this section we introduce the cyclotomic KLR algebras of type~$A$,
  which are isomorphic to the cyclotomic Hecke algebras of type~$A$.  We
  will use these algebras to define the modules $E^\bmu_{f,e}$
  from the introduction and to prove some combinatorial identities that
  we use in \autoref{SS:PositiveSum} to compute characters for
  the Jantzen filtration.

  \subsection{Graded modules and algebras}
  In this paper all modules are finitely generated modules over a
  commutative ring with one. A \textbf{graded module} $M$ is a module
  with a $\Z$-grading $M=\bigoplus_{d\in\Z}M_d$.  Similarly,
  \textbf{graded algebra} $A$ is a $\Z$-graded algebra
  $A=\bigoplus_{d\in\Z}A_d$ with $A_dA_e\subseteq A_{d+e}$, for all
  $d,e\in\Z$. If $A$ is a graded algebra then a graded $A$-module is a
  graded (right) $A$-module such that $M_dA_e\subseteq M_{d+e}$, for
  $d,e\in\Z$. If $F$ is a field and $M$ is a graded $F$-module then the
  \textbf{graded dimension} of $M$ is \[\dim_q M = \sum_{d\in\Z}(\dim
  M_d)q^d\quad\in\N[q,q^{-1}],\] where $q$ is an indeterminate
  over~$\Z$. Given an integer $s$ let $q^sM$ be the graded module
  obtained from~$M$ by \textbf{shifting} the grading of~$M$ by
  $s$ so that $(q^sM)_d = M_{d-s}$. More generally, given a Laurent
  polynomial $s(q)=\sum_d s_dq^d\in\N[q,q^{-1}]$ let \[     s(q)M =
  \bigoplus_{d\in\Z} (q^dM)^{\oplus s_d}.\] Then $\dim_q s(q)M =
  s(q)\dim_qM$. Let $A\grMod$ be the category of finitely generated
  graded (right) $A$-modules with degree preserving homomorphisms.

  If $A$ is a graded algebra let $\underline{A}$ be the (ungraded)
  algebra obtained by forgetting the grading on $A$. Similarly, if $M$
  is a graded $A$-module let $\underline{M}$ be the (ungraded)
  $\underline{A}$-module obtained by forgetting the grading on~$M$. This
  defines an exact functor from the category of finite dimensional
  graded $A$-modules to the category of finite dimensional
  $\underline{A}$-modules.

  \subsection{Graded cellar algebras}\label{S:Cellular}
   Cellular algebras, which were introduced by Graham and
   Lehrer~\cite{GL}, provide a convenient framework for
   constructing the simple modules of an algebra.

   Let $R$ be a commutative ring with one.

   \begin{Definition}[{Graham and Lehrer~\cite[Definition~1.1]{GL},
     Hu-Mathas~\cite[Definition~2.1]{HuMathas:GradedCellular}}]
     Let $A$ be a $\Z$-graded $R$-algebra that is free and of finite rank. A
     \textbf{graded cell datum} for $A$ is an ordered quadruple
     $(\Pcal, T, C, \deg)$ where $(\Pcal, \ge)$ is a finite poset, $T(\lambda)$
     is a finite set for $\lambda\in\Pcal$,
     \[
        C\map{\coprod_{\lambda\in\Pcal} T(\lambda)\times T(\lambda)} A;
        (\s,\t)\mapsto c_{\s\t}
        \qquad\text{and}\qquad
        \deg\map{\coprod_{\lambda\in\Pcal} T(\lambda)}\Z; \t\mapsto\deg\t;
     \]
     are maps such that $C$ is injective and the following hold:
     \begin{enumerate}[label=C${}_{\arabic*}$\upshape), start=0]
       \item For $\s,\t\in T(\lambda)$ and $\lambda\in\Pcal$, the
         element $c_{\s\t}\in A$ is homogeneous of degree $\deg\s+\deg\t$.
       \item The set $\set{c_{\s\t}|\s,\t\in T(\lambda), \lambda\in\Pcal}$ is
         a basis of $A$.
       \item If $a\in A$, $\s,\t\in T(\lambda)$ and $\lambda\in\Pcal$
       then there exist scalars $r_{\t\v}(a)\in R$ that do not depend
       on~$\s$, such that
       \[
            c_{\s\t}a \equiv \sum_{\v\in T(\lambda)} r_{\t\v}(a) c_{\s\v}
            \pmod{A^{>\lambda}},
       \]
       where $A^{>\lambda}$ is the $R$-submodule of $A$ spanned by
       $\set{c_{\u\v}|\u,\v\in T(\mu)\text{ for }\mu>\lambda}$.
      \item There is an unique algebra anti-isomorphism $*\map AA$ such
      that $c_{\s\t}^* = c_{\t\s}$, for all $\s,\t\in T(\lambda)$ and
      $\lambda\in\Pcal$.
     \end{enumerate}
     A \textbf{graded cellular algebra} is an algebra that has a graded
     cell datum.

     A \text{cellular algebra} is a graded cellular algebra
     that is concentrated in degree zero. That is, it has a graded cell
     datum $(\Pcal, T, C, \deg)$ with $\deg\s=0$. for all $\s$. In this
     case, $(\Pcal, T, C)$ is a \textbf{cell datum} for~$A$.
   \end{Definition}

   In particular, note that if $A$ is a graded cellular algebra then,
   forgetting the grading, $\underline{A}$ is a cellular algebra.

   Cellular algebras were defined by Graham and Lehrer~\cite{GL} with
   the natural extension of their definitions and results to the graded
   setting being given in~\cite{HuMathas:GradedCellular}. The proofs of
   all of the results in this section can be found in these two papers.

   Let $A$ be a (graded) cellular algebra with cell datum $(\Pcal,T, C,
   \deg)$. Using the definitions it is easy to see that for each $\lambda\in\Pcal$
   there exist a \textbf{(right) cell module} $C^\lambda$ with basis
   $\set{c_\t|\t\in T(\lambda)}$ and with $A$-action
   \[
     c_\t a = \sum_{\v\in T(\lambda)} r_{\t\v}(a) c_\v,\qquad\text{for }a\in A,
   \]
   where the scalar $r_{\t\v}(a)\in R$ is from (C$_2$). Similarly, let
   $C^{\lambda*}$ be the \textbf{left cell module} indexed by $\lambda$,
   which is isomorphic to $C^\lambda$ as a vector space and where the
   $A$-action is given by
   \[
     ac_\t = \sum_{\v\in T(\lambda)} r_{\t\v}(a^*) c_\v,\qquad\text{for }a\in A,
   \]
   Extending the notation from C$_2$, let
   $A^{\ge\lambda}=\<c_{\u\v}\mid\u,\v\in T(\mu)\text{ for }\mu\ge\lambda\>_R$ and
   define the $(A,A)$-bimodule map $c_\lambda$ by
   \[
       c_\lambda\map{C^{\lambda*}\otimes C^\lambda} A^{\ge\lambda}/A^{>\lambda};
            c_\s\otimes c_\t\mapsto c_{\s\t}+A^{>\lambda},
            \qquad\text{for }\s,\t\in T(\lambda).
   \]
   Using (C$_2$) and (C$_3$), it follows that
   $C^\lambda$ has a homogeneous symmetric bilinear form
   $\<\ ,\ \>\map{C^\lambda\times C^\lambda}R$ of degree zero such that
   \begin{equation}\label{E:InnerProduct}
            x\cdot c_\lambda(y\otimes z) = \<x,y\>z,\qquad\text{for }x,y,z\in C^\lambda.
   \end{equation}
   In particular, $c_\t c_{\u\v}=\<c_\t,c_\u\>c_\v$, for $\t,\u,\v\in T(\lambda)$.
   This form is associative in the sense that
   $\<xa,y\>=\<x,ya^*\>$, for all $x,y\in C^\lambda$ and $a\in A$. Hence,
   \[
        \rad C^\lambda = \set{x\in C^\lambda|\<x,y\>=0\text{ for all }y\in C^\lambda}
   \]
  is a (graded) $A$-submodule of $C^\lambda$. Set
  $D^\lambda=C^\lambda/\rad C^\lambda$. Then $D^\lambda$ is a (graded)
  $A$-module.

  Let $D$ be an irreducible graded $A$-module. Then $q^k D$ is a
  non-isomorphic irreducible graded $A$-module, for $k\in
  \Z\setminus\set{0}$. If $M$ is a graded $A$-module then the
  \text{graded decomposition multiplicity} of~$D$ in~$M$ is the Laurent
  polynomial
  \begin{equation}
    [M:D]_q = \sum_{k\in\Z} [M:q^kD] q^k\quad\in\quad\N[q,q^{-1}],
  \end{equation}
  where $[M:q^kD]$ is the multiplicity of $q^kD$ as a composition factor of~$M$.

  The main results from the theory of (graded) cellular algebras are:

  \begin{Theorem}[Graham and Lehrer~\cite{GL,HuMathas:GradedCellular}]
    Let $F$ be a field and suppose that $A$ is a graded cellular
    algebra. Then:
    \begin{enumerate}
      \item Let $\mu\in\Pcal$ and suppose that $D^\mu\ne0$. Then
      $D^\mu$ is a self-dual graded irreducible $A$-module
      \item Let $\Pcal_0=\set{\mu\in\Pcal|D^\mu\ne0}$. Then
      $\set{q^kD^\mu|\mu\in\Pcal_0\text{ and }k\in\Z}$
      is a complete set of pairwise non-isomorphic irreducible
      graded $A$-modules.
      \item If $\lambda\in\Pcal$ and $\mu\in\Pcal_0$. Then
      $[S^\mu:D^\mu]_q=1$ and
      $[S^\lambda:D^\mu]_q\ne0$ only if $\lambda\ge \mu$.
    \end{enumerate}
  \end{Theorem}

  We need the theory of graded cellular algebras only to explain the
  definition of the Specht modules of the cyclotomic Hecke algebras. The
  bilinear form on the Specht modules is the key to defining their
  Jantzen filtrations. Our main results also use the graded Specht
  modules, although we have to approach this indirectly because the
  bilinear form on the graded Specht modules is degenerate unless the
  graded Specht module is irreducible.

  \subsection{Cyclotomic KLR algebras}\label{SS:CyclotomicKLR}
  Fix an integer $e\in\set{2,3,4,\dots}\cup\set{\infty}$ and let
  $\Gamma=\Gamma_e$ be the quiver with vertex set $I=\Z/e\Z$ (we set
  $e\Z=\set{0}$ when $e=\infty$), and with edges $i\to i+1$, for $i\in I$.
  To the quiver $\Gamma$ we attach the standard Lie theoretic data
  of a Cartan matrix $(c_{ij})_{i,j\in I}$, fundamental weights
  $\set{\Lambda_i|i\in I}$, the positive weight lattice
  $P_e^+=\sum_{i\in I}\N\Lambda_i$, the positive root lattice
  $Q_e^+=\bigoplus_{i\in I}\N\alpha_i$ and
  we let $(\cdot,\cdot)$ be normalised invariant form determined by
  \[(\alpha_i,\alpha_j)=c_{ij}\qquad\text{and}\qquad
          (\Lambda_i,\alpha_j)=\delta_{ij},\qquad\text{for }i,j\in I.\]
  Let $Q^+_{e,n}=\set{\beta\in Q_e^+|\sum_{i\in I}(\Lambda_i,\beta)=n}$ and
  for $\beta\in Q^+_{e,n}$ let
  $I^\beta=\set{\bi\in I^n|\beta=\alpha_{i_1}+\dots+\alpha_{i_n}}$.

  Following Rouquier~\cite{Rouquier:QuiverHecke2Lie}, let $u$ and $v$ be indeterminates
  and for $i,j\in I$ define polynomials
  \[
       Q_{ij}(u,v) = \begin{cases}
                         (u- v)(v-u) & \text{if } i \leftrightarrows j,\\
                         u - v & \text{if } i \to j, \\
                         v -u & \text{if } i \leftarrow j, \\
                         0 & \text{if } i=j, \\
                         1 & \text{if } i \noarrow j.
       \end{cases}
  \]

  To define the cyclotomic KLR algebras of type $A$, fix a non-negative
  integer $n\ge0$ and a dominant weight $\Lambda\in P_e^+$ and let
  $\ell=\sum_{i\in I}(\Lambda,\alpha_i)$.

  \begin{Definition}[\protect{%
    Khovanov and Lauda~\cite{KhovLaud:diagI,KhovLaud:diagII}
    and Rouquier~\cite{Rouquier:QuiverHecke2Lie}}]\label{D:KLRAlgebra}
    Let $e\in\set{2,3,\dots}$ and fix  $\beta\in Q^+_{e,n}$.
    The \textbf{cyclotomic KLR algebra} $\Rbeta$ of type~$\Gamma$,
    determined by the weights $(\Lambda,\beta)$, is the unital associative
    $\Z$-algebra with generators
    $\set{\psi_1,\dots,\psi_{n-1}} \cup \set{ y_1,\dots,y_n}
                               \cup \set{e(\bi)|\bi\in I^\beta}
    $
    and relations
    \bgroup
      \setlength{\abovedisplayskip}{1pt}
      \setlength{\belowdisplayskip}{1pt}
      \begin{align*}
        y_1^{(\Lambda, \alpha_{i_1})}e(\bi)&=0
          & e(\bi) e(\bj) &= \delta_{\bi\bj} e(\bi),
          & {\textstyle\sum_{\bi \in I^\beta}} e(\bi)&= 1,\\
        y_r e(\bi) &= e(\bi) y_r,
          &\psi_r e(\bi)&= e(s_r{\cdot}\bi) \psi_r,
          &y_r y_s &= y_s y_r,
      \end{align*}
      \begin{align*}
        \psi_r \psi_s = \psi_s \psi_r \text{ if }|r-s|>1,
          \qquad
          \psi_r y_s  = y_s \psi_r \text{ if }s \neq r,r+1,\\
        \psi_r y_{r+1} e(\bi)=(y_r\psi_r+\delta_{i_ri_{r+1}})e(\bi),
          \qquad
          y_{r+1}\psi_re(\bi)=(\psi_r y_r+\delta_{i_ri_{r+1}})e(\bi),
      \end{align*}
      \begin{align*}
          \psi_r^2e(\bi)&= Q_{i_ri_{r+1}}(y_r,y_{r+1})e(\bi)\\
        (\psi_r\psi_{r+1}\psi_r-\psi_{r+1}\psi_r\psi_{r+1})e(\bi)
           &=\delta_{i_ri_{r+1}}\frac{Q_{i_r,i_{r+1}}(y_{r+2},y_{r+1}) -
              Q_{i_r, i_{r+1}} (y_r, y_{r+1})}{y_r-y_{r+2}} e(\bi)
      \end{align*}
    \egroup
    for $\bi,\bj\in I^\beta$ and all admissible $r$ and $s$.
    Set $\Rn = \bigoplus_{\beta\in Q^+_{e,n}}\Rbeta$.
  \end{Definition}

  The reader can check that the right hand side of the last relation is
  a polynomial in $y_r,y_{r+1}$ and~$y_{r+2}$.  Importantly, the
  algebras $\Rbeta$ and $\Rn$ are $\Z$-graded, where the grading is
  determined by
  \[
     \deg e(\bi)=0,\quad
     \deg y_r = 2 \quad\text{and}\quad
     \deg\psi_re(\bi)=-(\alpha_{i_r},\alpha_{i_{r+1}}),
  \]
  for all admissible $\bi\in I^n$ and $1\le r\le n$.  If $A$ is any ring
  let $\Rbeta(A)=\Rbeta\otimes_\Z A$ and $\Rn(A)=\Rn\otimes_\Z A$ be the
  corresponding cyclotomic KLR algebras over~$A$.

  From the relations, $\Rn$ has a unique homogeneous
  anti-automorphism~$*$ of degree~$0$ that fixes all of the generators
  of~$\Rn$.

  Let $\Sym_n$ be the symmetric group of degree~$n$, which we consider
  as a Coxeter group with its standard set of Coxeter generators
  $\set{s_1,\dots,s_{n-1}}$, where $s_r=(r,r+1)$ for $1\le r<n$.  For
  each $w\in\Sym_n$ fix a \textbf{reduced expression} $w=s_{r_1}\dots
  s_{r_k}$ for $w$ (that is, fix such a word with $k$ minimal), and
  define $\psi_w = \psi_{r_1}\dots\psi_{r_k}\in\Rn$.  In general, the
  element $\psi_w$ depends on the choice of reduced expression but any
  fixed choice of reduced expression suffices for the results that
  follow.

  \subsection{Formal characters}
  Formal characters are a useful tool in the representation theory of the
  affine Hecke algebra. We need an analogue of these characters for
  $\Rn$-modules. To this end, let $\Acal=\Z[q,q^{-1}]$ and let
  $\AIq=\Acal[\bi\mid \bi\in I^n]$ be the free $\Acal$-module with basis
  $\set{\bi|\bi\in I^n}$.

  If $F$ is a field let $\Rep\Rn(F)$ be the Grothendieck group of
  finitely generated graded $\Rn(F)$-modules and if $M$ is an
  $\Rn$-module let $[M]$ be its image in $\Rep\Rn(F)$. Then $\Rep\Rn(F)$
  is the free $\Acal$-module with basis $\set{[\Dmu]|\bmu\in\Klesh}$,
  where $q$ acts on $\Rep\Rn(F)$ by grading shift, so that $[qM]=q[M]$,
  for $M\in\Rep\Rn(F)$.

  Let $\Ln$ be the positively graded commutative subalgebra of $\Rn$
  that is generated by $y_1,\dots,y_n$ and the idempotents $e(\bi)$, for
  $\bi\in I^n$. It is easy to see that the irreducible representations
  of $\Ln$ are indexed by the $n$-tuples $\bi\in I^n$ such that
  $e(\bi)\ne0$. Let $ \Rep\Ln$ be the Grothendieck group of finitely
  generated graded $\Ln$-modules where $q$ acts by grading shift, which
  we view as a free $\Acal$-submodule of $\AIq$.

  \begin{Definition}\label{D:qcharacter}
    The \textbf{graded (formal) character} is the linear map
    $\Chq\map{\Rep\Rn}{\Rep\Ln}$ given by
    \[
        \Chq M = \sum_{\bi\in I^n}\bigl(\dim_qMe(\bi)\bigr)\bi,
        \qquad\text{for }M\in\Rep\Rn.
    \]
  \end{Definition}

  The map $\Chq$ can be viewed as an exact functor, namely the
  restriction functor, from the category of (graded) $\Rn$-modules to
  the category of graded $\Ln$-modules. For our purposes it is enough to
  think of $\Chq$ as a linear map $\Chq\map{\Rep\Rn}\AIq$.

  \begin{Proposition}
    \label{P:qCharacters}
    Let $F$ be a field. Then the graded character map
    $\Chq\map{\Rep\Rn(F)}\AIq$ is injective.
  \end{Proposition}

  \begin{proof}
    Khovanov and Lauda~\cite[Theorem~5.17]{KhovLaud:diagI} state this
    result for the Grothendieck group of the corresponding affine KLR
    algebra, however, this immediately implies the result for the
    quotient cyclotomic KLR Hecke algebra~$\Rn(F)$.
  \end{proof}

  \subsection{Tableau combinatorics}
  Before we can define a homogeneous basis for $\Rn$, and introduce the
  graded Specht modules for~$\Rn$, we need to introduce the
  combinatorics of standard tableaux. This combinatorics plays a key
  role in proofs of our main results. Recall that we fixed integers $n$
  and $\ell$ in the paragraph before \autoref{D:KLRAlgebra}.

  A \textbf{partition} of~$m$ is a weakly decreasing sequence
  $\lambda=(\lambda_1,\lambda_2,\dots)$ of non-negative integers such
  that $|\lambda|=\lambda_1+\lambda_2+\dots=m$. An
  \textbf{$\ell$-partition} of~$n$ is an $\ell$-tuple
  $\blam=(\lambda^{(1)},\dots,\lambda^{(\ell)})$ of partitions such that
  $|\lambda^{(1)}|+\dots+|\lambda^{(\ell)}|=n$. We identify the
  $\ell$-partition $\blam$ with its \textbf{diagram}, which is the set of
  \textbf{nodes}
  $\set{(l,r,c)|1\le c\le \lambda^{(l)}_r\text{ for }1\le l\le\ell}$.
  We think of the diagram of $\blam$ as an $\ell$-tuple of
  left-justified arrays of boxes. This allows us to talk of the
  \textbf{components}, \textbf{rows} and \textbf{columns} of~$\blam$.

  The set of $\ell$-partitions of~$n$ is a poset under the
  \textbf{dominance} order $\gedom$, where~$\blam\gdom\bmu$ if
  $$\sum_{k=1}^{l-1}|\lambda^{(k)}|+\sum_{j=1}^i\lambda^{(l)}_j
       \ge\sum_{k=1}^{l-1}|\mu^{(k)}|+\sum_{j=1}^i\mu^{(l)}_j,$$
  for $1\le l\le\ell$ and $i\ge1$. If $\blam\gedom\bmu$ and $\blam\ne\bmu$
  then write $\blam\gdom\bmu$. Let $\Parts=\Parts[\ell,n]$ be the set of
  $\ell$-partitions of~$n$, which we consider as a poset under dominance.

  Fix $\blam\in\Parts$. A \textbf{$\blam$-tableau} is a bijective map
  $\t\map{\blam}\set{1,2,\dots,n}$, which we identify with a labelling of
  (the diagram of) $\blam$ by $\{1,2,\dots,n\}$. Let $\Shape(\t)=\blam$
  be the \textbf{shape} of $\t$. For example,
  \begin{equation}
    \Tritab({1,2,3},{4}|{5},{6}|{7,8,9},)
    \hspace*{2mm}\text{and}\hspace*{2mm}
    \Tritab({3,5,8},{6}|{1},{7}|{2,4,9},)
  \end{equation}\label{E:Tableaux}
  are both $\blam$-tableaux, where $\blam=(3,1|1^2|3)\in\Parts[13]$,
  where $\ell=3$. A $\blam$-tableau is \textbf{standard} if, in each
  component, its entries increase along rows and down columns. For
  example, both of the tableaux above are standard. Let $\Std(\blam)$ be
  the set of standard $\blam$-tableaux. If $\Pcal$ is any set of
  $\ell$-partitions let
  $\Std(\Pcal)=\bigcup_{\blam\in\Pcal}\Std(\blam)$. Similarly set
  $\Std^2(\Pcal)=\set{(\s,\t)|\s,\t\in\Std(\blam) \text{ for
  }\blam\in\Pcal}$.

  If $\t$ is a $\blam$-tableau set $\Shape(\t)=\blam$ and let
  $\t_{\downarrow m}$ be the subtableau of~$\t$ that contains the
  numbers $\{1,2,\dots,m\}$. If~$\t$ is a standard $\blam$-tableau then
  $\Shape(\t_{\downarrow m})$ is a $\ell$-partition for all $m\ge0$.
  Extend the dominance ordering to $\Std(\Parts)$ by defining
  \[
  \s\gedom\t \quad\text{ if }
     \Shape(\s_{\downarrow m})\gedom\Shape(\t_{\downarrow m})
     \text{ for }1\le m\le n.
  \]
  As before, write $\s\gdom\t$ if $\s\gedom\t$ and $\s\ne\t$.

  Even though we write $\Parts$ for the set of $\ell$-partitions of~$n$,
  so far none of the definitions in this section depend on
  $\Lambda\in P_e^+$.  A \textbf{multicharge} for~$\Lambda$ is an
  $\ell$-tuple $\charge=(\kappa_1,\dots,\kappa_\ell)\in\Z^\ell$ such that
  \begin{equation}\label{E:multicharge}
    (\Lambda,\alpha_i)= \#\set{1\le l\le\ell|\kappa_l\equiv i\pmod e},
    \qquad\text{ for all }i\in I.
  \end{equation}
  The \textbf{residue} of a node $A=(l,r,c)\in\blam$ is $\res(A) =
  \kappa_l+c-r +e\Z\in I$. If $\t$ is a standard $\blam$-tableau and
  $\t(A)=m$, where $1\le m\le n$, then the residue of~$m$ in~$\t$ is
  $\res_m(\t) = \res(A)$. The \textbf{residue sequence} of~$\t$ is
  $\res(\t) =\(\res_1(\t),\res_2(\t),\dots,\res_n(\t)\)\in I^n$. Set
  \[
     \Std_\bi(\blam)=\set{\t\in\Std(\blam)|\res(\t)=\bi})
     \quad\text{and}\quad
     \Std_\bi(\Parts) =\bigcup_{\blam\in\Parts}\Std_\bi(\blam).
  \]

  A node $A$ is an \textbf{addable node} of~$\blam$ if $\res(A)=i$ and
  $A\notin\blam$ and $\blam\cup\{A\}$ is the (diagram of) a
  $\ell$-partition of~$n+1$.  Similarly, a node $B$ is a \textbf{removable
  node} of~$\blam$ if $\res(B)=i$ and $B\in\blam$ and
  $\blam\setminus\{B\}$ is a $\ell$-partition of~$n-1$. Let
  $\Add_i(\blam)$ and $\Rem_i(\blam)$ be the sets of addable and
  removable $i$-nodes of~$\blam$. Let $\le$ be the lexicographic order n
  the set of nodes. If $A$ is a removable node of~$\blam$ define
  \[
    d_{A,e}(\blam)=\#\set{B\in\Add_i(\blam)|A<B}-\#\set{B\in\Rem_i(\blam)|A<B}.
  \]
  Following Brundan, Kleshchev and Wang~\cite[\S3.5]{BKW:GradedSpecht},
  define the $e$-\textbf{degree} of a standard tableau $\t$ by
  \begin{equation}\label{E:TableauDegree}
      \deg_e(\t) = \begin{cases*}
        \deg_e(\t_{\downarrow(n-1)})+d_{A,e}(\bmu),& if $n>0$,\\
           0,& if $n=0$,
        \end{cases*}
  \end{equation}
  where $A=\t^{-1}(n)$ is the node in $\t$ that contains~$n$,
  For convenience, set $\deg_0(\t)=0$ for all $\t$.

  There is  unique standard $\blam$-tableaux $\tlam$ such that
  $\tlam\gedom\t$, for all $\t\in\Std(\blam)$.  The tableau $\tlam$ has
  the numbers $1,2,\dots,n$ entered in order from left to right along
  the rows of $\t^{\lambda^{(1)}}$, and then
  $\t^{\lambda^{(2)}},\dots,\t^{\lambda^{(\ell)}}$. For example, the
  first tableau in \autoref{E:Tableaux} is $\tlam$ for
  $\blam=(3,1|1^2|3)$. Given a standard
  $\blam$-tableau $\t$ define a permutation $d(\t)\in\Sym_n$ by
  $\t=\tlam d(\t)$.

  \subsection{Specht modules, simple modules and almost simple modules}
  \label{S:GradedCellular}
  The main results in this paper give explicit, cancellation free,
  descriptions of the Jantzen sum formula for the (ungraded) Specht
  modules. This section defines the modules that we use to describe
  these filtrations.

  As we now recall, the algebra $\Rn$ is a \textbf{graded cellular
  algebra} in the sense of \autoref{S:Cellular}.

  For $\blam\in\Parts$ set $\bi^{\blam}=\res(\tlam)=(i^\blam_1,\dots,i^\blam_n)$.
  For $1\le m\le n$ set
  $\blam_{\downarrow m}=\Shape(\tlam_{\downarrow m})$
  and let $\alpha^\blam_m\in\blam$ be the node such that
  $\tlam(\alpha^\blam_m)=m$.  Define the polynomial
  \[
      y^\blam = \prod_{k=1}^ny_r^{a_k(\blam)},\qquad
        \text{where }
        a_k(\blam)=\#\set{\alpha\in\Add_{i^\blam_k}(\blam_{\downarrow k})|%
                 \alpha<\alpha^\blam_k}.
  \]
   Then $y^\blam$ is homogeneous of degree $2\deg\tlam$. For
   $(\s,\t)\in\Std^2(\blam)$ define
   \[
                \psi_{\s\t} = \psi_{d(\s)}^*
                e(\bi^\blam)y^\blam\psi_{d(\t)},
   \]
   By construction, $\psi_{\s\t}$ is homogeneous and one can show that
   $\deg_e\psi_{s\t}=\deg_e\s+\deg\t$. Jun Hu and the
   author~\cite{HuMathas:GradedCellular} proved that these elements
   give a graded cellular basis for~$\Rn(F)$, when $F$ is a field.
   Building on this result, Ge Li~\cite{GeLi:IntegralKLR} proved that
   $\Rn$ is a $\Z$-free graded cellular algebra.

   \begin{Theorem}[Li~\cite{GeLi:IntegralKLR}]\label{T:GradedCellular}
     The algebra $\Rn$ is free as an $\Z$-module of rank $\ell^n n!$.
     Moreover, $\Rn$ is a graded cellular algebra with homogeneous cellular
     basis $\set{\psi_{\s\t}|(\s,\t)\in\Std^2(\Parts)}$.
  \end{Theorem}

  Subsequently, Kang and
  Kashiwara~\cite[Theorem~4.5]{KangKashiwara:CatCycKLR} generalised this
  result and proved that the cyclotomic KLR Hecke algebras indexed by
  symmetrisable Cartan matrices are $\Z$-free. For this paper, the fact
  that $\Rn(\Z)$ is a cellular algebra is very important.

  Since $\Rn$ is a cellular algebra, if $a\in\Rn$ and
  $(\s,\t)\in\Std^2(\blam)$, for $\blam\in\Parts$, then we can write
  \[
    \psi_{\s\t}a \equiv \sum_{\u\in\Std(\blam)}r_\u\psi_{\u\t} \pmod\Rlam,
  \]
  where the scalar $r_\u\in\Z$ depends only on $a$, $\s$ and $\u$ (and
  not~$\t$) and where $\Rlam$ is the two-sided ideal of~$\Rn$ spanned by
  $\set{\psi_{\u\v}|(\u,\v)\in\Std(\bmu)\text{ where }\bmu\gdom\blam}$.
  Define the \textbf{graded Specht module} $\Slam_\Z$ to the free $\Z$-module
  with homogeneous basis $\set{\psi_\s|\s\in\Std(\blam)}$, where
  $\deg\psi_\s=\deg_e\s$, and with $\Rn$-action
  \begin{equation}\label{E:}
    \psi_{\s}a = \sum_{\u\in\Std(\blam)}r_\u\psi_{\u},
    \qquad\text{ for $a\in\Rn$ and $\s\in\Std(\blam)$.}
  \end{equation}\label{E:GradedSpecht}
  In particular, it follows that if $\bi\in I^n$ then
  $\psi_\s e(\bi)= \delta_{\bi\bi^\s}\psi_\s$.

  For any ring $F$ let $\Slam_F=\Slam_\Z\otimes_\Z F$. Then $\Slam_F$ is a
  graded $\Rn(F)$-module. We abuse notation and write~$\psi_\s$, instead
  of $\psi_s\otimes_\Z1_F$, for the basis elements of~$\Slam_F$.

  The cellular algebra axioms imply that $\Slam_F$ has a
  bilinear form $\<\ ,\ \>\map{\Slam_F\times\Slam_F}F$ determined by
  \begin{equation}
    \<\psi_\s,\psi_\t\>\psi_\u = \psi_{\s\t}\psi_\u,
    \qquad\text{ for }\s,\t,\u\in\Std(\blam).
  \end{equation}\label{E:pechtForm}
  This form is associative in the sense that
  $\<ax,y\>=\<x,a^*y\>$, for all $a\in\Rn$ and $x,y\in\Slam_F$. Hence,
  \[
      \rad\Slam_F = \set{x\in\Slam_F|\<x,y\>=0\text{ for all }y\in\Slam_F}
  \]
  is a graded $\Rn(F)$-submodule of $\Slam_F$. Observe that these
  definitions make sense for any ring and, in particular, these
  definitions are valid even over $F=\Z$.

  \begin{Definition}\label{D:Simples}
      Let $\bmu\in\Parts$. Define
      $\Emu_\Z = \Slam[\bmu]_\Z/\rad \Slam[\bmu]_\Z$. For any ring $F$
      define
      \[ \Emu_F = \Emu_\Z\otimes_\Z F\qquad\text{and}\qquad
         \Dmu_F = \Slam_F/\rad \Slam_F.
      \]
  \end{Definition}

  By definition, $\Emu_\Z=\Dmu_\Z$ but in positive characteristic the
  $\Rn(F)$-modules $\Emu_F $ and $\Dmu_F$ are not necessarily
  isomorphic. By the theory of (graded) cellular algebras
  \cite[Corollary~2.11]{HuMathas:GradedCellular}, if $F$ is a field then
  the module~$\Dmu_F$ is  either zero or absolutely irreducible and self-dual.
  Moreover, up to grading shift, every irreducible $\Rn(F)$-module
  arises uniquely in this way.

  Let $F$ be a field and set
  $\Klesh=\Klesh[e,n]=\set{\bmu\in\Parts|\Dmu_F\ne0}$. Then
  $\set{\Dmu_F\<d\>|\bmu\in\Klesh\text{ d }d\in\Z}$ is a complete set of
  pairwise non-isomorphism irreducible graded $\Rn$-modules by
  \cite[Corollary~5.11]{HuMathas:GradedCellular}. (The graded irreducible
  ~$\Rn(F)$-modules were first classified in
  \cite[Theorem~5.10]{BK:GradedDecomp}.) Moreover, the set $\Klesh$
  is independent of the field~$F$. For $\blam\in\Parts$ and
  $\bmu\in\Klesh$ define the \textbf{graded decomposition number}
  \begin{equation}\label{E:GradedDecomp}
     d^{F,e}_{\blam\bmu}(q) = [\Slam_F:\Dmu_F]_q
         = \sum_{d\in\Z}[\Slam_F:\Dmu_F\<d\>]q^d\in\N[q,q^{-1}].
  \end{equation}
  In particular, $d^{\C,e}_{\blam\bmu}(q)$ is a
  graded decomposition number for $\Rn(\C)$ in characteristic zero. The
  \textbf{graded decomposition matrix} of $\Rn(F)$ is the unitriangular matrix
  $\mathsf{D}_{F,e}(q)=\bigl(d^{F,e}_{\blam\bmu}(q)\bigr)$, where
  $\blam\in\Parts$ and $\bmu\in\Klesh$ and where the rows are columns
  are ordered lexicographically.

  The following result plays a key role in this paper.  We give the
  proof, because we need the underlying ideas below and because the
  proof is quite short, with the original reference
  \cite{Mathas:Singapore} not being readily available. Given integers
  $d$ and $d'$ write $d\mid d'$ if $d'=ad$, for some $a\in\Z$. In
  particular, note that $d\mid 0$, for all $d\in\Z$.

  \begin{Theorem}[\protect{%
    \cite[Theorem~3.7.4 and Theorem~3.7.5]{Mathas:Singapore}}]
    \label{T:AdjustmemtMatrix}
      Suppose that $F$ is a field and $\bmu\in\Klesh$. Then:
      \begin{enumerate}
        \item The module $\Emu_\Z$ is a $\Z$-free $\Z$-graded $\Rn(\Z)$-module
        \item If $F=\Q$ then $\Emu_\Q\cong\Dmu_\Q$ is an absolutely
        irreducible self-dual graded $\Rn(\Q)$-module
        \item If $\blam\in\Parts$ then
        $[\Slam_F:\Dmu_F]_q=\sum_\bnu d^{\C,e}_{\blam\bnu}(q)[\Emu[\bnu]_F:\Dmu_F]_q$.
      \end{enumerate}
  \end{Theorem}

  \begin{proof}
    Let $G^\blam=\bigl(\<\psi_\s,\psi_\t\>\bigr)_{\s,\t\in\Std(\mu)}$ be the
    Gram matrix of the integral graded Specht module $\Slam$. Let
    $N=\#\Std(\bmu)$ so that $\Slam_\Z$ is a free $\Z$-module of rank~$N$.
    Since~$\Z$ is a principal ideal domain there exist positive integers
    $d_1,\dots,d_N$ such that $d_1\mid d_2\mid \dots\mid d_N$ and two homogeneous
    bases $\set{u_i}$ and $\set{v_i}$ of the graded Specht module
    $\Slam_\Z$ such that $\<u_i,v_j\>=\delta_{ij}d_i$, for $1\le i,j\le N$.
    That is, the diagonal matrix $\diag(d_1,\dots,d_N)$ is the
    \textit{Smith normal form} of~$G^\blam$. Therefore, $\rad\Slam_\Z$ is free
    as an $\Z$-module with basis $\set{u_i|d_i=0}$ and $\Emu_\Z$ is free
    as a $\Z$-module with basis $\set{u_i+\rad\Slam_\Z|d_i\ne0}$, proving~(a).
    (In the same way, observe that $\set{v_i|d_i=0}$ is a basis of $\rad\Slam_Z$ and
    $\set{v_i+\rad\Slam_\Z|d_i\ne0}$ is a basis of $\Emu_\Z$.)

    By definition, $\Emu_\Q=\Dmu_\Q$, so part~(b) follows from the
    general theory of (graded) cellular algebras; specifically, see
    \cite[Theorem~2.10]{HuMathas:GradedCellular}. The same result says
    that if $F$ is a field then $\Dmu_\F$ is a self-dual
    absolutely irreducible graded $\Rn$-module.

    Let $F$ be a field and consider~(c). Taking graded characters
    is exact and commutes with base change, so
    \begin{align*}
         \Chq\Slam_F & = \Chq\Slam_\Q
            = \sum_{\bnu\in\Klesh} d^{\C,e}_{\blam\bnu}(q)\Chq\Dmu[\bnu]_\Q
            = \sum_{\bnu\in\Klesh} d^{\C,e}_{\blam\bnu}(q)\Chq\Emu[\bnu]_\Q\\
           &= \sum_{\bnu\in\Klesh} d^{\C,e}_{\blam\bnu}(q)\Chq\Emu[\bnu]_F
            = \sum_{\bnu\in\Klesh} d^{\C,e}_{\blam\bnu}(q)
                \sum_{\bmu\in\Klesh} [\Emu[\bnu]_F:\Dmu_F]_q\Chq\Dmu_F\\
           &= \sum_{\bmu\in\Klesh} \Bigl(\sum_{\bnu\in\Klesh}
              d^{\C,e}_{\blam\bnu}(q)[\Emu[\bnu]_F:\Dmu_F]_q\Bigr) \Chq\Dmu_F.
    \end{align*}
    On the other hand,
    $\Chq\Slam_F=\sum_\bmu d^{F,e}_{\blam,\bmu}(q)\Chq\Dmu_F$, so part~(c)
    follows by \autoref{P:qCharacters}.
  \end{proof}

  Let $a^{F,e}_{\bnu\bmu}(q)=[\Emu[\bnu]_F:\Dmu_F]_q$, for
  $\bmu,\bnu\in\Klesh$.
  The matrix
  $\mathsf{A}_{F,e}(q)=\(a^{F,e}_{\bnu\bmu}(q)\)_{\bmu,\bnu\in\Klesh}$ is the
  (graded) \textbf{adjustment matrix} of~$\Rn(F)$. Notice that
  \autoref{T:AdjustmemtMatrix}(c) is equivalent to the following
  factorisation of the graded decomposition matrix:
  \begin{equation}\label{E:Adjustment}
    \mathsf{D}_{F,e}(q) = \mathsf{D}_{\C,e}(q)\mathsf{A}_{F,e}(q).
  \end{equation}
  The adjustment matrix $\mathsf{A}_{F,e}(q)$ coincides with the matrix
  defined by Brundan and Kleshchev \cite[Theorem~5.17]{BK:GradedDecomp}
  using different arguments. The arguments of \cite{BK:GradedDecomp}
  amount to carefully choosing a $\Z$-lattice for the simple module
  $\Dmu_\Q$. The proof of \autoref{T:AdjustmemtMatrix} shows that the
  $\psi$-basis, by virtue of \autoref{T:GradedCellular}, automatically
  chooses a $\Z$-lattice $\Emu_\Z$ for $\Dmu_\Q$.

  \begin{Corollary}\label{C:SlamEmuCharacters}
    Suppose that $\blam\in\Parts$. Then
    $\Chq\Slam_F = \Sum_{\bmu\in\Klesh}d^{\C,e}_{\blam\bmu}(q)\Chq\Emu_F$.
  \end{Corollary}

  Let $\overline{\phantom{X}}\map{\Z[q,q^{-1}]}{\Z[q,q^{-1}]}$ be the unique
  linear involution such that $\overline{q}=q^{-1}$. Abusing notation,
  extend this to a map $\overline{\phantom{X}}\map{\Rep\Ln}{\Rep\Ln}$.

  \begin{Corollary}\label{C:BarInvariance}
     Suppose that $\bmu,\bnu\in\Klesh$. Then
     $\overline{\Chq\Emu[\bnu]_F}=\Chq\Emu[\bnu]_F$ and
     $\overline{a^{F,e}_{\bnu\bmu}(q)}=a^{F,e}_{\bnu\bmu}(q)$.
  \end{Corollary}

  \begin{proof}
    By construction and \autoref{T:AdjustmemtMatrix}(b),
    $\Chq\Emu[\bnu]_F=\Chq\Emu[\bnu]_\Q=\Chq\Dmu[\bnu]_\Q$. Therefore,
    $\overline{\Chq\Emu[\bnu]_F}=\Chq\Emu[\bnu]_F$ since $\Dmu[\bnu]_\Q$ is a
    self-dual graded $\Rn(\Q)$-module. This implies the second claim
    because
    \[
        \Chq\Emu[\bnu]_F = \sum_{\bnu\in\Klesh}a^{F,e}_{\bnu\bmu}(q)\Chq\Dmu_F
    \]
    and $\Dmu_F$ is self-dual.
  \end{proof}

  Let $\bnu\in\Klesh$. Since $[\Emu[\bnu]_F:\Dmu[\bnu]_F]_q=1$ it
  follows that $\Emu[\bnu]_F$ is self-dual if and only if
  $\Emu[\bnu]_F=\Dmu[\bnu]_F$.

\section{Cyclotomic Hecke algebras}
  We are now ready to introduce the cyclotomic Hecke algebras of
  type~$A$, and their Specht modules, which are the main focus of this
  paper. We define several different bases for these algebras and use
  them to compute the characters of the sum of modules in the Jantzen
  filtration of the Specht modules.

  \subsection{Cyclotomic Hecke algebras of type~A}
  \label{S:Hecke}
  Let $R$ be a commutative ring with~$1$ and fix non-negative integers
  $n,\ell\ge0$. If $v\in R$ and $k\in\Z$ define the $v$-\textbf{quantum
  integer}
  \[       [k]_v = \begin{cases*}
                 1+v+\dots+v^{k-1},&if $k\ge0$,\\
                 -v^{-1}+v^{-2}-\dots-v^{-k},&otherwise.
            \end{cases*}
  \]
  Observe that if $v\ne1$ then $[k]_v=(v^k-1)/(v-1)$ and if $v=1$ then
  $[k]_v=k$.

  We use the following definition of the cyclotomic Hecke algebras of
  type~$A$, following \cite[\S2]{HuMathas:SeminormalQuiver}.

  \begin{Definition}\label{D:CyclotomicHecke}
      The \textbf{cyclotomic Hecke algebra}of type~$A$ with
      \textbf{Hecke parameter} $v\in R^\times$ and
      \textbf{$\ell$--charge}~$\charge=(\kappa_1,\dots,\kappa_\ell)\in\Z^\ell$ is
      the unital associative $R$-algebra $\Hkn=\Hkn(A)$ with
      generators $T_1,\dots,T_{n-1}$, $L_1,\dots,L_n$ and relations
      \begin{align*}
          \textstyle\prod_{l=1}^\ell(L_1-[\kappa_l]_v)&=0,  &
          (T_r-v)(T_r+1)&=0,  \\
          L_rL_t&=L_tL_r, &
        T_rT_s&=T_sT_r &\text{if }|r-s|>1,\\
        T_sT_{s+1}T_s&=T_{s+1}T_sT_{s+1}, &
        T_rL_t&=L_tT_r,&\text{if }t\ne r,r+1,\\
        \span\span L_{r+1}=T_rL_rT_r+T_r,\span\span\span
      \end{align*}
      where $1\le r<n$, $1\le s<n-1$ and $1\le t\le n$.
  \end{Definition}

  As is explained in \cite[\S2]{HuMathas:SeminormalQuiver}, if
  $v\ne1$ these algebras are isomorphic to the Ariki-Koike algebras,
  which were introduced in~\cite{AK}, and when $v=1$ they are
  isomorphic to the degenerate Ariki-Koike algebras.

  By definition, $\Hkn$ has a unique anti-isomorphism $*$ that fixes
  every generator.

  Recall from \autoref{SS:CyclotomicKLR} that $\Sym_n$ is the symmetric
  group of degree~$n$ with distinguished Coxeter generators
  $s_1,\dots,s_{n-1}$. If $w\in\Sym_n$ let $L(w)$ be the minimal length
  of $w$ as a product of the Coxeter generators. A word $w=s_{i_1}\dots
  s_{i_k}$ is a \textbf{reduced expression} for $w$ if $k=L(w)$, for
  $1\le i_j<n$. As the braid relations hold in $\Hn$, if $w=s_{i_1}\dots
  s_{i_k}$ is reduced then $T_w=T_{i_1}\dots T_{i_k}\in\Hn$ depends only
  on $w$ and not on the choice of reduced expression. These results are
  well-known and their proofs can be found, for example, in
  \cite[Chapter~1]{Mathas:ULect}.

  By the results of \cite{AMR:NazarovWenzl,DJM:cyc}, and in view of
  \cite[Theorem~2.8]{HuMathas:SeminormalQuiver}, $\Hkn$ is cellular
  algebra with cellular basis indexed by~$\Std^2(\Parts)$, as in
  \autoref{T:GradedCellular}. Hence, exactly as before (except that we
  no longer have a grading), for each $\ell$-partition $\blam\in\Parts$
  we have a Specht module $S^\blam$, which is a left $\Hkn$-module.

  Let $F$ be a field. A non-zero element $\xi\in F$ has \textbf{quantum
  characteristic}~$e$ if $e$ is the smallest positive integer such that
  $[e]_\xi=0$ --- set $e=\infty$ if $[k]_\xi\ne0$, for $k>0$.  Given
  $\charge\in\Z^\ell$ define a dominant weight
  $\Lambda=\Lambda_\charge\in P_e^+$ by using \autoref{E:multicharge}.
  Define $\Hn(F)=\Hkn[\xi](F)$.  Unravelling the notation, the
  cyclotomic relation in $\Hn(F)$ can be written as
  \begin{equation}\label{E:CyclotomicRelation}
            0 = \prod_{l=1}^\ell(L_1-[\kappa_l]_\xi)
              = \prod_{i\in I}(L_1-[i]_\xi)^{(\Lambda,\alpha_i)},
  \end{equation}
  where $[i]_\xi$ has the obvious meaning, for $i\in I$.

  We can now give the relationship between the cyclotomic KLR and Hecke
  algebras of type~$A$.

  \begin{Theorem}[Brundan and Kleshchev~\cite{BK:GradedDecomp},
    Rouquier~\cite{Rouquier:QuiverHecke2Lie}]\label{T:KLRIsomorphism}
      Let $F$ be a field. Suppose that $\xi\in F$ has quantum
      characteristic~$e$ and that $\Lambda\in P_e^+$. Then
      $\Hn(F)\cong\underline{\Rn(F)}$ as ungraded algebras.
  \end{Theorem}

  Notice that \autoref{T:KLRIsomorphism} implies that, up to
  isomorphism, $\Hn(F)$ depends only on~$e$, $\Lambda$ and $F$, and not
  on the particular choice of $\xi\in F$ with quantum
  characteristic~$e$. As explained in \autoref{S:SeminormalForms} below,
  we will always assume that~$e$ is finite.

  \subsection{The Murphy basis}
    Following \cite{DJM:cyc}, this section defines a cellular basis of
    $\Hn$ that we will use to define Specht modules for $\Hn$.

    Let $\blam\in\Parts$ and define elements
    \[
      m_\blam = u_\blam\sum_{w\in\Sym_\blam}T_w,
      \qquad\text{where }\qquad
        \prod_{l=2}^\ell\prod_{m=1}^{|\lambda^{(1)}|+\dots+|\lambda^{(l-1)}|}
             (L_m-[\kappa_l]).
    \]
    For $\s,\t\in\Std(\blam)$ define $m_{\s\t}=T_{d(\s)}^* m_\blam T_{d(\t)}$.
    By definition, $m_{\s\t}^*=m_{\t\s}$.

   \begin{Theorem}[{Dipper, James, Mathas\cite[Theorem~3.26]{DJM:cyc}}]
     \label{T:MurphyBasis}
     For any ring $R$, $\Hn$ is an $R$-free cellular algebra with cellular basis
     $\set[\big]{m_{\s\t}|\s,\t\in\Std(\blam)\text{ for }\blam\in\Parts}$.
   \end{Theorem}

   \begin{Remark}
      The paper \cite{DJM:cyc} only proves \autoref{T:MurphyBasis} in
      the case when $\xi\ne1$. However, the argument given in
      \cite{DJM:cyc} extends without change to the case when $\xi=1$
      because the result is bootstrapped from the Murphy basis of the
      Iwahori-Hecke algebra of the symmetric group~\cite{Murphy:basis},
      which is a cellular basis for all~$\xi$
      by~\cite[Chapter~3]{Mathas:ULect}. Based on \autoref{T:MurphyBasis},
      all of the results in this paper hold for the degenerate
      cyclotomic Hecke algebras of type~$A$, which correspond to the
      case when~$\xi=1$ (and $e$ is the characteristic of~$R$).
   \end{Remark}

   As in \autoref{S:Cellular}, for each $\blam\in\Parts$ let $S^\blam$
   be the Specht module determined by the cellular basis of
   \autoref{T:MurphyBasis}. Then $S^\blam$ is free as an $R$-module with basis
   $\set{m_\t|\t\in\Std(\blam)}$ and it comes equipped with an inner
   product $\<\, ,\ \>$ that is uniquely determined by
   \[
     m_\s m_{\t\v} =\<m_\s, m_\t\> m_\v,\qquad\text{for }\s,\t, \v\in\Std(\blam).
   \]

  Recall from \autoref{E:GradedDecomp} that
  $d^{F,e}_{\blam\bmu}(q)=[\Slam_F:\Dmu_F]_q$ is a graded decomposition
  number of~$\Rn(F)$, for $\blam\in\Parts$ and $\bmu\in\Klesh$.
  By \cite{BKW:GradedSpecht} or \cite[\S5.2]{HuMathas:GradedCellular},
  under the isomorphism of \autoref{T:KLRIsomorphism} the modules
  $S^\blam_F$ and $\Slam_F$ coincide, once we forget the grading.
  In more detail, we have:

  \begin{Corollary}\label{C:Ungrading}
    Let $\blam\in\Parts$ and $\bmu\in\Klesh$. Then
    $S^\blam_F \cong\underline{\mathbb{S}}^\blam_F$ and
    $D^\bmu_F  \cong\underline{\mathbb{D}}^\bmu_F$. Consequently,
    \[ [S^\blam_F:D^\bmu_F] = d^{F,e}_{\blam\bmu}(1). \]
  \end{Corollary}

  Let $\mathsf{D}_{F,e}=\bigl(d^{F,e}_{\blam\bmu}\bigr)$ be the
  \textbf{decomposition matrix} of~$\Hn(F)$. Then
  $\mathsf{D}_{F,e}=\mathsf{D}_{F,e}(1)$ by \autoref{C:Ungrading}, where
  $\mathsf{D}_{F,e}(1)=\mathsf{D}_{F,e}(q)|_{q=1}$ is the graded
  decomposition matrix of \autoref{SS:CyclotomicKLR} evaluated at $q=1$.

  \subsection{Seminormal forms}\label{S:SeminormalForms}
  This section recalls the construction of the seminormal basis of
  $\Hn(\Kcal)$ in the semisimple case, so all of these results in this
  section, ultimately, go back to the seminal work of Young~\cite{QSAI}.
  We follow the general framework from \cite{Mathas:Seminormal} but see
  also \cite{AK,HuMathas:SeminormalQuiver,Mathas:gendeg}.

  For the rest of this paper we assume that $e<\infty$ in
  \autoref{D:KLRAlgebra} or, equivalently, that $\xi$ is not a root of
  unity in \autoref{D:CyclotomicHecke}. By
  \cite[Corollary~2.10]{HuMathas:SeminormalQuiver}, up to
  isomorphism, we can always assume that~$e$ is finite, so there is no
  loss of generality in making this assumption.  The advantage of
  assuming that~$e$ is finite is that this is allows us to adjust the
  multicharge of~\autoref{E:multicharge} modulo~$e$ in
  \autoref{E:Separation} below.

  Let $F$ be a field. As in \autoref{E:multicharge} we fix an
  $\ell$--charge $\charge=(\kappa_1,\dots,\kappa_\ell)\in\Z^\ell$.  For
  the rest of this paper, \textit{we impose the additional assumption}
  that
  \begin{equation}\label{E:Separation}
    \kappa_1-2(\ell-1)n\ge\dots\ge\kappa_{\ell-1}-2n\ge\kappa_\ell\ge n.
  \end{equation}
  Implicitly, this assumption fixes a choice of lattice in a modular
  system for $\Hn(F)$ and, in principle, the Jantzen filtrations that we
  construct below depend on the choices in~\autoref{E:Separation}. In
  practice, we need \autoref{E:Separation} for the results on seminormal
  forms that we use below because \autoref{E:Separation} ensures that
  standard tableaux are uniquely determined by their content sequences.
  This said, assumption \autoref{E:Separation} is only a technical
  convenience because, by \autoref{E:multicharge}, the dominant
  weight $\Lambda=\Lambda_\charge$, and hence the algebra $\Hn(F)$,
  depend only on the \textit{multiset}
  $\set{\overline\kappa_1,\dots,\overline\kappa_\ell}$. In particular,
  $\Hn(F)$ does not depend on~\autoref{E:Separation}.

  As above, consider the cyclotomic Hecke algebra $\Hn=\Hn(F)$ with
  Hecke parameter $\xi\in F$ of quantum characteristic~$e<\infty$ and
  cyclotomic parameters $[\kappa_m]_\xi$, for $1\le m\le\ell$, where
  $\charge=(\kappa_1,\dots,\kappa_\ell)\in\Z^\ell$ is a fixed choice of
  $\ell$--charge.

  Let $x$ be an indeterminate over~$F$ and consider the localisation
  $\Ocal=F[x]_{(x)}$ of the polynomial ring $F[x]$ at the prime ideal
  $(x)=xF[x]$.  Let $\Kcal=F(x)$ be the field of fractions of~$F[x]$.
  Let $\HO$ be the cyclotomic Hecke algebra over~$\Ocal$ with Hecke
  parameter $z=x+\xi$ and with cyclotomic parameters
  $([\kappa_1]_z,\dots,[\kappa_\ell]_z)$. Let
  $\HK=\HO\otimes_\Ocal\Kcal$. Then $\HK$ is a split semisimple algebra.
  Consider the field~$F$ as an $\Ocal$-algebra where $x$ acts as~$0$.
  Then $\Hn(F)\cong\HO\otimes_\Ocal F$.

  For $\blam\in\Parts$ let $\SOlam$ be the (ungraded) Specht module for
  $\HO$ determined by the $\set{m_{\s\t}}$ cellular basis and $\blam$.
  Then $S^\blam_\Kcal=\SOlam\otimes_{\Ocal}\Kcal$ is the Specht module
  for $\HK$.

  Let $\t=(\t^{(1)}|\dots|\t^{(\ell)})\in\Std(\blam)$, for
  $\blam\in\Parts$.  Suppose that the integer $m$, with $1\le m\le n$,
  appears in row~$r$ and column~$c$ of $\t^{(k)}$. The \textbf{content}
  of~$m$ in~$\t$ is $\con_m(\t)=\kappa_m+c-r$. Observe that the residue
  of~$m$ is~$\res_m(\t)=\con_m(\t)+e\Z\in I$. The \textbf{content sequence}
  of~$\t$ is~$\con(\t)=\bigl(\con_1(\t),\dots,\con_n(\t)\bigr)$. More
  generally, the content of the node $(m,r,c)$ is
  $\con(m,r,c)=\kappa_m+c-r$. The key point of these definitions is that
  $\s=\t$ if and only if~$\con(\s)=\con(\t)$, for
  $\s,\t\in\Std(\Parts)$. This is easily proved by induction on~$n$
  using \autoref{E:Separation}.\footnote{This is the only place where
  assumption \autoref{E:Separation} is needed in the construction of the
  seminormal basis. That is, in this section \autoref{E:Separation} is
  only used to ensure that the content sequences separate the standard
  tableaux.  The stronger form of \autoref{E:Separation} is used to
  define sequences of charged beta numbers in \autoref{S:classical}.}

  James and the author~\cite[Proposition~3.7]{JamesMathas:cycJantzenSum}
  proved the following fundamental property of this basis:
  \begin{equation}\label{E:msLk}
        m_\s L_k = [c_k(\s)] m_\s+\sum_{\u\gdom\s}r_\u m_\u,\qquad
            \text{ for $\s\in\Std(\blam)$, $1\le k\le n$ and some }r_\u\in R.
  \end{equation}

  Let $\s,\t\in\Std(\blam)$, for $\blam\in\Parts$. Define
  \[
  f_{\s\t} = F_\s m_{\s\t} F_\t,\qquad\text{where}\quad
        F_\t = \prod_{k=1}^n\prod_{\substack{\s\in\Std(\Parts)\\ \con_k(\s)\ne \con_k(\t)}}
               \frac{L_k-[\con_k(\s)]_z}{[\con_k(\t)]_z-[\con_k(\s)]_z}.
  \]
  Equation \autoref{E:msLk} implies that
  \begin{equation}\label{E:fstExpansion}
    f_{\s\t} \equiv m_{\s\t} + \sum_{\u\gdom\s, \v\gdom\t}
    r_{\u\v}m_{\u\v}\pmod\Hlam,\qquad\text{for some }r_{\u\v}\in \Kcal.
  \end{equation}
  In particular, $\set{f_{\s\t}}$ is a basis of $\HK$. Since the
  transition matrix between the two bases $\set{m_{\s\t}}$ and
  $\set{f_{\s\t}}$ is unitriangular it is not hard to see that
  $\set{f_{\s\t}}$ is a cellular basis of $\HK$.  Moreover, if
  $\blam\in\Parts$  and $\s,\t\in\Std(\blam)$ then
  $S^\blam_\Kcal\cong f_{\s\t}\HK$.

  In fact, $\set{f_{\s\t}}$ is a \textbf{seminormal basis} of $\HK$ in
  the sense of \cite{HuMathas:SeminormalQuiver}, which means that
  $f_{\s\t}$ are simultaneous eigenvectors for $L_1,\dots,L_n$.
  Explicitly, in view of \cite[Proposition~3.7]{JamesMathas:cycJantzenSum},
  \begin{equation}\label{E:fstLk}
        L_kf_{\s\t} = [\con_k(\s)]f_{\s\t}\quad\text{and}\quad
        f_{\s\t}L_k = [\con_k(\t)]f_{\s\t},
        \qquad\text{for }1\le k\le n.
  \end{equation}
  These two formulas are equivalent since $f_{\s\t}^*=f_{\t\s}$. Using
  \autoref{E:fstLk} and the definitions it follows that
  $F_\u f_{\s\t}F_\v=\delta_{\s\u}\delta_{\t\v}f_{\s\t}$ and hence that
  there exist scalars $\gamma_\t\in\Kcal$ such that
  \begin{equation}\label{E:gamma}
  f_{\s\t}f_{\u\v} = \delta_{\t\u}\gamma_\t f_{\s\v},\qquad
     \text{for all $\s,\t\in\Std(\blam)$ and $\u,\v\in\Std(\bmu)$.}
  \end{equation}
  In particular, it follows that $F_\t=\frac1{\gamma_\t}f_{\t\t}$, for all
  $\t\in\Std(\Parts)$.

  Suppose that $\s,\t\v\in\Std(\blam)$ and that $\v=\t(r,r+1)$. Then by
  \cite[Proposition~3.11]{HuMathas:SeminormalQuiver},
  \[
  f_{\s\t}T_r = \begin{cases*}
    f_{\s\v} - \frac{1}{[\rho_r(\t)]}f_{\s\t},& if $\t\gdom\v$,\\
    \frac{[1+\rho_r(\t)][1-\rho_r(\t)]}{[\rho_r(\t)]^2}f_{\s\v}
         - \frac{1}{[\rho_r(\t)]}f_{\s\t},& if $\v\gdom\t$,\\
  \end{cases*}
  \]
  where $\rho_r(\t)=\con_r(\t)-\con_{r+1}(\t)$. Using this to compute
  $\gamma_\v f_{\v\v}=f_{\v\v}f_{\v\v}$, when $\s=\v$, shows that:

  \begin{Corollary}[{\cite[Corollary~3.10]{HuMathas:SeminormalQuiver}}]\label{C:GammaRecurrence}
    Suppose that $\t,\v\in\Std(\blam)$, $\t\gdom\v$ and that
    $\v=\t(r,r+1)$, for some $r$ with $1\le r<n$. Then
    \[
       \gamma_\v=\frac{[1+\rho_r(\t)][1-\rho_r(\t)]}{[\rho_r(\t)]^2}\gamma_\t.
    \]
  \end{Corollary}

  \subsection{Idempotents and characters}\label{S:idempotents}
    In view of \autoref{T:KLRIsomorphism}, over a field the cyclotomic
    Hecke algebra $\Hn$ has analogues of the KLR idempotents
    $e(\bi)\in\Rn$, for $\bi\in I^n$. Rather than working over a field,
    we need these idempotents over $\Ocal$. Fortunately, these
    idempotents are easy to describe using the seminormal basis. For
    $\bi\in I^n$ define
    \[
         f_\bi = \sum_{\t\in\Std_\bi(\Parts)} F_\t =
         \sum_{\t\in\Std_\bi(\Parts)}\frac1{\gamma_\t}f_{\t\t}.
    \]
    By definition, $f_\bi\in\HK$. In fact, inspired by results of
    Murphy, we have:

    \begin{Lemma}\label{L:Idempotents}
      Suppose that $\bi\in I^n$. Then $f_\bi=f^*_\bi$ is an idempotent in $\HO$.
      Moreover, if $\bj\in I^n$ then $f_\bi f_\bj=\delta_{\bi\bj}f_\bi$.
    \end{Lemma}

    \begin{proof}
      First observe that if $\t$ is a standard tableau then
      $F_\t^*=F_\t$ since, by definition, $L_k^*=L_k$ for $1\le k\le n$. Therefore,
      $f_\bi^*=f_\bi$, for all~$\bi\in I^n$.

      By \cite[Lemma~4.3]{HuMathas:SeminormalQuiver}, $f_\bi\in\HO$. In order to apply this
      lemma from~\cite{HuMathas:SeminormalQuiver} we need to first show
      that $(\Ocal, t)$ is an ``idempotent subring of $\Kcal$'' in the
      sense of \cite[Definition~4.1]{HuMathas:SeminormalQuiver}. In
      fact, there is nothing to do here because this is already
      established in \cite[Example~4.2(b)]{HuMathas:SeminormalQuiver}.

      Finally, $f_\bi f_\bj=\delta_{\bi\bj}f_\bi$ because
      $\set{F_\t|\t\in\Std(\Parts)}$ is a set of pairwise orthogonal idempotents.
    \end{proof}

    \begin{Remark}\label{R:idempotents}
      By definition, $f_\bi\ne0$ if and only if $\bi=\res(\t)$ for some
      standard tableau $\t\in\Std(\Parts)$. In contrast, it is not clear
      from \autoref{D:KLRAlgebra} when the KLR idempotent $e(\bi)$ is
      non-zero.  In fact, in view of \autoref{T:GradedCellular},
      $e(\bi)\ne0$ if and only if $f_\bi\ne0$.
    \end{Remark}

    Since $\Hn\cong\HO\otimes_\Ocal F$ we obtain pairwise orthogonal
    idempotents $f_\bi\otimes1$ in $\Hn$. We abuse notation and write
    $f_\bi$ for these idempotents in both $\HO$ and $\Hn$. The meaning
    should always be clear from context.

    As a first consequence of the existence of the idempotents $f_\bi$
    in~$\HO$ and $\Hn$, we can define an analogue of the
    graded characters of \autoref{D:qcharacter} for these algebras. If $M$ is an
    (ungraded) $\Hn$-module then $M=\bigoplus_\bi M_\bi$ as an
    $\Ocal$-module, when $M_\bi=Mf_\bi$.

    \begin{Definition}\label{D:character}
    Let $M$ be an $\Hn(F)$-module. the \textbf{character}
    of~$M$ is
    \[
            \ch M = \sum_{\bi\in I^n} (\dim M_\bi)\, \bi\in\Z[I^n].
    \]
    \end{Definition}

    Just as with the graded character $\Chq$, we can view $\ch$ as an
    exact functor from the category of $\Hn$-modules to the category of
    $\Ln'$-modules, where $\Ln'=\<L_1,\dots,L_n\>$.  Let $\Rep\Hn$ be
    the Grothendieck group of finitely generated $\Hn$-modules. Then we
    can view $\ch$ as the induced map
    \[\ch\map{\Rep\Hn}\Rep\Ln'\hookrightarrow\Z[I^n].\]
    Like the graded characters and \autoref{P:qCharacters}, we have:

  \begin{Proposition}
    \label{P:Characters}
    Let $F$ be a field. The character map
    $\ch\map{\Rep\Rn(F)}\Rep\Z[I^n]$ is injective.
  \end{Proposition}

  \begin{proof}
    This is part of the folklore for $\Hn$. It is proved in exactly the
    same way as \autoref{P:qCharacters} using the corresponding results
    for the affine Hecke algebra; see, for example,
    \cite[Theorem~3.3.1]{Klesh:book}.
  \end{proof}

  Finally, note that if $M$ is a graded $\Rn(F)$-module then, in view of
  \autoref{T:KLRIsomorphism}, $\underline{M}$ is naturally an
  $\Hn(F)$-module. The proof \autoref{T:KLRIsomorphism} (or, more
  correctly, \cite[Theorem~A]{HuMathas:SeminormalQuiver}), identifies
  the KLR idempotent $e(\bi)\in\Rn(F)$ with $f_\bi\in\Hn(F)$, for all
  $\bi\in I^n$. Therefore, $\ch\underline{M}=\Chq M\big|_{q=1}$. Hence,
  \autoref{P:Characters} implies \autoref{P:qCharacters}.

  \subsection{Gram determinants of Specht modules}

    In this section we compute the Gram determinant of the Specht module
    with respect to the Murphy basis. Gram determinants are only
    well-defined up to a sign that depends on the ordering of the rows
    and columns. In this and later sections we fix a total ordering of
    the tableaux in $\Std(\blam)$ that refines the dominance ordering
    and we use this ordering for both the rows and the columns of the
    Gram matrices. In fact, \autoref{C:GramDet} below implies that if
    the same ordering is used for the rows and columns of the Gram
    matrix then the Gram determinant is independent of this choice of
    ordering.

    Throughout this section fix an $\ell$-partition $\blam\in\Parts$.
    The \textbf{Gram matrix} of the Specht module $S^\blam_\Kcal$ is
    \[ G^\blam = \Bigl(\<m_\s,m_\t\>\Bigr)_{\s,\t\in\Std(\blam)}.\]
    For $\t\in\Std(\blam)$ define $f_\t=m_\t
    F_\t$.  Then $f_\t=m_\t+\sum_{\u\gdom\t}s_\u m_\u$ by
    \autoref{E:fstExpansion}, for some $s_\u\in\Kcal$. Therefore,
    $\set{f_\t|\t\in\Std(\blam)}$ is a ``seminormal'' basis of
    $S^\blam_\Kcal$.

    The following result is essentially a restatement of
    \cite[Proposition~3.19]{JamesMathas:cycJantzenSum}, although the
    proof below is considerably easier because it exploits the defining
    property of the inner product on a cell module.

    \begin{Lemma}
      Suppose that $\s,\t\in\Std(\blam)$, for $\blam\in\Parts$. Then
      $\<f_\s,f_\t\>=\delta_{\s\t}\gamma_\t$.
    \end{Lemma}

    \begin{proof}
      In view of \autoref{E:InnerProduct} and \autoref{E:gamma},
      $\<f_\s,f_\t\>f_\t=f_\s f_{\t\t}=\delta_{\s\t}\gamma_\s f_\s$.
      Hence, $\<f_\s,f_\t\>=\delta_{\s\t}\gamma_\t$ as claimed.
    \end{proof}

    The transition matrix between the Murphy and the
    seminormal bases of $S^\blam_\Kcal$ is unitriangular, so:

    \begin{Corollary}\label{C:GramDet}
      Suppose that $\blam\in\Parts$. Then
      $\det G^\blam = \det\Bigl(\<f_\s,f_\t\>\Bigr)
                    = \Prod_{\t\in\Std(\blam)}\gamma_\t$.
    \end{Corollary}

    The Gram determinant is quite a crude statistic. To approach the
    Jantzen sum formula we need to embellish $\det G^\blam$ by adding
    enough data so that it determines a formal character. To do this we
    need a new cellular basis of~$\Hn$.

    Let $\s,\t\in\Std(\blam)$ and set $\bi=\res(\s)$ an $\bj=\res(\t)$.
    Define $b_{\s\t} = f_\bi m_{\s\t} f_\bj$. We can view $b_{\s\t}$ as
    an element of~$\HO$ or of~$\Hn$. By \autoref{L:Idempotents},
    $b_{\s\t}=f_\bi b_{\s\t} f_\bj$, where $\bi=\res(\s)$ and
    $\bj=\res(\t)$.

    \begin{Proposition}\label{P:bbasis}
      The set $\set{b_{\s\t}|\s,\t\in\Std(\blam), \blam\in\Parts}$
      is a cellular basis of $\HO$.
    \end{Proposition}

    \begin{proof}
      By \autoref{E:msLk}, if $\u\in\Std_\bi(\Parts)$ then there exist
      scalars $r_{\t\u\v}\in\Ocal$ such that
      \[
        m_{\s\t}F_\u\equiv\delta_{\t\u}m_{\s\t}+\sum_{\v\gdom\t}r_{\t\u\v}m_{\s\v}
             \pmod{\Hlam}.
      \]
      Hence, $m_{\s\t}f_\bi$ is equal to $m_{\s\t}$ plus a linear
      combination of more dominant terms, so the transition matrix from
      the basis~$\set{m_{\s\t}}$ to the set $\set{b_{\s\t}}$ is unitriangular. In
      particular, $\set{b_{\s\t}}$ is a basis of $\HO$. It is now
      straightforward to check that $\set{b_{\s\t}}$ is a cellular basis
      of $\HO$. As the cellularity of this basis is not used in what
      follows we leave these details to the reader.
    \end{proof}

    \begin{Remark}
      The $b$-basis is compatible with the block decomposition of~$\Hn$.
      This basis first appeared in a more general context in
      \cite[Theorem~4.5]{Mathas:Seminormal}.  Instead of the $b$-basis
      we could use an analogue of the $\psi$-basis from
      \autoref{T:GradedCellular} for $\HO$, which is constructed in
      \cite{HuMathas:SeminormalQuiver}. There is no real gain in doing
      this, however, because the $\psi$-basis of $\HO$ takes
      considerably more effort to define.
    \end{Remark}

    The $b$-basis $\set{b_{\s\t}}$ of $\Hn$ gives a $b$-basis
    $\set{b_\t}$ of the Specht module $S^\blam_\Ocal$ and, by extension
    of scalars, a $b$-basis of $S^\blam_\Kcal$. The easiest way to see this
    is to notice that $b_{\tlam}=m_{\tlam}$, so that
    $b_{\tlam\t}=\sum_{\v\gedom\t}a_{\t\v}m_{\tlam\v}$ for some
    $a_{\t\v}\in\Ocal$ by the proof of \autoref{P:bbasis}. Hence, we can
    set $b_\t=\sum_\v a_{\u\v}m_\v\in S^\blam$, for $\t\in\Std(\blam)$.
    For each $\bi\in I^n$ define Gram matrices:
    \[
        G^\blam_\bi  = \Bigl(\<b_\s,b_\t\>\Bigr)_{\s,\t\in\Std_\bi(\blam)}.
    \]
    We adopt the convention that $\det G^\blam_\bi=1$ if
    $\Std_\bi(\blam)=\emptyset$. Then we have:

    \begin{Lemma}\label{L:GramiDet}
      Let $\blam\in\Parts$. Then
      $\det G^\blam = \prod_{\bi\in I^n}\det G^\blam_\bi$. Moreover, if
      $\bi\in I^n$ then
      \[\det G^\blam_\bi = \prod_{\t\in\Std_\bi(\blam)}\gamma_\t.\]
    \end{Lemma}

    \begin{proof}
      The transition matrix from the Murphy basis to the $b$-basis is
      unitriangular, so $\det G^\blam = \det\bigl(\<b_\s,b_\t\>\bigr)$.
      Next, observe that if $\s\in\Std_\bi(\blam)$ and
      $\t\in\Std_\bj(\blam)$ then
      \[
         \<b_\s,b_\t\> = \<b_\s f_\bi, b_\t f_\bj\>= \<b_\s, b_\t  f_\bj f_\bi^*\>
              = \delta_{\bi\bj}\<b_\s,b_\t\>,
      \]
      where the last inequality follows by \autoref{L:Idempotents}.
      Hence, $\det G^\blam = \prod_{\bi\in I^n}\det G^\blam_\bi$. Finally,
      \[\det G^\blam_\bi=\Bigl(\<f_\s,f_\t\>\Bigr)=\prod_{\t\in\Std_\bi(\blam)}\gamma_\t\]
      because, by going through the Murphy basis,  the transition matrix
      from the $b$-basis to the seminormal basis of~$S^\blam_\Kcal$ is
      unitriangular.
    \end{proof}

  \section{Jantzen characters}

  Jantzen filtrations of Specht modules underpin all of the results of
  this paper. Extending what we said in the introduction, Jantzen
  filtrations can be defined whenever we are given a module over a
  principal ideal domain that has a non-degenerate bilinear form.
  Therefore, to define the Jantzen filtrations of the Specht modules we
  need to work in the ungraded setting because the form on the graded
  Specht module $\Slam[\bmu]_\Q$ is degenerate unless
  $\Slam[\bmu]_\Q=\Dmu_\Q$, as is evident from
  \autoref{T:AdjustmemtMatrix}. On the other hand, we need to
  incorporate the idempotents $e(\bi)$ from \autoref{SS:CyclotomicKLR}
  into the Jantzen filtrations, which requires a careful choice of
  ground-ring. This section sets up the necessary machinery and then
  proves our main results.

  \subsection{Jantzen filtrations}
  We recall Jantzen's construction of filtrations of modules that come
  equipped with a non-degenerate bilinear form in the special case when
  the ground ring is $F[x]_{(x)}$. We then apply Jantzen's construction
  to the Specht modules of~$\Hn(F)$.

  As above, let $\Ocal=F[x]_{(x)}$ and let $\Kcal=F(x)$. Let
  $M_\Ocal$ be a free $\Ocal$-module of finite rank and let
  $M_\Kcal=M_\Ocal\otimes_\Ocal\Kcal$. Suppose that there exists a
  non-degenerate bilinear form
  \[
            \<\ ,\ \>\map{M_\Ocal\times M_\Ocal}\Ocal.
  \]
  By definition, $\Ocal$ is a discrete valuation ring with maximal ideal
  $x\Ocal$. Let $\nu_x$ be the associated valuation. Explicitly, if
  $0\ne a\in\Ocal$ then $\nu_x(a) = \max\set{k\ge0|a\in x^k\Ocal}$. In
  particular, $\nu_x(x+\xi)=0$.

  The \textbf{Jantzen filtration} of $M_\Ocal$ is the filtration
  \[
     M_\Ocal = J_0(M_\Ocal)\supseteq J_1(M_\Ocal)\supseteq J_2(M_\Ocal)\supseteq\dots
  \]
  where $J_k(M_\Ocal) = \set{m\in M_\Ocal|\nu_x(m)\ge k}$.

  Let $M_F = M_\Ocal/xM_\Ocal$. Then the \textbf{Jantzen filtration}
  of $M_F$ is the filtration
  \[
     M_F = J_0(M_F)\supseteq J_1(M_F)\supseteq J_2(M_F)\supseteq\dots
  \]
  where $J_k(M_F) = \bigl(J_k(M_\Ocal)+xM_\Ocal\bigr)/xM_\Ocal
         \cong J_k(M_\Ocal)/\bigl(J_k(M_\Ocal)\cap xM_\Ocal\bigr)$, for $k\ge0$.
 Since $M_F$ is finite dimensional, $J_k(M_F)=0$ for $k\gg0$. The
 following easy observation of Jantzen's is the key to computing the
 Jantzen sum formula.

 \begin{Lemma}[Jantzen]\label{L:Jantzen}
      Suppose that $M_\Ocal$ has a non-degenerate bilinear form. Then
      \[
          \sum_{k>0}\dim J_k(M_F) = \nu_x(G_{M_\Ocal}).
      \]
  \end{Lemma}

  The proof is a straightforward calculation using the Smith normal
  form of~$G_M$, which exists because~$\Ocal$ is a principal ideal
  domain by assumption. For a proof that uses very similar language and
  notation to what we are using see \cite[Lemma~5.30]{Mathas:ULect}.

  The cellular basis of \autoref{T:MurphyBasis} gives us
  a non-degenerate bilinear form on $\SOlam$, so we can apply
  Jantzen's constructions to the Specht modules $\SOlam$, for
  $\blam\in\Parts$. Rather than looking at the Jantzen filtration
  directly we consider their direct sum.

  \begin{Definition}
    Suppose that $\blam\in\Parts$ and let
    $\displaystyle J^\blam_F = \oplus_{k>0} J_k(S^\blam_F)$.
    The \textbf{Jantzen character} of~$S^\blam_F$ is
    \[  \ch J^\blam_F = \Sum_{k>0}\ch J_k(S^\blam_F).\]
  \end{Definition}

  The main results of this paper follow from an explicit description of the
  Jantzen character $\ch J^\blam_F$. The next result is our first step
  towards describing this character.

  \begin{Proposition}\label{P:JantzenDetCharacter}
    Let $\blam\in\Parts$. Then
    $\ch J^\blam_F=\Sum_{\bi\in I^n}\nu_x(\det G^\blam_\bi)\bi$.
  \end{Proposition}

  \begin{proof}
    As an $\Ocal$-module, $S^\blam_\Ocal = \bigoplus_{\bi\in I^n}
    S^\blam_\Ocal f_\bi$.  By construction,
    $\set{b_\s|\s\in\Std_\bi(\blam)}$ is a basis of $S^\blam_\Ocal f_\bi$
    and~$G^\blam_\bi$ is the Gram matrix of the bilinear form
    on~$S^\blam_\Ocal$ restricted to~$S^\blam_\Ocal f_\bi$. Moreover, in
    view of \autoref{L:GramiDet}, if $\bi\ne\bj\in I^n$ then the
    summands $S^\blam_\Ocal f_\bi$ and $S^\blam_\Ocal f_\bj$ are
    orthogonal with respect to the inner product on~$S^\blam_\Ocal$.
    Therefore, $J_k(S^\blam_\Ocal)f_\bi=J_k(S^\blam_\Ocal f_\bi)$ and
    $J_k(S^\blam_F)f_\bi=J_k(S^\blam_Ff_\bi)$, for all $\bi\in I^n$. (As in
    \autoref{S:idempotents}, we are abusing notation and identifying
    $f_\bi\in\HO$ and $f_\bi\otimes1_F\in\Hn$.) Hence,
    \[
      \ch J^\blam_F = \sum_{k>0}\ch J_k(S^\blam_F)
           =\sum_{\bi\in I^n}\sum_{k>0}\dim J_k(S^\blam_F)f_\bi \cdot \bi
           =\sum_{\bi\in I^n}\sum_{k>0}\dim J_k(S^\blam_Ff_\bi)\cdot \bi
           =\sum_{\bi\in I^n}\nu_x(\det G^\blam_\bi) \bi,
    \]
    where the last equality follows by applying \autoref{L:Jantzen} to
    $S^\blam_\Ocal f_\bi$ and $S^\blam_F f_\bi$. Note that \autoref{L:Jantzen}
    applies because the bilinear form on~$S^\blam_\Kcal f_\bi$ is
    non-degenerate by \autoref{L:GramiDet}. This completes the proof.
  \end{proof}

  \autoref{L:GramiDet} gives an explicit formula for $\det G^\blam_\bi$
  in terms of the $\gamma$-coefficients. The next step is to calculate
  $\nu_x(\gamma_\t)$, for $\t\in\Std(\blam)$. First, a small
  interlude on cyclotomic polynomials.

  \subsection{Cyclotomic polynomials}
    This section contains some brief reminders on cyclotomic polynomials
    that we will need to the compute the Jantzen characters. All of the
    facts that we quote are standard and can be found, for example, in
    \cite[VI \S3]{Lang:Algebra}.

    Let $f\ge1$. The \textbf{$f$th cyclotomic polynomial} in $\C[x]$ is
    \[
            \Phi_f(x) = \prod_{\substack{1\le d\le
            f\\\gcd(d,f)=1}}\Bigl(x-\exp(\tfrac{2\pi id}{f})\Bigr)
    \]
    Since $\Phi_f(x)$ is fixed by complex conjugation it follows that
    $\Phi_f(x)\in\R[x]$.  In fact, $\Phi_f(x)\in\Z[x]$. Hence, by
    base change, we consider $\Phi_f(x)$ as an element of $F[x]$.

    For this paper, a trivial but important observation is that
    $x^f-1=\prod_{d\mid f}\Phi_d(x)$. Therefore,
    \begin{equation}\label{E:QIntFactorisation}
      [f]_x= \frac{x^f-1}{x-1} = \prod_{\substack{1<d\le f\\d\mid f}}\Phi_d(x).
    \end{equation}
    In particular, if $p$ is a prime integer then
    $\Phi_p(x)=[p]_x=1+x+\dots+x^{p-1}$. Moreover, if $f\in\Z$ then
    \begin{equation}\label{E:CycPolyFactorise}
            \Phi_{fp}(x) = \begin{cases*}
              \dfrac{\Phi_f(x^p)}{\Phi_f(x)},& if $p\nmid f$,\\
              \Phi_f(x^p),& if $p\mid f$.\\
            \end{cases*}
    \end{equation}
    This recurrence relation determines $\Phi_d(x)$ uniquely as a
    polynomial in $F[x]$.

    Recall that $z=x+\xi$, that $e>0$ is minimal such that $[e]_\xi=0$,
    and that $F$ is a field of characteristic $p\ge0$.  We can now state
    and prove the main fact that we need about cyclotomic polynomials.

    \begin{Lemma}\label{L:nuPhi}
      Suppose that $F$ is an algebraically closed field of
      characteristic $p\ge0$ and that $f\ge1$. Then
      \[
        \nu_x\bigl(\Phi_f(z)\bigr) = \begin{cases*}
          1,& if $f=e$,\\
          (p-1)p^{r-1}, & if $p>0$ and $f=ep^r$ for some $k>0$,\\
          0,& otherwise.
        \end{cases*}
      \]
    \end{Lemma}

    \begin{proof}
        First suppose that $p=0$. Since $F$ is algebraically closed,
        it contains a complete set
        $\set{\omega_1,\dots,\omega_{\phi_f}}$ of primitive $f$th roots
        of unity in~$F$ and $\Phi_f(x)=\prod_{k=1}^{\phi_f}(x-\omega_k)$.
        The constant term of $\Phi_f(z)$ is~$\Phi_f(\xi)$, which is
        non-zero if and only if $e\ne f$ since~$\xi$ is a primitive $e$th
        root of unity. In particular, $\nu_x(\Phi_f(z))=0$ if~$f\ne e$.
        On the other hand, if $f=e$ then $\xi=\omega_s$, for some~$s$,
        so the coefficient of $x$ in $\Phi_e(z)$ is $\prod_{k\ne
        s}(\xi-\omega_k)\ne0$, so $\nu_x(\Phi_e(z))=1$.  Therefore,
        if~$p=0$ then $\nu_x(\Phi_f(z))=\delta_{ef}$.

        Now suppose that $p>0$. Write $f=f'p^r$, for integers $f'$ and
        $r\ge0$ such that $\gcd(f',p)=1$. If $r=0$, so that $f=f'$, then
        the constant term of $\Phi_f(z)$ is $\Phi_{f}(\xi)$ which is
        non-zero if and only if $e=f=f'$. Hence, by the argument of the
        last paragraph, $\nu_x\bigl(\Phi_f(z)\bigr)=\delta_{ef}$ when
        $r=0$. Finally, suppose that $r>0$. By
        \autoref{E:CycPolyFactorise} and the fact that we are in
        characteristic~$p>0$,
        \[
            \Phi_f(z) = \bigl(\Phi_{f'}(z)\bigr)^{(p-1)p^{r-1}}.
        \]
        Therefore,
        $\nu_x\bigl(\Phi_f(z)\bigl)=(p-1)p^{r-1}\nu_x\bigl(\Phi_{f'}(z)\bigr)
                =(p-1)p^{r-1}\delta_{ef'}$,
        completing the proof.
    \end{proof}

  \subsection{Computing Jantzen characters}

    In \cite[Theorem~3.13]{HuMathas:SeminormalQuiver}, Hu and the author
    gave a closed formula for $\det G^\blam$ using the KLR degree
    function \autoref{E:Tableaux} on standard tableaux
    \[ \deg_f\map{\Std(\blam)}\Z; \t\mapsto\deg_f(\t).\]
    We need to generalise this result to give a formula for $\det
    G^\blam_\bi$, for $\bi\in I^n$. First, for $\blam\in\Parts$ define
    \[
      [\blam]_z^!=\prod_{l=1}^\ell\prod_{r\ge0}[\lambda^{(l)}_r]_z^!,
    \]
    where for $k>0$ the \textbf{$t$-quantum factorial} of $k$ is
    $[k]_z^!=[k]_z[k-1]_z\dots[1]_z$. For convenience, set $[0]^!_z=1$.

    An analogue of the next result is implicit in the proof
    of~\cite[Theorem~3.13]{HuMathas:SeminormalQuiver}.

    \begin{Lemma}\label{L:gammaDegree}
      Suppose that $\t\in\Std(\blam)$, for $\blam\in\Parts$. Then there
      exists an integer $g_\t\in\N$ such that
      \[
         \gamma_\t = z^{g_\t}\prod_{f>1}\Phi_f(x)^{\deg_f(\t)}.
      \]
    \end{Lemma}

    \begin{proof}
       The proof is by induction on the dominance order on~$\Std(\blam)$.
       By \cite[Proposition~3.11]{HuMathas:SeminormalQuiver}, which is a
       straightforward calculation using \autoref{E:msLk},
       \begin{equation}\label{E:gammatlam}
            \gamma_{\tlam} = [\blam]_z^!\prod_{1\le l<m\le\ell}\prod_{(l,r,c)\in[\blam]}
                            \bigl([\kappa_l+c-r]_z-[\kappa_m]_z\bigr)
              =\prod_{(l,r,c)\in[\blam]} [c]_z\prod_{l<m\le\ell}
              z^{\kappa_m}[\kappa_l+c-r-\kappa_m]_z.
       \end{equation}
       Let $A=(l,r,c)\in[\blam]$ be a node in $\blam$ and set
       $r=\tlam(A)$ and $\bmu=\Shape(\tlam_{\downarrow r})$.
       Recalling \autoref{E:QIntFactorisation}, and the definition of the
       integer $d_{A,f}(\bmu)$ from \autoref{E:TableauDegree}, the
       contribution of the node $A$ to~$\gamma_{\tlam}$ is
       \[
            [c]_z\prod_{l<m\le\ell}z^{\kappa_m}[\kappa_l+c-r-\kappa_m]_z\\
                = z^{\kappa_{l+1}+\dots+\kappa_\ell}\prod_{f>1}
                    \Phi_f(z)^{d_{A,f}(\bmu)}.
       \]
       Hence, the lemma holds when $\t=\tlam$.

       Now suppose that $\t\ne\tlam$. Then there exists a tableau
       $\s\in\Std(\blam)$ such that $\s\gdom\t=\s(r,r+1)$, where $1\le r<n$.
       By \autoref{C:GammaRecurrence}, and induction,
       \[
       \gamma_\t = \frac{[1+\rho_r(\t)]_z[1-\rho_r(\t)]_z}{[\rho_r(\t)]_z^2}
                        \gamma_\s
                 = \frac{[1+\rho_r(\t)]_z[1-\rho_r(\t)]_z}{[\rho_r(\t)]_z^2}
                      z^{g_\s}\prod_{f>1}\Phi_f(x)^{\deg_f(\s)}.
       \]
       Let $\bi=\res(\t)$. In the graded Specht module $\Slam$,
       $\psi_\t=\psi_\s\psi_r$ with the degrees adding. Therefore,
       \[
       \deg_f(\t) = \deg_f(\s) - (\alpha_{i_r}, \alpha_{i_{r+1}})
           = \deg_f(\s) + \begin{cases*}
                1,  & if $1\pm\rho_r(\t)\equiv0\pmod f$,\\
                -2, & if $\rho_r(\t)\equiv0\pmod f$,\\
                0,  & otherwise.
           \end{cases*}
       \]
       Hence, by induction, $\gamma_\t$ can be written in the required form.
       This completes the proof.
    \end{proof}

    We are one definition away from our first description of the Jantzen
    characters.  First, recall that $e$ is the quantum characteristic
    of~$\xi\in F$, so that $e$ is the minimal positive integer such that
    $[e]_\xi=0$. That is, either $\xi=1$ and $e=\Char F$, or $\xi$ is a
    primitive $e$th root of unity in $F$. Consequently, either $e=p$ or
    $\gcd(e,p)=1$.

    \begin{Definition}\label{D:epDegree}
      Suppose that $\blam\in\Parts$ and $\bi\in I^n$. Define the
      \textbf{$(e,\bi)$-degree} of $\blam$ to be the integer
      \[
          \deg_{e,\bi}(\blam) = \sum_{\t\in\Std_\bi(\blam)}\deg_e(\t).
      \]
      Let $p$ be the characteristic of the field $F$. Define the
      \textbf{$(e, p,\bi)$-degree} of $\blam$ to be
      \[
      \pdeg_{e,\bi}(\blam)=
          \deg_{e,\bi}(\bi)+\sum_{r\ge1}(p-1)p^{r-1}\deg_{ep^r,\bi}(\blam).
      \]
    \end{Definition}

    Recall from \autoref{E:TableauDegree} that $\deg_0(\t)=0$, for all
    $\t\in\Std(\Parts)$. Therefore,
    $\pdeg_{e,0,\bi}(\blam)=\deg_{e,\bi}(\blam)$.

    The next result gives a complete description of the Jantzen
    characters. This formula has the advantage of being easy to compute in
    examples but it is not particularly useful in practice!

    \begin{Proposition}\label{P:JantzenDegCharacter}
      Suppose that $\blam\in\Parts$. Then
        $\ch J^\blam_F = \Sum_{\bi\in I^n} \pdeg_{e,p,\bi}(\blam)\,\bi$.
    \end{Proposition}

    \begin{proof}
      Since $\Hn(F)$ is a cellular algebra by \autoref{T:MurphyBasis},
      every field is a splitting field for so without loss of generality
      we can and do assume that $F$ is algebraically closed.
      By \autoref{P:JantzenDetCharacter} and \autoref{L:gammaDegree},
      \begin{align*}
            \ch J^\blam_F &= \sum_{\bi\in I^n}\nu_x(\det G^\blam_\bi)\bi
                 && \text{by \autoref{P:JantzenDetCharacter}},\\
              & = \sum_{\bi\in I^n}\nu_x\Bigl(\prod_{f>1}
                \Phi_f(z)^{\deg_{f,\bi}(\blam)}\Bigr)\bi
                 && \text{by \autoref{L:gammaDegree}},\\
              & = \sum_{f>1} \Bigl(\sum_{\bi\in I^n}
               \deg_{f,\bi}(\blam)\nu_x\bigl(\Phi_f(z)\bigr)\Bigr)\bi\\
                & = \sum_{\bi\in I^n}\pdeg_{e,p,\bi}(\blam)\,\bi,
      \end{align*}
      where the last equality follows by \autoref{L:nuPhi} and \autoref{D:epDegree}.
    \end{proof}

    As the coefficient of $\bi$ in the Jantzen character is a
    non-negative integer we obtain the following combinatorial
    statements, refining similar results from \cite[\S3.3]{HuMathas:SeminormalQuiver}.

    \begin{Corollary}
      Suppose that $\blam\in\Parts$ and $\bi\in I^n$. Then
      $\deg_{e,\bi}(\blam)\ge0$ and $\pdeg_{e,p,\bi}(\blam)\ge0$.
    \end{Corollary}

  \subsection{The positive Jantzen sum formula}\label{SS:PositiveSum}

  Building on \autoref{P:JantzenDegCharacter}, we are now ready to prove
  our main results, which describe the Jantzen characters as explicit linear
  combinations of the characters of $\Hn(F)$-modules.

  The next definition extends \autoref{C:Ungrading}.

  \begin{Definition}\label{D:Emu}
    Suppose that $F$ is a field.  Set
    $E^\bmu_F=\underline{\mathbb{E}}^\bmu_F$, for $\bmu\in\Klesh$.
  \end{Definition}

  Using \autoref{T:KLRIsomorphism}, we view~$E^\bmu_F$ as an
  $\Hn(F)$-module. The ``traditional'' way to define a module
  like~$E^\bmu_F$ is to construct a ``decomposition map'' by first
  choosing a modular system.  We have sidestepped these additional
  complications by using the KLR algebra $\Rn(\Z)$ to define~$\Emu_\Z$
  over~$\Z$. In effect, we used the triple $(\Q,\Z, F)$ as a ``modular
  system''. We return to this theme in \autoref{P:fAdjustment} below.

  If $f(q)\in\Acal=\Z[q,q^{-1}]$ is a Laurent polynomial let $f'(1)$ be the
  derivative of $f(q)$ evaluated at $q=1$. Define a linear map
  $\partial\map{\AIq}{\Z[I^n]}$ by
  \[
    \partial\Bigl(\sum_{\bi\in I^n}f_\bi(q)\bi\Bigr)
            = \sum_{\bi\in I^n}f_\bi'(1)\,\bi.
  \]
  By \autoref{R:idempotents},  we can restrict $\partial$ to a linear
  map of Grothendieck groups $\partial\map{\Rep\Ln}{\Rep\Ln'}$.

  For $\blam\in\Parts$ and $\bmu\in\Klesh$, recall from
  \autoref{E:GradedDecomp} that
  $d^{\C,e}_{\blam\bmu}(q)=[\Slam_\C:\Dmu_\C]_q\in\N[q,q^{-1}]$ is a
  characteristic zero graded decomposition number of $\Rn(\C)$.

  \begin{Lemma}\label{L:positivity}
    Suppose that $f=ep^r$, for $r\ge0$, and let $\blam\in\Parts$ and
    $\bmu\in\Kleshf$. Then $(d^{\C,f}_{\blam\bmu})'(1)\ge0$ with
    equality if and only if $\blam=\bmu$.
  \end{Lemma}

  \begin{proof}
    This is immediate because
    $d^{\C,e}_{\blam\bmu}(q)\in\delta_{\blam\bmu}+q\N[q]$ by
    \cite[Corollary~5.15]{BK:GradedDecomp}. (However, the proof of this
    fact from \cite{BK:GradedDecomp} is highly non-trivial.)
  \end{proof}

  According to the lemma, the coefficients in the character formula
  below are non-negative integers.

  \begin{Theorem}\label{T:PositiveJantzen}
    Suppose that $F$ is a field of characteristic zero and let $\blam\in\Parts$. Then
    \[
        \ch J^\blam_F = \Sum_{\bmu\in\Klesh} (d^{\C,e}_{\blam\bmu})'(1)\ch E^\bmu_F.
    \]
  \end{Theorem}

  \begin{proof}
    Let $\bmu\in\Klesh$ and observe that $\partial\Chq\Emu_F=0$ because
    $\overline{\Chq\Emu_F}=\Chq\Emu_F$ by \autoref{C:BarInvariance}. Therefore,
    \[
        \partial\Bigl(d^{\C,e}_{\blam\bmu}(q)\Chq\Emu_F\Bigr)
              =(d^{\C,e}_{\blam\bmu})'(1)\ch E^\bmu_F
    \]
    by the chain rule. Next, notice that \autoref{C:SlamEmuCharacters}
    implies that
    \[
          \sum_{\bi\in I^n}\sum_{\t\in\Std_\bi(\blam)}q^{\deg_e(\t)}\bi
             =\Chq \Slam_F = \sum_{\bmu\in\Klesh}d^{\C,e}_{\blam\bmu}(q)\Chq\Emu_F.
    \]
    Applying $\partial$ to both sides of this equation, using the
    previous remark, shows that
    \[
          \sum_{\bi\in I^n}\sum_{\t\in\Std_\bi(\blam)}\deg_e(\t)\,\bi
             = \sum_{\bmu\in\Klesh}d'_{\blam\bmu}(1)\ch E^\bmu_F.
    \]
    Since $\pdeg_{e,0,\bi}(\blam)=\deg_{e,\bi}(\bi)=\sum_{\t\in\Std_\bi(\blam)}\deg_e(\t)$,
    for $\bi\in I^n$, comparing the left hand side of the last displayed
    equation with \autoref{P:JantzenDetCharacter} completes the proof.
  \end{proof}

  The characters of the graded Specht modules are crucial to the proof
  of \autoref{T:PositiveJantzen}, however, the proof is purely
  combinatorial and it does not use a representation theoretic
  connection between the Jantzen filtration of $S^\blam_F$ and the
  graded $\Rn(F)$-modules $\Emu_F$. It would be interesting to find a
  direct connection between the modules $J_k(S^\blam_F)$ and $\Emu_F$.

  The \textit{Jantzen sum formula} is usually stated in the Grothendieck
  group $\Rep\Hn$. In view of \autoref{P:Characters},
  \autoref{T:PositiveJantzen} is equivalent to the following more
  ``traditional'' statement:

  \begin{Corollary}
    Suppose that $p=0$ and let $\blam\in\Parts$. Then
    \[
      [J^\blam_F] = \sum_{k>1} [J_k(S^\blam_F)] = \sum_{\bmu\in\Klesh}
                (d^{\C,e}_{\blam\bmu})'(1)[E^\bmu_F].
    \]
  \end{Corollary}

  \begin{Remark}
    In characteristic zero, Ryom-Hansen \cite[Theorem
    1]{RyomHansen:Schaper} proved an analogue of this result for
    $\Lambda=\Lambda_0$ by assuming a conjecture of
    Rouquier's~\cite[\S9]{LLT} that, in the post KLR world, says that
    the Jantzen filtration of~$\Slam_\C$ coincides with its
    \textit{grading filtration}.
    Yvonne~\cite[Theorem~2.11]{Yvonne:Conjecture} extended Ryom-Hanson's
    work to arbitrary $\Lambda\in P_e^+$.  \autoref{T:PositiveJantzen}
    shows that the description of the Jantzen sum formula in terms of
    the integers $(d^{\C,e}_{\blam\bmu})'(1)$ is independent of
    Rouquier's conjecture.
  \end{Remark}

  We next compute the Jantzen characters for fields of positive
  characteristic $p>0$. As \autoref{P:JantzenDegCharacter} and
  \autoref{D:epDegree} suggest, this involves looking at Hecke algebras
  at $ep^r$th roots of unity for $r\ge0$. We first need to set up the
  combinatorial and representation theoretic machinery to do this. When
  $r=0$ the tool that we need is provided by
  \autoref{T:AdjustmemtMatrix} but when $r>0$ we need to work harder.

  If $f=ep^r$, for $r\ge0$, let $\zeta_f=\exp(\tfrac{2\pi}{f})\in\C$.
  Then $\zeta_f$ is a primitive $f$th root of unity in~$\C$. Let
  $I_f=\Z/f\Z$. (In particular, $I=I_e$.) Since $e$ divides $f$, there is
  a well-defined surjective map $I_f\longrightarrow I=I_e; i\mapsto i_f$
  given by ``reducing modulo~$e$''. More explicitly, this map sends
  $a+f\Z$ to $a+e\Z$, for $a\in\Z$.  Extending this notation, if $\bi\in I_f^n$ let
  $\bi_f$ be the corresponding sequence in $I^n$. This map induces an
  abelian group homomorphism
  \[
        \Lf\map{\Z[I^n_f]}{\Z[I^n]};
            \sum_{\bi\in I^n_f}c_\bi\bi\mapsto
               \sum_{\bi\in I_f^n}c_\bi\bi_f,
               \qquad\text{for }c_\bi\in\Z.
  \]

  Recall from \autoref{SS:CyclotomicKLR} that $P^+_f$ and $Q^+_f$ are
  the positive weight lattice and positive root lattice the attached to
  the quiver $\Gamma_f$, for $f\ge2$. Abusing notation, let
  $\Lambda_i$ and $\alpha_i$ be fundamental weights and simple roots in
  the corresponding weight lattices for any $f\ge 2$, with the meaning being clear from
  the choice of index set~$I_f$. In \autoref{SS:CyclotomicKLR} we
  fixed a dominant weight $\Lambda=\Lambda_e\in P^+=P^+_e$. For $f=ep^r$ fix a dominant
  weight $\Lambda_f\in P^+_f$ such that
  \begin{equation}\label{E:Lambdaf}
    (\Lambda,\alpha_i) = \sum_{\substack{j\in I_f\\j_f=i}}(\Lambda_f,\alpha_j),
    \qquad\text{for all }i\in I.
  \end{equation}
  These conditions determine the dominant weight $\Lambda_f$ uniquely if
  and only if $r=0$. If $r>0$ then the results below are independent of
  the choice of~$\Lambda_f$.

  Let $\Hnf(\C)$ be the cyclotomic Hecke algebra over~$\C$ with Hecke
  parameter $\zeta_f$ and dominant weight~$\Lambda_f$. For
  $\blam\in\Parts$ let $\Slam_{\C,f}$ be the graded Specht module for
  $\Rnf(\C)$ indexed by $\blam$. Similarly, given an $\ell$-partition
  $\bmu\in\Kleshf=\Kleshf[f,n]$ let $\Emu_{\C,f}=\Dmu_{\C,f}$ be the
  (almost) simple $\Hnf(\C)$-module defined by \autoref{D:Simples} and let
  $S^\blam_{\C,f}$ and $D^\blam_{\C,f}$ be the corresponding Specht and
  simple modules for $\Hnf(\C)$. By \autoref{E:GradedDecomp}, we have
  polynomials
  \[
      d^{\C,f}_{\blam\bmu}(q) = [\Slam_{\C,f}: \Emu_{\C,f}]_q,\qquad
            \text{for $\blam\in\Parts$ and $\bmu\in\Kleshf$}.
  \]
  Armed with this notation we can state the result that we need. This
  result is partly motivated by the discussion after
  \cite[Conjecture~6.37]{Mathas:ULect}.

  \begin{Proposition}\label{P:fAdjustment}
      Suppose that $f=ep^r$, where $r\ge0$, and let $K$ be a field of
      characteristic zero that contains~$\zeta_f$. Then there is a
      unique homomorphism $\Df\map{\Rep\Hnf(K)}\Rep\Hn(F)$ such that the
      following diagram commutes:
      \[
          \begin{tikzpicture}[>=stealth,->,shorten >=2pt,looseness=.5,auto]
            \matrix (M)[matrix of math nodes,row sep=1cm,column sep=16mm]{
              \Rep\Hnf(K)& \Rep\Hn(F)\\
              \Z[I^n_f] & \Z[I^n]\\
             };
             \draw[-Latex] (M-1-1)--node[above]{$\Df$}(M-1-2);
             \draw[-Latex] (M-2-1)--node[below]{$\Lf$}(M-2-2);
             \draw[-Latex] (M-1-1)--node[left]{$\ch$}(M-2-1);
             \draw[-Latex] (M-1-2)--node[right]{$\ch$}(M-2-2);
          \end{tikzpicture}
      \]
      Moreover, $\Df[S^\blam_{K,f}]=[S^\blam_F]$, for all $\blam\in\Parts$.
  \end{Proposition}

  \begin{proof}
    The two character maps are injective by \autoref{P:Characters}, so
    once we show that $\Df$ exists it is automatically unique.

    The cyclotomic Hecke algebras are cellular, so every field is a
    splitting field. Hence, without loss of generality, we can assume
    that $K=\Q[\zeta_f]$ and $F=\mathbb{F}_p[\xi]$.  Let
    $\Zcal_f=\Z[\zeta_f]$. Then $\Hnf(\Zcal_f)$ is a
    $\Zcal_f$-subalgebra of $\Hnf(K)$ and
    $\Hnf(K)\cong\Hnf(\Zcal_f)\otimes_{\Zcal_f}K$. Since $e$ divides $f$
    there is a unique surjective ring homomorphism
    $\pi_f\colon\Zcal_f\to F$ such that $\pi_f(\zeta_f)=\xi$.
    Hence, we can consider $F$ as a
    $\Zcal_f$-algebra. We claim that, as $F$-algebras,
    \[
       \Hn(F)\cong\Hnf(\Zcal_f)\otimes_{\Zcal_f}F.
    \]
    Since the cyclotomic Hecke algebras are free over any ring by
    \autoref{T:MurphyBasis}, to prove this it is enough to check that
    tensoring with $\Fp[\xi]$ respects the relations. The only relation
    that is not obviously preserved is the cyclotomic relation, which
    can be written in the form
    $\prod_{i\in I_f}(L_1-[i]_{\zeta_f})^{(\Lambda_f,\alpha_i)}=0$ by
    \autoref{E:CyclotomicRelation}. This relation is preserved by
    tensoring with $\Zcal_f$ by \autoref{E:Lambdaf} because
    $[i_f]_\xi=\pi_f([i]_{\zeta_f})$, for all $i\in I_f$.

    Set $\mathfrak{p}_f=\ker\pi_f$, a prime ideal in~$\Zcal_f$, and let
    $\Ocal_f=(\Zcal_f)_{\mathfrak{p}_f}$ be the localisation of
    $\Zcal_f$ at $\mathfrak{p}_f$.  Then $\Ocal_f$ is a discrete
    valuation ring with maximal ideal $\mathfrak{p}\Ocal_f$, quotient
    field $K$ and residue field $F\cong\Ocal_f/\mathfrak{p}\Ocal_f$.
    Hence, the triple $(K,\Ocal_f,F)$ is a $p$-modular system in the
    sense of \cite[\S16A]{CurtisReiner:VolII}. Therefore, there is a
    well-defined decomposition map $\Df\map{\Rep\Hnf(K)}\Rep\Hn(F)$ that
    is independent of the choice of~$\Ocal_f$ by
    \cite[Proposition~16.17]{CurtisReiner:VolII}. To describe $\Df$
    explicitly, if~$V$ is an $\Hnf(K)$-module then $\Df(V)$ is defined
    by choosing a full $\Ocal_f$-lattice~$V_{\Ocal_f}$ in~$V$ and then
    setting $\Df(V)=V_{\Ocal_f}\otimes_{\Ocal_f}F$. Consequently,
    because the cellular basis~$\set{b_{\s\t}}$ of $\Hnf(K)$ from
    \autoref{P:bbasis} is defined over any ring, if~$\blam\in\Parts$
    then $\Df(S^\blam_K)=S^\blam_F$.

    It remains to show that the diagram in the proposition commutes.
    Since $\Lf(\ch S^\blam_K)=\ch(S^\blam_F)$, we have commutativity on
    the Specht modules. The decomposition matrix $\mathsf{D}_{K,f}$ is
    unitriangular, so $\set{S^\blam_F|\blam\in\Kleshf}$ is a basis for
    $\Rep\Hnf(K)$. Hence, this completes the proof.
  \end{proof}

  \begin{Definition}\label{D:Efbnu}
    Suppose that $f=ep^r$, for $r\ge0$, and let $\bnu\in\Kleshf$. Set
    $E^\bnu_{f,e} = \Df\bigl(E^\bnu_{\C,f}\bigr)$.
  \end{Definition}

  Strictly speaking, $E^\bnu_{f,e}$ is the image of a module in the
  Grothendieck group of $\Hn(F)$ but we abuse notation refer to it as
  the modular reduction of the almost simple module~$E^\bmu_\C$
  for~$\Hnf(\C)$. (By the proof of \autoref{P:fAdjustment}, the module
  $E^\bnu_{f,e}$ is well-defined only up to a choice of lattice
  whereas~$[E^\bnu_{f,e}]$ is independent of this choice.) In the
  special case when $e=f$ or, equivalently, $r=0$, comparing characters
  and applying \autoref{P:fAdjustment} yields the following:

  \begin{Corollary}\label{C:Emu}
    Let $\blam\in\Parts$ and $\bmu\in\Klesh$. Then
    $[S^\blam_F]=\Df(S^\blam_{\C,e})$ and $[E^\bmu_F]=[E^\bmu_{F,e}]$.
  \end{Corollary}

  This gives the following generalisation of the (ungraded analogue of)
  \autoref{T:AdjustmemtMatrix}. Define the \textbf{$(f,e)$-adjustment
  matrix} $\mathsf{A}_{f,e}=\bigl(a^{f,e}_{\bnu\bmu}\bigr)$, for
  $\bnu\in\Kleshf$ and $\bmu\in\Klesh$, by
  \[
     a^{f,e}_{\bnu\bmu} = [ E^\bnu_{f,e}: D^\bmu_F].
  \]
  \autoref{P:fAdjustment} gives the following generalisation of the
  ungraded analogue of \autoref{E:Adjustment} at $q=1$.

  \begin{Corollary}\label{C:eAdjustment}
    Suppose that $f=ep^r$, for $r\ge0$. Then
    $\mathsf{D}_{F,e}=\mathsf{D}_{\C,f}\mathsf{A}_{f,e}$.
  \end{Corollary}

  \begin{proof}
    If $\blam\in\Parts$ then
    $[S^\blam_{\C,f}]=\sum_\bnu d^{\C,f}_{\blam\bnu}[E^\bnu_{\C,f}]$
    in $\Rep\Hnf(\C)$, where the sum is over $\bnu\in\Kleshf$. Applying the linear map~$\Df$,
    \begin{align*}
          [S^\blam_F]
              =\sum_{\bnu\in\Kleshf}d^{\C,f}_{\blam\bnu}[E^\bnu_{f,e}]
              =\sum_{\bnu\in\Kleshf}d^{\C,f}_{\blam\bnu}\Bigl(
                 \sum_{\bmu\in\Klesh}a^{f,e}_{\bnu\bmu}[D^\bmu_F]\Bigr)
             = \sum_{\bmu\in\Klesh}\Bigl(\,\sum_{\mathclap{\bnu\in\Kleshf}}
                 d^{\C,f}_{\blam\bnu} a^{f,e}_{\bnu\bmu}\Bigr)[D^\bmu_F].
    \end{align*}
    Since $[S^\blam_F]=\sum_\bmu d^{F,e}_{\blam\bmu}[D^\bmu_F]$ the result
    follows.
  \end{proof}

  In the special case when $e=f$, \autoref{C:Emu} implies that
  $\mathsf{A}_{F,e}=\mathsf{A}_{F,e}(1)$, where
  $\mathsf{A}_{F,e}(1)$ is the adjustment matrix obtained by setting
  $q=1$ in the graded adjustment matrix $\mathsf{A}_{F,e}(q)$ of
  \autoref{E:Adjustment}.

  \begin{Remark}
    By \autoref{T:AdjustmemtMatrix},  \autoref{C:eAdjustment} extends to
    give graded adjustment matrices in the special case when $e=f$ or,
    equivalently, $r=0$. If  $r>0$ then $\Rnf(\C)\cong\Hnf(\C)$ by
    \autoref{T:KLRIsomorphism}. By \autoref{D:KLRAlgebra} we have a
    cyclotomic KLR algebra $\Rnf(F)$ defined over~$F$, however,
    \autoref{D:CyclotomicHecke} does not define a cyclotomic Hecke
    algebra $\Hnf(F)$ because a field of characteristic~$p>0$ cannot
    contain a primitive $f$th root of unity when~$p$ divides~$f$.
    Consequently, \autoref{T:KLRIsomorphism} does not apply
    to~$\Rnf(F)$ and, \textit{a priori}, the algebra $\Rnf(F)$ is not
    isomorphic to a cyclotomic Hecke algebra. Even so,
    \autoref{T:GradedCellular} shows that $\Rnf(F)$ is a graded cellular
    algebra that comes equipped with Specht modules~$\Slam_F$, simple
    modules $\Dmu_F$ and almost simple modules $\Emu_F$,
    for~$\blam\in\Parts$ and~$\bmu\in\Kleshf$. Let $K$ be the field used
    in the proof of \autoref{P:fAdjustment}. Then $\Hnf(K)\cong\Rnf(K)$
    and reduction modulo~$p$ induces a decomposition map from
    $\Rep\Rnf(K)$ to $\Rep\Rnf(F)$, so the decomposition map of
    \autoref{P:fAdjustment} factors through $\Rep\Rnf(F)$, however, it
    is not clear how to define a map between the Grothendieck groups
    of~$\Rep\Rnf(F)$ and~$\Rep\Rn(F)$.

    Using formal characters, it is straightforward to compute the following graded
    decomposition matrices corresponding to the choices $e=p=2$,
    $\Lambda=\Lambda_0$, $f=4$ and $\Lambda_f=\Lambda_0$ (so
    $F=\mathbb{F}_2$ is the field of order~$2$):
    \[
    \begin{array}{c@{\hspace*{20mm}}c}
          \begin{array}{r|*{3}c}
                  &1^5&2,1^3&2^2,1\\\hline
            1^5   &    1&     &     \\
            2,1^3 &    .&    1&     \\
            2^2,1 &    1&    .&    1\\
            3,1^2 &  2q&    .&    q\\
            3,2   &  q^2&    .&  q^2\\
            4,1   &    .&    q&    .\\
            5     &  q^2&    .&    .\\
          \end{array}
        &
          \begin{array}{r|*{6}c}
                  &1^5&2,1^3&2^2,1&3,1^2&  3,2&  4,1\\\hline
            1^5   &    1&     &     &     &     &     \\
            2,1^3 &    .&    1&     &     &     &     \\
            2^2,1 &    q&    .&    1&     &     &     \\
            3,1^2 &    .&    .&    .&    1&     &     \\
            3,2   &    .&    .&    q&    .&    1&     \\
            4,1   &    .&    .&    .&    .&    .&    1\\
            5     &    .&    .&    .&    .&    q&    .\\
          \end{array}\\[2mm]
          \text{$\Rn[5]=\mathscr{R}^{\Lambda_0}_5\cong \mathbb{F}_2\Sym_5$ with $(e,p)=(2,2)$}
          &
          \text{$\Rnf[5]=\mathscr{R}^{\Lambda_0}_5$ with $(e,p)=(4,2)$}
        \end{array}
    \]
    Comparing these two matrices shows that there is no graded
    adjustment matrix with entries in $\N[q,q^{-1}]$ in this case. (In
    contrast, in agreement with \autoref{P:fAdjustment}, there is an
    adjustment matrix when we set~$q=1$, which corresponds to forgetting
    the grading.) Hence, when $r>0$ we cannot expect the naive
    generalisations of \autoref{P:fAdjustment} and
    \autoref{C:eAdjustment} to hold in the graded setting.  It would be
    interesting to find a connection between the algebras~$\Rnf(F)$
    and~$\Rn(F)$.
  \end{Remark}

  We are now able to prove what is really the main result of this paper.

  \begin{Theorem}\label{T:FullyPositiveJantzen}
    Suppose that $F$ is a field of characteristic $p>0$ and let
    $\blam\in\Parts$. Then
    \[
        \ch J^\blam_F = \sum_{\bmu\in\Klesh} (d^{\C,e}_{\blam\bmu})'(1)
             \ch E^\bmu_F
        + \sum_{r>0}\,(p-1)p^{r-1}\sum_{\substack{f=ep^r\\\bnu\in\Kleshf}}
          (d^{\C,f}_{\blam\bnu})'(1)\ch E^\bnu_{f,e}.
    \]
    Moreover, all of the coefficients on the right hand side are
    non-negative integers.
  \end{Theorem}

  \begin{proof}
    Arguing as in the proof of \autoref{T:PositiveJantzen}, if $f=ep^r$
    then
    \begin{align*}
          \sum_{\bi\in I^n_f}\deg_f(\blam)\,\bi
            &= \sum_{\bnu\in\Kleshf}(d^{\C,f}_{\blam\bnu})'(1)\ch E^\bnu_\C.
        \intertext{Applying the map $\Lf$, and using \autoref{P:fAdjustment}
        and \autoref{D:Efbnu}, this becomes}
          \sum_{\bi\in I^n_f}\deg_f(\blam)\,\bi_f
          &= \sum_{\bnu\in\Kleshf}(d^{\C,f}_{\blam\bnu})'(1)\ch E^\bnu_{f,e}.
       \intertext{Multiplying by $(p-1)p^{r-1}$ if $r>0$ and then summing over~$r$, as in
       \autoref{D:epDegree},}
        \sum_{\bi\in I^n}\pdeg_{e,p,\bi}(\blam)\,\bi
          &= \sum_{\bmu\in\Klesh}d'_{\blam\bmu}(1)\ch E^\bmu_{F}
          + \sum_{r>0}(p-1)p^{r-1}\sum_{\mathclap{\bnu\in\Kleshf}}(d^{\C,f}_{\blam\bnu})'(1)
                \ch E^\bnu_{f,e}.
    \end{align*}
    Hence, the character formula for $\ch J^\blam_F$ follows by
    \autoref{P:JantzenDegCharacter}. Finally, the coefficients on the
    right hand side of the character formula for $\ch J^\blam_F$ are
    non-negative integers by \autoref{L:positivity}.
  \end{proof}

  Applying \autoref{P:Characters}, just as we did after
  \autoref{T:PositiveJantzen}, we obtain the
  \hyperlink{T:MainTheorem}{Main Theorem} from the introduction:

  \begin{Corollary}\label{C:FullyPositiveJantzen}
    Suppose that $F$ is a field of characteristic $p>0$ and let $\blam\in\Parts$. Then
    \[
      \sum_{k>0}[J_k(S^\blam_F)] = \sum_{\bmu\in\Klesh}
          (d^{\C,e}_{\blam\bmu})'(1)[E^\bmu_F]
      + \sum_{r>0}(p-1)p^{r-1}\sum_{\substack{f=ep^r\\\bnu\in\Kleshf}}
                    (d^{\C,f}_{\blam\bnu})'(1)[E^\bnu_{f,e}].
    \]
  \end{Corollary}

  One of the applications of the classical Jantzen sum formula is to
  determined when $S^\blam_F=D^\blam_F$, for $\blam\in\Parts$. If $F$ is
  a field of characteristic zero then \autoref{T:PositiveJantzen}
  implies that $S^\blam_F$ is irreducible if and only if
  $d^{\C,e}_{\blam\bmu}=0$ for $\bmu\ne\blam$, which is true but not
  very insightful. If $p>0$ we obtain new information.

  \begin{Corollary}
    Suppose that $F$ is a field of characteristic $p>0$ and let $\blam\in\Parts$. Then
    $S^\blam_F=D^\blam_F$ if and only $S^\blam_{\C,f}=D^\blam_{\C,f}$,
    for all $f=ep^r$ and $r\ge0$.
  \end{Corollary}

  Although this result does not appear to be the literature it can be
  deduced from \cite[Theorem~4.7(iii)]{JamesMathas:cycJantzenSum}.

  In a similar vein, de Boeck et al~\cite{deBoeckEvseevLyleSpeyer} have
  shown that in positive characteristic the dimensions of the simple
  modules indexed by two column partitions for the Hecke algebra of the
  symmetric with $\xi$ of quantum characteristic~$e$ depend upon the
  dimensions of the simple modules  at $ep^r$th roots of unity in
  characteristic zero. It would be interesting to find stronger
  connections between the representation theories of the Hecke algebras
  with quantum characteristic~$e$ in characteristic~$p>0$ and the Hecke
  algebras in characteristic zero at $ep^r$th roots of unit, for
  $r\ge0$.

  A second application of the Jantzen sum formula is to give upper bound
  for decomposition numbers. For $\blam\in\Parts$ and $\bmu\in\Klesh$
  define
  \[
    j^{F,e}_{\blam\bmu} = (d^{\C,e}_{\blam\bmu})'(1)a_{\blam\bmu}
    + \sum_{r>0}(p-1)p^{r-1}\sum_{\substack{f=ep^r\\\bnu\in\Kleshf}}
      (d^{\C,f}_{\blam\bnu})'(1)a^{f,e}_{\bnu\bmu}.
  \]
  Then comparing the coefficient of $\ch D^\bmu_F$ on both sides of
  \autoref{T:FullyPositiveJantzen} shows that the following holds.

  \begin{Corollary}
    Suppose that $\blam\in\Parts$ and $\bmu\in\Klesh$, with
    $\bmu\ne\blam$. Then $d^{F,e}_{\blam\bmu}\le j^{F,e}_{\blam\bmu}$.
  \end{Corollary}

  In particular, $d^{F,e}_{\blam\bmu}\ne0$ only if $j^{F,e}_{\blam\bmu}\ne0$.

  We find it striking that the upper bound on the decomposition numbers
  of $\Hn(F)$ increases dramatically as soon as
  $d^{\C,f}_{\blam\bnu}(q)$ is non-zero, for some $f=ep^r$. In the
  special case when $\Lambda=\Lambda_0$, so that $\Hn(F)$ is an
  Iwahori-Hecke algebra of the symmetric group,
  \autoref{T:FullyPositiveJantzen} provides some evidence in favour of
  (the now disproved) James' conjecture~\cite[\S4]{James:DecII}. Both
  the statement and proof use the standard terminology of the weight,
  hooks and cores of partitions; see, for example,
  \cite[Chapter~4]{Mathas:ULect}.

  \begin{Corollary}\label{C:James}
    Suppose that $\Lambda=\Lambda_0$ and that $\lambda\in\Parts$ is
    partition $e$-weight $w<p$. Then
   \[
      \ch J^\lambda_F=\ch J^\lambda_\C
           = \sum_{\mu\in\Klesh} (d^{\C,e}_{\lambda\mu})'(1)\ch
           E^\mu_F.
    \]
  \end{Corollary}

  \begin{proof}
    Let $f=ep^r$, where $r>0$. Since $w<p$, the partition $\lambda$
    cannot contain a removable $f$-hook. Therefore, $\lambda$ is an
    $f$-core and the Specht module $S^\lambda_{\C,f}=E^\lambda_{\C,f}=D^\lambda_{\C,f}$
    is irreducible.  Hence, $d^{\C,f}_{\lambda\bmu}(q)=\delta_{\lambda\mu}$
    and \autoref{T:FullyPositiveJantzen} implies the result.
  \end{proof}

  In general we would not expect the character $\ch J^\lambda_F$ to
  determine the Jantzen filtration uniquely.  Combining
  \autoref{C:James} with Williamson's counterexample to the James
  conjecture~\cite{Williamson:JamesLusztig}, shows that $\ch
  J^\lambda_F$ does not determine the Jantzen filtration uniquely even
  when $w<p$.

  \section{The classical sum formula}\label{S:classical}
  For completeness, this section shows how to use formal characters to
  prove a more ``classical'' version of the Jantzen sum
  formula~\cite{JamesMathas:JantzenSchaper,JamesMathas:cycJantzenSum}.
  In the cyclotomic case, the sum formula that we obtain is equivalent
  to, but slightly nicer than, that given by
  \cite[Theorem~4.3]{JamesMathas:cycJantzenSum}.  Mostly this is because
  \autoref{D:CyclotomicHecke} is implicitly invoking the Morita
  equivalence of \cite{DipperMathas:Morita} to restrict the ``cyclotomic
  parameters'' of $\Hn(F)$ to a single $\xi$-orbit. Two other reasons
  why the results in this section are more elegant than the
  corresponding results in \cite{JamesMathas:cycJantzenSum} is because
  we have already fixed the quantum characteristic~$e$ and because the
  combinatorial framework that we introduce directly links beta numbers
  and contents for $\ell$-partitions.

  As in \autoref{S:SeminormalForms}, throughout this section we fix a
  field~$F$ and a primitive $e$th root of unity $\xi\in F$, where
  $e<\infty$, and study the Hecke algebra $\Hn(F)$. We continue assume
  that the charge $\charge$ satisfies \autoref{E:Separation} and we
  consider the algebras $\HO$ and $\HK$ with Hecke parameter $z=x+\xi$,
  where~$\Ocal=F[x]_{(x)}$ and~$\Kcal=F(x)$.

  \subsection{Charged beta numbers}

  We want an analogue of the branching rules for the Gram
  determinants~$G^\blam_\bi$. This requires a mild dose of new notation
  that is motivated by \cite[\S3]{JamesMathas:cycJantzenSum}.

  The \textbf{length} of a partition
  $\lambda=(\lambda_1\ge\lambda_2\ge\dots\ge0)$ is the smallest non-negative
  integer~$L(\lambda)$ such that $\lambda_{L(\lambda)}=0$. The
  \textbf{length} of an $\ell$-partition
  $\blam\in\Parts[m]$ is $L(\blam)=\max\set{L(\lambda^{(l)})|1\le l\le\ell}$.
  If $\blam\in\Parts$ then $L(\blam)\le n$ so, in what follows, we will
  always work with sequences of length~$n$.

  Following Littlewood~\cite{betanumbers}, a sequence of \textbf{beta
  numbers} of length $n$ is a strictly decreasing sequence
  $\beta_1>\dots>\beta_n\ge0$ of $n$ non-negative integers. It is
  straightforward to show that the set of beta numbers of length $n$ are
  in bijection with the partitions of length at most~$n$, where the
  partition $\lambda$ corresponds to the sequence of \textbf{beta
  numbers} $\beta^\lambda_1>\dots>\beta^\lambda_n\ge0$, where
  $\beta^\lambda_r=\lambda_r+n-r$, for $r\ge1$. (If $\lambda$ is a
  partition of length~$L=L(\lambda)$ then the integers
  $(\beta^\lambda_1+L-n>\dots>\beta^\lambda_l+L-n)$ are the first column
  hook lengths in~$\lambda$.) Notice that if $(\beta_1,\dots,\beta_n)$
  is a sequence of beta numbers of length~$n$ then
  $(\beta_1+1,\dots,\beta_n+1,0)$ is a sequence of beta numbers of
  length $n+1$.

  Recall that we are assuming that the charge $\charge$ satisfies
  \autoref{E:Separation}.
  A sequence of \textbf{charged beta numbers} (of length $\ell n$) is a
  sequence $\bbeta=(\beta_1,\dots,\beta_{\ell n})\in\N^{\ell n}$ such
  that
  \begin{equation}\label{E:ChargedBetaNumbers}
    \kappa_l+n>\beta_{n(l-1)+1}>\dots>\beta_{nl}\ge\kappa_{l}-n,
    \qquad\text{for }1\le l\le\ell.
  \end{equation}
  That is, for $1\le l\le n$ the subsequences
  $(\beta_{n(l-1)+1},\dots,\beta_{nl})$ are beta numbers for the
  partitions $\lambda^{(l)}$ \textit{shifted} by $\kappa_l-n\ge0$.
  Assumption~\autoref{E:Separation} ensures that
  $\beta_{1}>\dots>\beta_{\ell n}\ge0$ and, in particular, that
  the charged beta numbers are distinct. When the order is not
  important, call $\set{\beta^\blam_1,\dots,\beta^\blam_{\ell n}}$
  a set of (charged) beta numbers for~$\blam$.

  We frequently need to translate between $\ell$-partitions
  and charged beta numbers. To facilitate this, given an integer
  $s$ such that $1\le s\le\ell n$ write $s=(l_s-1)n+r_s$, where
  $1\le l_s\le\ell$ and $1\le r_s\le n$. Similarly, if $1\le t\le\ell n$
  then we write $t=(l_t-1)n+r_t$ etc. In this way, we refer to
  \textbf{row}~$s$ of~$\blam$ as row~$r_s$ of $\lambda^{(l_s)}$. Of
  course, row~$s$ of~$\blam$ is empty if $r_s>L(\lambda^{(l_s)})$ or,
  equivalently, $\lambda^{(l_s)}_{r_s}=0$.

  As with the beta numbers of partitions, if $\blam\in\Parts[m]$ and
  $n\ge L(\blam)$ then $\blam$ can be associated with a unique sequence
  of charged beta numbers of a given length. Explicitly,
  $\blam\in\Parts[m]$ has sequence of charged beta numbers
  $\bbeta^\blam=(\beta^\blam_1,\dots,\beta^\blam_{\ell n})$, where
  \begin{equation}\label{E:betablam}
     \beta^\blam_s = \kappa_{l_s}+\lambda^{(l_s)}_{r_s}-r_s, \qquad
        \text{where $s=(l_s-1)n+r_s$ as above}.
  \end{equation}
  Even though our notation does not reflect this,
  $\bbeta^\blam=\bbeta^\blam(\charge,n)$ depends
  on $\blam$ and the choice of $n$ and $\charge$.
  Via \autoref{E:betablam}, the beta number~$\beta^\blam_s$ is naturally
  associated with the node $(l_s,r_s,\lambda^{(l_s)}_{r_s})$ since, by
  definition, the content of this node is
  $\con(l_s,r_s,\lambda^{(l_s)}_{r_s})=\beta^\blam_s$. Notice that this
  node belongs to $\blam$ only if $L(\lambda^{(l_s)})\ge r_s$.

  If~$\blam\in\Parts[m]$ then
  $\sum_{r=1}^{\ell n}\beta^\blam_r
     =n(\kappa_1+\dots+\kappa_\ell)-\frac12\ell n(n+1)+m$. Even though
  we do not need this (although compare \autoref{D:betaPerm}), it is not
  difficult to see that there is a bijection between the
  $\ell$-partitions in $\Parts[m]$ of length at most $n$ and the
  sequences of charged beta numbers with sum
  $n(\kappa_1+\dots+\kappa_\ell)-\frac12\ell n(n+1)+m$.

  \begin{Example}\label{Ex:BetaNumbers}
    Let $e=3$ and $\Lambda=\Lambda_0+\Lambda_2$, with $\res(1,1,1)=0$.
    Let $\blam=(3,2|1)$ so that $m=|\blam|=6$. The following table lists
    some of the infinitely many choices of sequence of charged beta
    numbers for $\blam$, corresponding to different choices of~$\charge$
    and $n\ge L(\blam)=2$.
    \[
    \begin{array}{lll} \toprule n & \charge & \bbeta^\blam\\\midrule
             2 & (6,2)  & (8,6,2,0) \\
             3 & (12,5) & (14,12,9,5,3,2) \\
             4 & (15,5) & (17,15,12,11,5,3,2,1) \\
             5 & (15,5) & (17,15,12,11,10,5,3,2,1,0) \\
             6 & (21,8) & \betanumbers \\
             \bottomrule
        \end{array}
    \]
    Notice that $\charge$ must be chosen in accordance with
    \autoref{E:Separation}, so it depends on $e$, $n$ and $\Lambda$, and
    that $\beta^\blam$ depends on $\blam$ and $\charge$. Below we will
    assume that $n\ge m$, as in the last two example beta sequences,
    because this ensures that all partitions in $\Parts[m]$ have a
    charged beta sequence of length~$\ell n$. Consequently, these
    sequences of beta numbers can describe all possible
    ways of moving hooks (see \autoref{SS:hooks}), between the
    components of~$\blam$.
  \end{Example}

  The main results in this chapter are proved by induction on the
  \textit{size} $m$ of the $\ell$-partitions using charged beta
  sequences of \textit{length}~$\ell n$. Accordingly, we work with an
  $\ell$-partition $\blam\in\Parts[m]$, for $n\ge m$, and use the
  sequence of charged beta numbers $\bbeta^\blam$ defined in
  \autoref{E:betablam}, which has length~$\ell n$.

  If $\bi=(i_1,\dots,i_n)\in I^m$ and $1\le k\le m$ let
  $\bi_{\downarrow k}=(i_1,\dots,i_k)\in I^k$.  If $\bbeta\in\N^{\ell n}$
  and $w\in\Sym_{\ell n}$ then let $w\cdot\bbeta=(\beta_{w(1)},
  \dots,\beta_{w(\ell n)})$.  Since sequence of charged beta numbers are
  monotonically decreasing, $w\cdot\bbeta$ is a sequence of charged beta
  numbers if and only if $w=1$.  Abusing notation slightly, let $L(w)$
  be the minimal length of $w\in\Sym_{\ell n}$ as a product of the
  standard Coxeter generators of $\Sym_{\ell n}$.

  \begin{Definition}\label{D:betaPerm}
    Suppose that $\bi\in I^m$ and
    $\bbeta=(\beta_1,\dots,\beta_{\ell n})\in\Z^{\ell n}$, where $n\ge m$.
    Define
    \[
    d_\bi(\bbeta) = \begin{cases*}
      (-1)^{L(w)}\dim S^\blam_\bi,
           & if $\bbeta=w\cdot\bbeta^\blam$ for some $\blam\in\Parts[m]$
               and some $w\in\Sym_{\ell n}$,\\
      0,   & otherwise,
    \end{cases*}
    \]
    where $S^\blam_\bi = S^\blam_F f_\bi$. Similarly, set
    $d_\bi(\blam)=\dim S^\blam_\bi$.
  \end{Definition}

  In particular, $d_\bi(\beta)=0$ if $\bbeta=w\cdot\bbeta^\blam$ and
  $S^\blam_\bi=0$ for some $w\in\Sym_{\ell n}$, or if $\beta_s=\beta_t$
  for $s\ne t$, or if~$\#\set{1\le r\le\ell
  n|\kappa_l+n>\beta_r\ge\kappa_l-n}\ne n$, for any $l$ with $1\le l\le\ell$.

  \subsection{Beta numbers and branching rules}
  The classical branching rules for the Specht modules of the symmetric
  groups are one of the cornerstones of modern representation theory.
  Analogues of these branching rules exist for the Specht modules
  of~$\Hn$~\cite{HuMathas:GradedInduction,Mathas:SpechtRestriction}.
  For this paper we only need the following simple consequence of these
  branching rules.

  Let $\blam\in\Parts$ and $\bnu\in\Parts[n-1]$. Write $\bnu\to\blam$ if
  $\blam=\bnu\cup\set{\alpha}$ for some addable node $\alpha$ of~$\bnu$.
  If $\alpha$ is an $i$-node, for $i\in I$ then we write
  $\bnu\xrightarrow{i}\bmu$ to emphasise the residue. If $\bj\in
  I^{n-1}$ and $i\in I$ let $\bj\vee i=(j_1,\dots,j_{n-1},i)\in I^n$.
  Extend $\_\vee i$ to a linear map $\Z[I^{n-1}]\hookrightarrow\Z[I^n]$.

    \begin{Lemma}\label{L:SpechtBranching}
      Suppose that $\blam\in\Parts$. Then
      $
      \ch S_F^\blam = \Sum_{i\in I}\sum_{\bnu\xrightarrow{i}\blam}\ch S_F^\bnu\vee i.
      $
    \end{Lemma}

    \begin{proof}
      By the main theorem of \cite{Mathas:SpechtRestriction}, the Specht
      module $S_F^\blam$ has a filtration as an $\Hn[n-1](F)$-module where
      the quotients are exactly the Specht modules $S_F^\bnu$, where
      $\bnu\to\blam$, with each Specht module appearing exactly once. As the
      character map $\ch\map{\Rep\Hn[n-1]}\Rep\Ln[n-1]$ is exact the
      result follows.
    \end{proof}

    The advantage of using charged beta numbers, instead of
    $\ell$-partitions, is that it is much easier to describe induction
    and restriction. The next corollary uses charged beta numbers to
    give an equivalent, but easier to apply, analogue of
    \autoref{L:SpechtBranching}:

    \begin{Corollary}\label{C:Branching}
      Suppose that $\bi\in I^m$ and
      $\bbeta=(\beta_1,\dots,\beta_{\ell n})\in\N^{\ell n}$, where $n\ge m$. Then
      \[
          d_\bi(\bbeta)
             = \sum_{\substack{1\le s\le \ell n\\\beta_s\equiv i_m\pmod e}}
      d_{\bi'}(\beta_1,\dots,\beta_{s-1},\beta_s-1,\beta_{s+1},\dots,\beta_{\ell n}),
      \]
      where $\bi'=\bi_{\downarrow(m-1)}$.
    \end{Corollary}

    \begin{proof}
      First suppose that $\bbeta=\bbeta^\blam$, for some
      $\ell$-partition $\blam\in\Parts[m]$.  It is well-known and easy
      to see that removing a node from a $\blam$ corresponds to
      decreasing one of the charged beta numbers of~$\blam$ by~$1$. Now,
      $\bbeta'=(\beta_1,\dots,\beta_{s-1},\beta_s-1,\beta_{s+1},\dots,\beta_{\ell n})$
      is a sequence of charged beta numbers if and only if either
      $\beta_s>\kappa_l$ and $n\mid s$, or
      $\beta_s>\beta_{r+1}+1$ and $n\nmid s$. As above, write $s=(l_s-1)n+r_s$, where
      $1\le l_s\le\ell$ and $1\le r_s\le n$, then this corresponds to
      removing the node $(l_s,r_s,\lambda^{(l_s)}_{r_s})$ from the
      $\ell$-partition corresponding to~$\bbeta$. When
      $\bbeta=\bbeta^\blam$ the result now follows in view of
      \autoref{L:SpechtBranching}. Finally, if
      $\bbeta=w\cdot\bbeta^\blam$, for $w\in\Sym_{\ell n}$, then it is
      easy to see that $\bbeta' $ is a sequence charged beta numbers if and
      only if $\bbeta'=w\cdot\bbeta^\bnu$, for some $\ell$-partition
      $\bnu\in\Parts[m-1]$, which implies the result.
    \end{proof}

  \subsection{Beta numbers and hooks}\label{SS:hooks}
  This section shows that adding and subtracting numbers in a sequence
  of charged beta numbers corresponds to adding and removing rim hooks
  from the corresponding $\ell$-partitions. In order to describe this
  recall that if $\lambda$ is a partition then
  $\lambda'=(\lambda'_1,\lambda'_2,\dots)$ is the partition that is
  \textbf{conjugate} to~$\lambda$, where
  $\lambda_c'=\#\set{r\ge1|\lambda_r\ge c}$.  Equivalently, $\lambda_c'$ is
  the length of column~$c$ of~$\lambda$.

  If $\alpha=(l,a,b)\in\blam$ then the \textbf{$\alpha$-rim hook} of $\blam$
  is the following set of nodes:
  \begin{align*}
    R^\blam_\alpha &= \set{(l,c,d)\in\blam|a\le c\le\lambda^{(l)\prime}_b,
    b\le d\le\lambda^{(l)}_a, \text{ and }(l,c+1,d+1)\notin\blam}.
  \end{align*}
  By definition, $\blam\setminus R^\blam_\alpha$ is (the diagram of) an
  $\ell$-partition.  The
  \textbf{$\alpha$-hook length}
  of~$R^\blam_\alpha$ is $h_\alpha=|R^\blam_\alpha|$ and
  \[
      l_\alpha = \#\set{(l,c,d)\in R^\blam_\alpha|c>a}
               = \lambda^{(l)\prime}_b-a
  \]
  is the \textbf{$\alpha$-leg length} of $R^\blam_\alpha$.
  Let $r_\alpha=n(l-1)+a$ be the row of $\blam$ that contains~$\alpha$.
  The \textbf{foot} of $R^\blam_\alpha$ is the node
  $f_\alpha=(l,a+l_\alpha,b)=(l,\lambda^{(l)\prime}_b,b)\in\blam$, which
  is the leftmost node in the last row of~$R^\blam_\alpha$.  Finally, an
  \textbf{$h$-hook} is any rim hook of length~$h$.

  \begin{Example}
    Let $\blam=(4,2,1|6,5,3,2|,3,1,1)$, $\alpha=(2,1,3)$ and
    $n=m=28$. Then
    $f_\alpha=(2,3,3)$,
    $h_\alpha=6$,
    $l_\alpha=2$ and
    $r_\alpha=n+1=29$. The rim hook $R^\blam_\alpha$ can be pictured
    as follows:
    \[
        \Bigger[12](
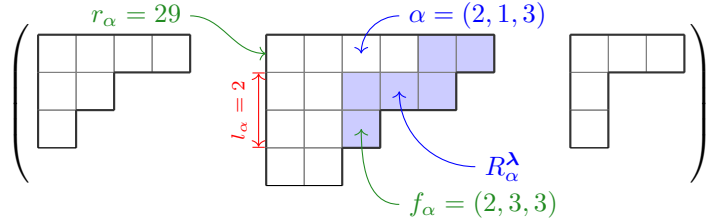
        \begin{tikzpicture}[baseline={(current bounding box.center)}]
            \draw[fill=blue!20, scale=0.5](8,-3) --++(1,0) --++(0,1)
               --++(2,0)--++(0,1)--++(1,0) --++(0,1) --++(-2,0)
               --++(0,-1)--++(-2,0)--cycle;
            \pic (A) at (0,0) {diagram={4,2,1}};
            \pic (B) at (3,0) {diagram={6,5,3,2}};
            \pic (C) at (7,0) {diagram={3,1,1}};
            \draw[<-,blue] ([shift={(-0.25,0.25)}]B-1-3) to[out=90,in=180]++(0.5,0.5)
                 node[right]{$\alpha=(2,1,3)$};
            \draw[<-,ForestGreen] ([shift={(-0.25,0.25)}]B-3-3) to[out=-90,in=180]++(0.5,-1.0)
                 node[right]{$f_\alpha=(2,3,3)$};
            \draw[<-,blue] ([shift={(-0.25,0.25)}]B-2-4) to[out=-90,in=180]++(1,-1)
                 node[right]{$R^\blam_\alpha$};
            \draw[red,thin,|<->|] ([shift={(-0.6,0)}]B-3-1)
                 --node[rotate=90,above=0.8]{\scriptsize$l_\alpha=2$}++(0,1);
            \draw[thin,ForestGreen,<-] ([shift={(-0.5,0.25)}]B-1-1) to[out=180, in=0]
                 ++(-1.0,0.5)node[left,ForestGreen]{$r_\alpha=29$};
        \end{tikzpicture}
        \Bigger[12])
    \]
  \end{Example}
  \begin{Definition}\label{D:Wrapping}
    Let $\blam\in\Parts[m]$ and suppose that $\alpha=(l_s,r_s,c_s)\in\blam$
    and that $t>r_\alpha+l_\alpha=(l-1)n+\lambda^{(l)'}_b$, where $\lambda^{(l)'}$ is the
    partition conjugate to $\lambda^{(l)}$. Define $\blam_{\alpha,t}$ to
    be the $\ell$-partition of~$m$ obtained by wrapping a rim hook of
    length $h_\alpha$ onto $\blam\setminus R^\blam_\alpha$ so that its foot node
    is the leftmost addable node in row~$t$.  If $\blam_{\alpha,t}$ is not
    an $\ell$-partition, set $\blam_{\alpha,t}=\emptyset$. If
    $\blam_{\alpha,t}\ne\emptyset$ then set
    \[\varepsilon_{\alpha,t}=(-1)^{l_\alpha+l_{\alpha'}+1},\]
    where $\alpha'\in\blam_{\alpha,t}$ is the unique node such that
    $\blam_{\alpha,t}\setminus
    R^{\blam_{\alpha,t}}_{\alpha'}=\blam\setminus R^\blam_\alpha$.
    Finally, set $h_\alpha'=\beta^\blam_s-\beta^\blam_t-h_\alpha$, where
    $s=(l_s-1)n+r_s=r_\alpha$.
  \end{Definition}

  By definition, if $\alpha\in\blam$ then the integer $h_\alpha'$
  depends on $\blam$, $\charge$ and $n\ge m$. Notice also that if
  $\blam_{\alpha,t}=\emptyset$ then $\ch S_F^{\blam_{\alpha,t}}=0$.

  \begin{Example}\label{Ex:Wrapping}
    \newcommand\m{\phantom{-}}
     Let $\blam=(3,2|1)\in\Parts[6]$. As in \autoref{Ex:BetaNumbers},
     take $n=m=6$ and $\charge=(20,6)$ so that
     $\bbeta^\blam=\betanumbers$. Using the notation of
     \autoref{D:Wrapping}, the following table gives the complete list
     of pairs $(\alpha,t)$ where $s\ge t$. Here, as in
     \autoref{D:Wrapping}, if $\alpha=(l_s,r_s,c_s)$ then
     $s=(l_s-1)n+r_s=r_\alpha$ and
     $h_{\alpha}'=\beta^\blam_s-\beta^\blam_t-h_\alpha$.
     \begin{NumberedArray}{ccr@{\hspace*{1.3mm}|\hspace*{1.3mm}}lcccclT}
      \alpha & t &\multicolumn{2}{c}{\blam_{\alpha,t}} & s & h_\alpha &h_\alpha'
             & \varepsilon_{\alpha,t}&
       \\\midrule
         (1,1,3) & 3 & (2^2,1&1)  & 1 & 1 & \halpha{1}{3}{1} & (-1)^{0+0+1}=-1 &h\\
         (1,1,3) & 7 & (2^2&2)    & 1 & 1 & \halpha{1}{7}{1} & (-1)^{0+0+0}=\m1& \\
         (1,1,3) & 8 & (2^2&1^2)  & 1 & 1 & \halpha{1}{8}{1} & (-1)^{0+0+0}=\m1& \\
         (1,1,2) & 5 & (1^5&1)    & 1 & 3 & \halpha{1}{5}{3} & (-1)^{1+2+1}=\m1&i\\
         (1,1,2) & 7 & (1^2&4)    & 1 & 3 & \halpha{1}{7}{3} & (-1)^{1+0+0}= -1& \\
         (1,1,2) & 8 & (1^2&2^2)  & 1 & 3 & \halpha{1}{8}{3} & (-1)^{1+1+0}=\m1& \\
         (1,1,2) &10 & (1^2&1^4)  & 1 & 3 & \halpha{1}{10}{3}& (-1)^{1+2+0}= -1& \\
         (1,1,1) & 3 & (2^2,1&1)  & 1 & 4 & \halpha{1}{3}{4} & (-1)^{1+2+1}=\m1&a\\
         (1,1,1) & 5 & (1^5&1)    & 1 & 4 & \halpha{1}{5}{4} & (-1)^{1+3+1}= -1&d\\
         (1,1,1) & 7 & (1&5)      & 1 & 4 & \halpha{1}{7}{4} & (-1)^{1+0+0}= -1& \\
         (1,1,1) & 8 & (1&3,2)    & 1 & 4 & \halpha{1}{8}{4} & (-1)^{1+1+0}=\m1& \\
         (1,1,1) & 9 & (1&2^2,1)  & 1 & 4 & \halpha{1}{9}{4} & (-1)^{1+2+0}= -1& \\
         (1,1,1) &11 & (1&1^5)    & 1 & 4 & \halpha{1}{11}{4}& (-1)^{1+3+0}=\m1& \\
         (1,2,2) & 3 & (3,1^2&1)  & 2 & 1 & \halpha{2}{3}{1} & (-1)^{0+0+1}= -1&q\\
         (1,2,2) & 7 & (3,1&2)    & 2 & 1 & \halpha{2}{7}{1} & (-1)^{0+0+0}=\m1& \\
         (1,2,2) & 8 & (3,1&1^2)  & 2 & 1 & \halpha{2}{8}{1} & (-1)^{0+0+0}=\m1& \\
         (1,2,1) & 3 & (3,1^2&1)  & 2 & 2 & \halpha{2}{3}{2} & (-1)^{0+1+1}=\m1&n\\
         (1,2,1) & 7 & (3&3)      & 2 & 2 & \halpha{2}{7}{2} & (-1)^{0+0+0}=\m1& \\
         (1,2,1) & 9 & (3&1^3)    & 2 & 2 & \halpha{2}{9}{2} & (-1)^{0+1+0}= -1& \\
    \end{NumberedArray}%
    The right hand column identifies when a $2$-partition appears more
    than once in the table as $\blam_{\alpha,t}$, for some~$\alpha$
    and~$t$. \autoref{P:BetaNumbers} below shows that no $2$-partition
    appears more than twice and that if
    $\blam_{\alpha,t}=\blam_{\gamma,u}$, for some
    $\alpha,\gamma\in\blam$ and some integers $t,u$, then
    $\set{h_\alpha,h_\alpha'}=\set{h_\gamma,h_\gamma'}$ and
    $\varepsilon_{\alpha,t}=-\varepsilon_{\gamma,u}$.
  \end{Example}

  Given $\blam\in\Parts$ and two integers $1\le s, t\le \ell n$,
  define the sequence:
  \begin{equation}\label{E:blamst}
    \blamst = (\beta^\blam_1,\dots,\beta^\blam_{s-1},
        \beta^\blam_s-h,\beta^\blam_{s+1},\dots,\beta^\blam_{t-1},
        \beta^\blam_t+h,\beta^\blam_{t+1},\dots,\beta^\blam_{\ell n})\in\Z^{\ell n}.
  \end{equation}
  These sequences will be used heavily in what follows. The next result
  records the combinatorial properties of sequence of charged beta
  numbers that we need.

  \begin{Proposition}\label{P:BetaNumbers}
    Let $\blam\in\Parts[m]$, where $n\ge m$ and fix $\bi\in I^m$ and
    integers $1\le s<t\le\ell n$.  Suppose that $d_\bi(\blamst)\ne0$,
    where $h>0$, and set $h'=\beta^\blam_s-\beta^\blam_t-h$.
    \begin{enumerate}
      \item We have $h\ne h'$ and $d_\bi(\blamst(h'))=-d_\bi(\blamst)$.
      In particular, $d_\bi(\blamst(h'))\ne0$.
      \item If $h\le m$ then $h=h_\alpha$, $h'=h'_\alpha$ and
      $d_\bi(\blamst)=-\varepsilon_{\alpha,t}d_\bi(\blam_{\alpha,t})$,
      for some $\alpha\in\blam$.
      \item If $h>m$ then $l_s<l_t$ and $h'\le m$. Moreover,
      $h'=h_{\alpha'}$, $h=h'_{\alpha'}$ and
      $d_\bi(\blamst)=\varepsilon_{\alpha',t}d_\bi(\blam_{\alpha',t})$,
      for some $\alpha'\in\blam$.
      \item If $\bmu\in\Parts[m]$ and $d_\bi(\bmu)=\pm d_\bi(\blamst)$ then
      $\#\set{(\alpha,t)|\bmu=\blam_{\alpha,t}}=\delta_{l_sl_t}+1$.
    \end{enumerate}
  \end{Proposition}

  \begin{proof}
    By definition, $\beta^\blam_s-h'=\beta^\blam_t+h$ and
    $\beta^\blam_t+h'=\beta^\blam_t-h$. That is,
    $\blamst(h')=(s,t)\cdot\blamst(h)$, which implies that
    $d_\bi(\blamst(h')) =-d_\bi(\blamst)\ne0$.  If $h=h'$ this forces
    $d_\bi(\blamst)=0$, which is a contradiction.  Hence, $h\ne h'$
    completing the proof of part~(a).

    Before considering (b) and (c) recall that going back to Littlewood
    (cf.  \cite[Lemma~5.26]{Mathas:ULect}), it is well-known that
    removing an $h$-rim hook from a partition corresponds to subtracting
    $h$ from one of the beta numbers of the partition. This result
    easily translates into an analogous result for $\ell$-partitions.
    Explicitly, if $\alpha\in\blam$ then $\bmu=\blam\setminus
    R^\blam_\alpha$ if and only if
    \[
      \beta^\bmu_k = \begin{cases*}
         \beta^\blam_{k+1},              & if $r_\alpha<k<r_\alpha+l_\alpha$,\\
         \beta^\blam_{r_\alpha}-h_\alpha,& if $r=r_\alpha+l_\alpha$,\\
         \beta^\blam_k,                  & otherwise.
      \end{cases*}
    \]
    Hence, if $\bj\in I^{m-h_\alpha}$ then
    $d_\bj(\bmu)=(-1)^{l_\alpha}d_\bj(\beta^\blam_1,\dots,
         \beta^\blam_{r_\alpha}-h_\alpha,\dots,\beta^\blam_{\ell n})$.
    In turn, this implies that if~$\bnu\setminus R^{\bnu}_{\alpha'}=\blam$
    then
    $d_\bk(\bnu)=(-1)^{l_{\alpha'}}d_\bk(\beta^\blam_1,\dots,
    \beta^\blam_{r_\alpha}+h_{\alpha'},\dots,\beta^\blam_{\ell n})$,
    for $\bk\in I^{m+h_{\alpha'}}$. Consequently,
    if~$\alpha\in\blam$ and $\bi\in I^m$ then
    $
           d_\bi(\blamst(h_\alpha))
           =(-1)^{l_\alpha+l_{\alpha'}}d_\bi(\blam_{\alpha,t}).
           =-\varepsilon_{\alpha,t}d_\bi(\blam_{\alpha,t}).
    $
    In particular, it follows that  $h\le m$ if and only if~$h=h_\alpha$
    for some $\alpha\in\blam$, so we have proved (b) since
    $h'_\alpha=\beta_s-\beta_t-h_\alpha=h'$.

    Next consider (c) and suppose that $h>m$. Then
    $\set{\beta^\blam_{l_s+1},\dots,\beta^\blam_s-h,\dots,\beta^\blam_{l_s+n-1}}$
    cannot be a set of (shifted) beta numbers for~$\lambda^{(l_s)}$ so
    $l_s<l_t$ and, as in the first paragraph,
    $\set{\beta^\blam_{l_s+1},\dots,\beta^\blam_t+h',\dots,\beta^\blam_{l_s+n-1}}$
    and
    $\set{\beta^\blam_{l_t+1},\dots,\beta^\blam_s-h',\dots,\beta^\blam_{l_t+n-1}}$
    are (shifted) sets of beta numbers for~$\lambda^{(l_s)}$
    and~$\lambda^{(l_t)}$, respectively. Hence, by the last paragraph, $h'=h_{\alpha'}$ for
    some $\alpha'\in\blam$. In particular, $h'\le m$.  Moreover, by parts~(a) and~(b),
    $d_\bi(\blamst)=-d_\bi(\blamst(h'))
                   =\varepsilon_{\alpha',t}d_\bi(\blam_{\alpha',t})$.
    This completes the proof of~(c).

    Finally, as in (d), suppose that $\bmu\in\Parts[m]$ and that
    $d_\bi(\bmu)=\pm d_\bi(\blamst)$.
    Then, there exist $s'<t'$ and~$h>0$ such that
    $\set{\beta^\blam_r|r\ne s,t}=\set{\beta^\bnu_r|r\ne s',t'}$ and
    \[  \set{\beta^\bnu_{s'},\beta^\bnu_{t'}}
             =\set{\beta^\blam_s-h,\beta^\blam_t+h}
             =\set{\beta^\blam_s+h',\beta^\blam_t-h'},
    \]
    where $h'=\beta^\blam_s-\beta^\blam_t-h$. In particular, there are
    exactly two ways to write $\bmu$ in the form $\blam_{\alpha,t}$ if
    $h\le m$ and $h'\le m$ and, as in the last paragraph, there is only
    one to write $\bmu$ in this form if $h>m$ or $h'>m$, in which case
    $l_s<l_t$. Note that~part~(c) implies that if $h>m$ then
    $\bmu=\blam_{\alpha',t}$ and $h'=h_{\alpha'}\le m$, for some
    $\alpha'\in\blam$. Similarly, if $h'>m$ then $\bmu=\blam_{\alpha,t}$
    and $h=h_\alpha\le m$, for $\alpha\in\blam$.
  \end{proof}

  \subsection{Beta numbers and Gram determinants}
    We are now ready to start proving our ``classical'' Jantzen sum
    formula. We start by giving a closed formula for the
    $\gamma$-coefficients from~\autoref{E:gamma}.

    Given two nodes $\alpha=(m,r,c)$ and $\alpha'=(m',r',c')$ define
    $\alpha$ to be \textbf{below} $\alpha'$ if either~$m'<m$ or~$m'=m$ and
    $c'<c$.  Let $\Add_k(\t)$ to be the set of \textbf{addable} nodes
    of the tableau~$\t_{\downarrow k}$ that are below the
    node~$\t^{-1}(k)$.  Similarly, let $\Rem_k(\t)$ be the set of
    \textbf{removable} nodes of the tableau~$\t_{\downarrow k}$ that are
    below the node~$\t^{-1}(k)$.

    \begin{Lemma}\label{L:GammaClosedForm}
      Suppose that $\t\in\Std(\blam)$, for $\blam\in\Parts[m]$. Then
      \[
        \gamma_\t = z^{h_\t}\prod_{k=1}^m
             \dfrac{\Prod_{\alpha\in\Add_k(\t)} [\con_k(\t)-\con(\alpha)]_z}
                   {\Prod_{\beta\in\Rem_k(\t)}[\con_k(\t)-\con(\beta)]_z},
      \]
      for some $h_\t\in\Z$.
    \end{Lemma}

    \begin{proof}
      This follows by induction on the dominance order. The base case is
      given by the closed formula for~$\gamma_{\tlam}$ given
      in~\autoref{E:gammatlam} and the inductive step follows using
      \autoref{C:GammaRecurrence}. Compare with
      \cite[Definition~3.15 and Proposition~3.19]{JamesMathas:cycJantzenSum}.
    \end{proof}

    The next result is a ``branching rule'' for $\det G^\blam_\bi$ using
    the charged beta numbers.

    \begin{Lemma}\label{L:Inductive}
    Suppose that $\blam\in\Parts[m]$ and $\bi\in I^m$. Let
    $\bi'=\bi_{\downarrow(m-1)}$. Then there exists
    an integer $h^\blam_i\in\Z$ such that
    $\det G^\blam_\bi = z^{h^\blam_\bi}A^\blam_1 A^\blam_2
          \Prod_{\bnu\xrightarrow{i_m}\blam} G^\bnu_{\bi'}$, where
    \begin{align*}
         A^\blam_1&=\prod_{h\ge1}
        \prod_{\substack{1\le s<t\le\ell n\\\beta^\blam_s-\beta^\blam_t-1=h\\
                              \beta^\blam_s\equiv i_m}} [h]_z^{
                -d_{\bi'}(\beta^\blam_1,\dots,\beta^\blam_s-h-1,
                     \dots,\beta^\blam_t+h,\dots,\beta^\blam_{\ell n})
         }\\
    \intertext{and}
         A^\blam_2&=\prod_{h\ge1}
             \prod_{\substack{1\le s<t\le\ell n\\\beta^\blam_s-\beta^\blam_t=h\\
                              \beta^\blam_s\equiv i_m}} [h]_z^{
               d_{\bi'}(\beta^\blam_1,\dots,\beta^\blam_s-h,
                     \dots,\beta^\blam_t+h-1,\dots,\beta^\blam_{\ell n})
            }.
    \end{align*}
  \end{Lemma}

  \begin{proof}
    If $m=1$ then there is nothing to prove, so assume $m>1$.
    By \autoref{L:GramiDet} and \autoref{L:GammaClosedForm},
    \begin{align*}
        \det G^\blam_\bi
           &= \prod_{\t\in\Std_\bi(\blam)}\gamma_\t
            = \prod_{\t\in\Std_\bi(\blam)}
              z^{h_\t}\prod_{k=1}^m
                  \dfrac{\Prod_{\alpha\in\Add_k(\t)}[\con_k(\t)-\con(\alpha)]_z}
                        {\Prod_{\beta\in\Rem_k(\t)}[\con_k(\t)-\con(\beta)]_z}\\
          &= z^{a}\prod_{\bnu\xrightarrow{i_m}\blam}
              z^{h^\bnu_{\bi'}}G^\bnu_{\bi'}\cdot \prod_{\t\in\Std_\bi(\blam)}
              \dfrac{\Prod_{\alpha\in\Add_m(\t)}[\con_m(\t)-\con(\alpha)]_z}
                    {\Prod_{\beta\in\Rem_m(\t)}[\con_m(\t)-\con(\beta)]_z}
    \end{align*}
    for some $a\in\Z$, since $\t_{\downarrow(m-1)}\in\Std(\bnu)$ for
    some $\bnu\xrightarrow{i_m}\blam$, for all $\t\in\Std_\bi(\blam)$.
    Set $h^\blam_\bi=a+\sum_\bnu h^\bnu_{\bi'}$.

    Fix $\t\in\Std(\blam)$ and consider its contribution to $\det
    G^\blam_\bi$ in the product above.  As we read down the rows in each
    component, the addable and removable nodes in
    $\Add_m(\t)\cup\Rem_m(\t)$ alternate between addable and removable
    nodes, finishing with the addable node $(k,L(\lambda^{(j)})+1,1)$ at
    the bottom of the $j$th component, for $1\le j\le\ell$. Set
    $\gamma=\t^{-1}(m)=(l_s,r_s,a)$ and let $\beta=(l_v,r_v,c)\in\Rem_m(\t)$
    be a removable node and let $\alpha=(l_u,r_u,b)$ be the addable
    node directly above it. Write $s=(l_s-1)n+r_s$, $u=(l_u-1)n+r_u$ and
    $v=(l_v-1)n+r_v$. Note that $r_v>r_u$ since $\beta$ is below
    $\alpha$, so $v>u$. The consecutive beta numbers between rows $u$ and $v$
    all differ by~$1$, so he contribution from these three nodes to
    $G^\bmu_\bi$ is
    \[
        \frac{[c_m(\t)-c(\alpha)]_z}{[c_m(\t)-c(\beta)]_z}
        = \frac{[\beta^\blam_s-\beta^\blam_{u}-1]_z}
               {[\beta^\blam_s-\beta^\blam_{v}]_z}
        =\prod_{t=u}^{v}
             \frac{[\beta^\blam_s-\beta^\blam_{t}-1]_z}
                  {[\beta^\blam_s-\beta^\blam_{t}]_z}.
    \]
    The addable nodes $(j,L(\lambda^{(j)})+1,1)\in\Add_m(\t)$ at the
    bottom of $\lambda^{(j)}$, for $l\le j\le\ell$, have not yet been
    included in this calculation. Via~$\t$, the contribution of such
    nodes to $\det G_\bi^\blam$ is:
    \[
        \prod_{j=l}^\ell[c_m(\t)-\kappa_j+L(\lambda^{(j)})]
        = \prod_{j=l+1}^\ell [\beta^\blam_s-\beta^\blam_{jn}]_z
        \prod_{t=(j-1)n+L(\lambda^{(j)})+1}^{jn}
           \frac{[\beta^\blam_s-\beta^\blam_t-1]_z}
                {[\beta^\blam_s-\beta^\blam_t]_z}.
    \]
    For each~$j$, all but one of the terms in the right hand product
    cancels out because $\beta^\blam_t=\kappa_j-t$, for
    $(j-i)L(\lambda^{(j)})<t\le jn$. The point of adding the extra
    terms to the last two displayed equations is that this allows us to give
    a generic formula for these determinants.

    If $\s,\t\in\Std_\bi(\blam)$ and
    $\bnu=\Shape(\s_{\downarrow(m-1)})=\Shape(\t_{\downarrow(m-1)})$
    then $\Add_m(\s)=\Add_m(\t)$ and $\Rem_m(\s)=\Rem_m(\t)$, so such
    tableaux contribute the same factor to $\det G^\blam_\bi$. The
    number of these tableaux is
    $\dim S^\bnu_{\bi'}
    =d_{\bi'}(\beta^\blam_1,\dots,\beta^\blam_s-1,\dots,\beta^\blam_{\ell n})$,
    since $\t^{-1}(m)=(l_s,r_s,a)$. Therefore, combining all of contributions
    to $G^\blam_\bi$ from all $\t\in\Std_\bi(\blam)$ shows that
    \[
        \det G^\blam_\bi
        =z^{h^\blam_\bi}\prod_{\bnu\xrightarrow{i_m}\blam}G^\bnu_{\bi'}\cdot
        \prod_{\substack{t=1\\\beta^\blam_t\equiv i_m}}^{\ell n}\left(
        \dfrac{\displaystyle\prod_{t=s+1}^{\ell n}
                    [\beta^\blam_s-\beta^\blam_t-1]_z}
         {\displaystyle\prod_{\substack{t=s+1\\n\nmid t}}^{\ell n}
                  [\beta^\blam_s-\beta^\blam_t]_z}
      \right)^{
        d_{\bi'}(\beta^\blam_1,\dots,\beta^\blam_s-1,\dots,
                   \beta^\blam_t,\dots,\beta^\blam_{\ell n}),
      }
    \]
    In fact, in addition to the terms identified above the right hand
    product above introduces extra terms when $\beta^\blam_t\not\equiv i_m$
    or when the node corresponding to $\beta^\blam_t$ is not removable,
    in which case the exponent
    $d_{\bi'}(\beta^\blam_1,\dots,\beta^\blam_s-1,\dots,\beta^\blam_{\ell n})=0$.
    Hence, all of these ``extra terms'' either cancel out or are equal
    to~$1$, so adding them does not change the determinant.  In
    addition, the condition $n\nmid t$ in the denominator can be omitted
    because all of the extra terms that this introduces to the
    denominator are equal to~$1$ since they all have exponent zero.

    Write the right hand product as $A^\blam_1/A^\blam_2$. Setting
    $h=\beta^\blam_s-\beta^\blam_t-1$, the exponent
    of $[h]_z$ in $A^\blam_1$ is
    \begin{align*}
      d_{\bi'}(\beta^\blam_1,\dots,\beta^\blam_s-1,\dots,
                         \beta^\blam_t,\dots,\beta^\blam_{\ell n})
       &= d_{\bi'}(\beta^\blam_1,\dots,\beta^\blam_t+h,\dots,
                         \beta^\blam_s-h-1,\dots,\beta^\blam_{\ell n})\\
       &= -d_{\bi'}(\beta^\blam_1,\dots,\beta^\blam_s-h-1,\dots,
                         \beta^\blam_t+h,\dots,\beta^\blam_{\ell n}).
    \end{align*}
    Similarly, setting $h=\beta^\blam_s-\beta^\blam_t$ in $A^\blam_2$,
    the exponent  $[h]_z$ in $A^\blam_2$ is
    \[
        -d_{\bi'}(\beta^\blam_1,\dots,\beta^\blam_t+h-1, \dots,
                         \beta^\blam_s-h,\dots,\beta^\blam_{\ell n})
         = d_{\bi'}(\beta^\blam_1,\dots,\beta^\blam_s-h,\dots,
                         \beta^\blam_t+h-1,\dots,\beta^\blam_{\ell n}).
    \]
    Combining the terms in the numerator and denominator in this way
    completes the proof.
  \end{proof}

  Armed with the inductive formula of \autoref{L:Inductive} we can now
  give our first closed formula for the Gram determinant $G^\blam_\bi$.
  This result should be compared with
  \cite[Theorem~3.35]{JamesMathas:cycJantzenSum}. The statement of the
  result, and its proof, require only that $1\le s,t\le \ell n$ and
  $\beta^\blam_s-\beta^\blam_t<h$. In particular, in spite of what the
  notation might suggest, we do not assume that $s$ and $t$ are ordered,
  so both~$s<t$ and~$s>t$ are possible.

  Recall the definition of the sequence $\blamst$ \autoref{E:blamst}.
  Given $1\le s, t\le\ell n$ define $\dkst=(1-\delta_{l_sl_t})n$. That
  is, $\dkst=0$ if $l_s=l_t$ and $\dkst=n$ otherwise. Notice that the
  product formula in the next result does not require $s<t$ although,
  implicitly, it does require $l_s\le l_t$.

  \begin{Proposition}\label{P:BetaGramDet}
     Suppose that $\bi\in I^m$ and $\blam\in\Parts[m]$, and that $n\ge m$.
     Then
     \[
     \det G^\blam_\bi=z^{h^\blam_\bi}\prod_{h\ge1}
       \prod_{\substack{1\le s,t\le\ell n\\
                    \dkst<h\le\beta^\blam_s-\beta^\blam_t}}
              [h]_z^{d_\bi(\blamst)},
        \qquad\text{where }h^\blam_\bi\in\Z.
     \]
  \end{Proposition}

  \newcommand\prodspace{\phantom{\Prod_{\substack{1\le r\le \ell n\\\beta^\blam_r\equiv i_m}}\prod_{h\ge1}\Biggr\{\qquad}}

  \begin{proof}
    We argue by induction on $m$. If $m=0$ then $\blam=(0|\dots|0)$ and
    \[
       \bbeta^\blam=(\kappa_1-1,\kappa_1-2,\dots,\kappa_1-n,
             \kappa_2-1,\dots,\kappa_2-n,\dots,
             \dots,\kappa_\ell-n).
    \]
    By convention empty products are~$1$, so $\det G^\blam_\bi=1$.
    Similarly, since $\dkst=0$ for all $(s,t$) the right hand product is
    equal to $1$ in view of \autoref{P:BetaNumbers}. Hence, the
    proposition holds when $m=0$, giving us the base case of our
    induction.

    Now suppose that $m>0$ and, as in \autoref{L:Inductive}, let
    $\bi'=\bi_{\downarrow(m-1)}$.  Recalling the notation and statement
    of \autoref{L:Inductive} and applying induction,
    \begin{align*}
      z^{-h^\blam_\bi}\det G^\blam_\bi &= A^\blam_1 A^\blam_2
      \prod_{\bnu\xrightarrow{i_m}\blam} \det G^\bnu_{\bi'}\\
         &= A^\blam_1 A^\blam_2 \prod_{\bnu\xrightarrow{i_m}\blam} \prod_{h\ge1}
            \prod_{\substack{1\le s,t\le \ell n\\\dkst<h<\beta^\bnu_s-\beta^\bnu_t}}
             [h]_z^{
             d_{\bi'}(\beta^\bnu_1,\dots,\beta^\bnu_s-h,\dots,\beta^\bnu_t+h,
                  \dots\beta^\bnu_{\ell n})}\\
         &= \prod_{\substack{1\le r\le \ell n\\\beta^\blam_r\equiv i_m}}
         \prod_{h\ge1} \Biggl\{
            \prod_{\substack{1\le s,t\le \ell n, r\ne s,t\\
                     \dkst<h<\beta^\blam_s-\beta^\blam_t}} [h]_z^{
           d_{\bi'}(\beta^\blam_1,\dots\dots,\beta^\blam_s-h,\dots,\beta^\blam_r-1,
                       \dots,\beta^\blam_t+h,\dots\beta^\blam_{\ell n})}\\
         &\prodspace\times A^\blam_1
            \prod_{\substack{1\le t\le \ell n\\
               \dkst<h<\beta^\blam_r-1-\beta^\blam_t}} [h]_z^{
             d_{\bi'}(\beta^\blam_1,\dots,\beta^\blam_r-1-h,\dots,\beta^\blam_t+h,
                \dots\beta^\blam_{\ell n})}\\
         & \prodspace\times A^\blam_2
            \prod_{\substack{1\le s\le\ell n\\
               \dkst<h<\beta^\blam_s-\beta^\blam_r+1}} [h]_z^{
             d_{\bi'}(\beta^\blam_1,\dots,\beta^\blam_s-h,\dots,\beta^\blam_r-1+h,
             \dots\beta^\blam_{\ell n})}
         \quad\Biggr\}\\
         &= \prod_{\substack{1\le r\le \ell n\\\beta^\blam_r\equiv i_m}}
         \prod_{h\ge1} \Biggl\{
            \prod_{\substack{1\le s,t\le \ell n, r\ne s,t\\
                     \dkst<h<\beta^\blam_s-\beta^\blam_t}} [h]_z^{
            d_{\bi'}(\beta^\blam_1,\dots\dots,\beta^\blam_s-h,\dots,\beta^\blam_r-1,
                       \dots,\beta^\blam_t+h,\dots\beta^\blam_{\ell n})}\\
         &\prodspace\prod_{\substack{1\le t\le \ell n\\
               \dkst<h<\beta^\blam_r-1-\beta^\blam_t}} [h]_z^{
             d_{\bi'}(\beta^\blam_1,\dots,\beta^\blam_r-1-h,\dots,\beta^\blam_t+h,
                \dots\beta^\blam_{\ell n})}\\
         & \prodspace\prod_{\substack{1\le s\le\ell n\\
               \dkst<h<\beta^\blam_s-\beta^\blam_r}} [h]_z^{
             d_{\bi'}(\beta^\blam_1,\dots,\beta^\blam_s-h,\dots,\beta^\blam_r-1+h,
             \dots\beta^\blam_{\ell n})}
         \quad\Biggr\},
    \end{align*}
    where the last equality follows using the formulas for $A^\blam_1$
    and $A^\blam_2$ from \autoref{L:Inductive}.  Applying
    \autoref{C:Branching} completes the proof.
  \end{proof}

  \begin{Remark}
    The proof of \autoref{P:BetaGramDet}, and the notation that it uses,
    is different from the results in
    \cite{JamesMathas:JantzenSchaper,JamesMathas:cycJantzenSum} because
    our definition of charged beta numbers corresponds to using first
    \textit{column} hook lengths whereas these references use first
    \textit{row} hook lengths. In the mathematics, this difference
    manifests itself through subtle sign changes. For example,
    \cite{JamesMathas:JantzenSchaper} requires that
    $h>\beta^\lambda_s-\beta^\lambda_t$ whereas we require that
    $h<\beta^\lambda_s-\beta^\lambda_t$. We use column hook
    lengths because of their direct connection with beta numbers.
  \end{Remark}

  \begin{Example}\label{Ex:BetaDet}
    Let $e=3$, $\Lambda=\Lambda_0+\Lambda_2$ and $\blam=(3,2|1)$. As in
    \autoref{Ex:BetaNumbers} and \autoref{Ex:Wrapping}, take $n=6$ and
    $\charge=(20,6)$ so that $\bbeta^\blam=\betanumbers$.  The following
    table gives the different contributions to $\det G^\blam_\bi$, for
    $\bi\in I^n$, given by \autoref{P:BetaGramDet}. The entries in each
    row are complete;y determined by $\bbeta^\blam$ and the integers
    $s,t$ and~$h$. In particular, $h'=\beta^\blam_s-\beta^\blam_t-h$ and
    the $\ell$-partition~$\bnu$ and the \textit{sign} $\varepsilon=\pm1$
    are determined by
    $d_\bi(\beta^\blam_1,\dots,{\color{red}\beta^\blam_s-h},\dots,
         {\color{blue}\beta^\blam_t+h},\dots\beta^\blam_{\ell n})
         =\varepsilon d_\bi(\bnu)$.
         \begin{NumberedArray}[3]{ll@{\hspace*{1.2mm}}*{11}{c@{,\hspace*{1.5mm}}}
               c@{\hspace*{1.5mm}}l*{5}cr@{\hspace*{1.2mm}|\hspace*{1.2mm}}lT}
         & \multicolumn{14}{l}{\blamst=(\beta^\blam_1,\dots,{\color{red}\beta^\blam_s-h},\dots,
                 {\color{blue}\beta^\blam_t+h},\dots\beta^\blam_{\ell n})}
      & s & t & h & h'& \varepsilon & \multicolumn2c{\bnu} &\\\midrule
         \blamsth{1}{ 3}{ 4} & (2^2,1 & 1)     &h\\
         \blamsth{1}{ 7}{14} & (2^2   & 2)     & \\
         \blamsth{1}{ 8}{16} & (2^2   & 1^2)   & \\
         \blamsth{1}{ 5}{ 4} & (1^5   & 1)     &i\\
         \blamsth{1}{ 7}{12} & (1^2   & 4)     & \\
         \blamsth{1}{ 8}{14} & (1^2   & 2^2)   & \\
         \blamsth{1}{10}{16} & (1^2   & 1^4)   & \\
         \blamsth{1}{ 3}{ 1} & (2^2,1 & 1)     &a\\
         \blamsth{1}{ 5}{ 3} & (1^5   & 1)     &d\\
         \blamsth{1}{ 7}{11} & (1     & 5)     & \\
         \blamsth{1}{ 8}{13} & (1     & 3,2)   & \\
         \blamsth{1}{ 9}{14} & (1     & 2^2,1) & \\
         \blamsth{1}{11}{16} & (1     & 1^5)   & \\
         \blamsth{2}{ 3}{ 2} & (3,1^2 & 1)     &q\\
         \blamsth{2}{ 7}{12} & (3,1   & 2)     & \\
         \blamsth{2}{ 8}{14} & (3,1   & 1^2)   & \\
         \blamsth{2}{ 3}{ 1} & (3,1^2 & 1)     &n\\
         \blamsth{2}{ 7}{11} & (3     & 3)     & \\
         \blamsth{2}{ 9}{14} & (3     & 1^3)   & \\
    \end{NumberedArray}%
    (The only purpose of the colours is to highlight which beta numbers
    change. The row labels in the left hand column are to
    facilitate comparison with \autoref{Ex:Wrapping}.) Hence, by \autoref{P:BetaGramDet},
    \begin{align*}
        G_\bi^\blam &=
          \Bigl(\frac{[1]_z}{[4]_z}\Bigr)^{d_\bi(2^2,1|1)}
          [14]_z^{d_\bi(2^2|2)}
          [16]_z^{d_\bi(2^2|1^2)}
          \Bigl(\frac{[4]_z}{[3]_z}\Bigr)^{d_\bi(1^5|1)}
          [12]_z^{-d_\bi(1^2|4)}
          [14]_z^{d_\bi(1^2|2^2)}
          [16]_z^{-d_\bi(1^2|1^4)}
          [11]_z^{-d_\bi(1|5)}
        \\ &\qquad\times
          [13]_z^{d_\bi(1|3,2)}
          [14]_z^{d_\bi(1|2^2,1)}
          [16]_z^{d_\bi(1|1^5)}
          \Bigl(\frac{[1]_z}{[2]_z}\Bigr)^{d_\bi(3,1^2|1)}
          [12]_z^{d_\bi(3,1|2)}
          [14]_z^{-d_\bi(3,1|1^2)}
          [11]_z^{d_\bi(3|3)}
          [14]_z^{-d_\bi(3|1^3)}\,.
    \end{align*}
    By construction, $\det G^\blam_\bi\in\Ocal$.  By \autoref{L:nuPhi},
    $[h]_z$ is a unit in~$\Ocal$ whenever $h\not\equiv0\pmod e$.
    Therefore,
    \[
      G_\bi^\blam
      =u [3]_z^{-d_\bi(1^5|1)}[12]_z^{d_\bi(3,1|2)-d_\bi(1^2|4)}
      =u [3]_z^{d_\bi(3,1|2)-d_\bi(1^2|4)-d_\bi(1^5|1)}
        [4]_{z^3}^{d_\bi(3,1|2)-d_\bi(1^2|4)},
    \]
    for some unit $u\in\Ocal$.
  \end{Example}

  Recall from \autoref{SS:hooks} that $r_\alpha$ is the row index of the
  foot of the $\alpha$-hook $R^\blam_\alpha$.

  \begin{Corollary}\label{C:ClassicalDet}
    Suppose that $\blam\in\Parts$ and $\bi\in I^n$. Then
    $\det G^\blam_\bi = \Prod_{\alpha\in\blam}\prod_{t>r_\alpha}
            [h_\alpha']_z^{\varepsilon_{\alpha,t}d_\bi(\lambda_{\alpha,t})}$.
  \end{Corollary}

  \begin{proof}
    By \autoref{P:BetaGramDet},
    $G^\blam_\bi=\prod_{h,s,t}[h]_z^{d_\bi(\blamst)}$ where,
    in the product, $\dkst<h\le\beta^\blam_s-\beta^\blam_t$.
    By parts~(b) and~(c) of \autoref{P:BetaNumbers}, if
    $\dkst<h\le\beta^\blam_s-\beta^\blam_t$ then there exists
    a node $\alpha\in\blam$ such that~$h=h_{\alpha}'$, $t>r_\alpha$
    and $d_\bi(\blamst)=\varepsilon_{\alpha,t}d_\bi(\blam_{\alpha,t})$.
    Hence, the result follows.
  \end{proof}

  We can now give the promised ``classical'' higher level Jantzen
  character formula.

  \begin{Theorem}\label{T:ClassicalJSF}
    Suppose that $\blam\in\Parts$. Then
    $\ch J^\blam_F = \Sum_{\alpha\in\blam}\sum_{t>r_\alpha}
            \varepsilon_{\alpha,t}\nu_x([h_\alpha']_z)\ch S_F^{\blam_{\alpha,t}}$.
  \end{Theorem}

  \begin{proof}
    Applying \autoref{P:JantzenDetCharacter} and \autoref{C:ClassicalDet},
    \begin{align*}
       \ch J^\blam_F
         &= \sum_{\bi\in I^n} \nu_x\bigl(\det G^\blam_\bi\bigr)\bi
          = \sum_{\bi\in I^n}\Bigl(\sum_{\alpha\in\blam}\sum_{t>r_\alpha}
              \varepsilon_{\alpha,t}\nu_x([h_\alpha']_z)d_\bi(\blam_{\alpha,t})\Bigr)\bi\\
         &= \sum_{\alpha\in\blam}\sum_{t>r_\alpha}
         \varepsilon_{\alpha,t}\nu_x([h_\alpha']_z)\ch S_F^{\blam_{\alpha,t}},
    \end{align*}
    where the last equality follows by swapping the order of summation since, by definition,
    $\ch S_F^{\blam_\alpha,t} = \sum_{\bi\in I^n}d_\bi(\blam_{\alpha,t})\bi$.
    This completes the proof.
  \end{proof}

  \begin{Remark}
    Using \autoref{L:nuPhi}, if $h'\in\N$ then we can explicitly compute
    $\nu_x([h']_z)=\sum_{f\mid h'}\nu_x(\Phi_f(z))$. In particular,
    $\ch S_F^{\blam_{\alpha,t}}$ appears in $\ch J^\blam_F$ only if~$e$
    divides $h_\alpha'$. Further, in view of \autoref{P:BetaNumbers}(d),
    if~$1\le h_\alpha,h_\alpha'\le n$ then
    $\ch S_F^{\blam_{\alpha,t}}$ appears in $\ch J^\blam_F$ if and only
    if~$\nu_x([h_\alpha]_z)\ne\nu_x([h'_\alpha]_z)$.
  \end{Remark}

\bibliography{papers}

\end{document}